\documentclass[12pt]{article}
\usepackage{amsmath,amssymb}
\usepackage{graphicx}
\def\date{\hfill \number\day.  \number\month. \number\year }
\input cyracc.def
\font\tencyr=wncyr10 scaled1200
\def\cyr{\tencyr\cyracc}
\font\tencyss=wncyss10 scaled1200
\def\cyss{\tencyss\cyracc}
\def\sm{\hskip-1pt$\raise1pt\hbox{$\scriptstyle\setminus$}$\hskip-1pt}

\def\ledot{\le\!\!\!\raise2pt\hbox{$\scriptscriptstyle\bullet$}\;}
\def\ldot{\!<\!\!\!\raise1.5pt\hbox{$\scriptscriptstyle\bullet$}}
\def\and{\hskip2pt{\scriptstyle \land}\hskip2pt}
\def\smor{\hskip2pt{\scriptstyle \lor}\hskip2pt}
\def\smcap{\kern 2pt {\scriptstyle \cap}\kern 2pt}
\def\smcup{\kern 2pt {\scriptstyle \cup}\kern 2pt}
\def\ssm{{\smallsetminus}}
\def\nr#1. {{\noindent{\bf #1. }}}

\def\Qed{\hglue 0pt plus 1filll $\square$}

\def\cA{{\cal A}}
\def\cB{{\cal B}}
\def\cC{{\cal C}}
\def\cD{{\cal D}}
\def\cE{{\cal E}}
\def\cF{{\cal F}}

\def\cH{{\cal H}}

\def\cK{{\cal K}}

\def\cM{{\cal M}}

\def\cO{{\cal O}}
\def\cP{{\cal P}}

\def\cR{{\cal R}}
\def\cS{{\cal S}}
\def\cT{{\cal T}}

\def\1{1\kern-2.5pt {\rm l}}    

\def\Cs#1#2{{\rm Cs}_{\hskip1pt #1}{\hskip.5pt#2}} 

\def\Aut{\mathop{{\rm Aut}}}
 
\def\GL#1#2{{\rm GL}_{#1}{#2}}

\def\SL#1#2{{\rm SL}_{#1}{#2}}  
\def\PSL#1#2{{\rm PSL}_{#1}{#2}} 
\def\O#1#2{{\rm O}_{#1}{#2}}
\def\Opr#1#2{{\rm O}^\prime_{#1}{#2}}  
\def\SO#1#2{{\rm SO}_{#1}{#2}}

\def\Spin#1#2{{\rm Spin}_{#1}{#2}}   
\def\U#1#2{{\rm U}_{#1}{#2}}
\def\SU#1#2{{\rm SU}_{#1}{#2}} 
\def\PU#1#2{{\rm PU}_{#1}{#2}}
\def\PSU#1#2{{\rm PSU}_{#1}{#2}}

\def\Gtwo{{\rm G}_2}

\def\lie{\mathop{\strut\frl}}  

\font\Bbb=msbm10 at12 pt
\font\scriptBbb=msbm7  at12 pt
\font\scriptscriptBbb=msbm5 at12 pt

\newfam\Bbbfam
\textfont\Bbbfam=\Bbb
\scriptfont\Bbbfam=\scriptBbb
\scriptscriptfont\Bbbfam=\scriptscriptBbb

\def\CC{{\fam=\Bbbfam C}}

\def\HH{{\fam=\Bbbfam H}}

\def\KK{{\fam=\Bbbfam K}}

\def\OO{{\fam=\Bbbfam O}}

\def\RR{{\fam=\Bbbfam R}}

\def\Ss{{\fam=\Bbbfam S}}
\def\TT{{\fam=\Bbbfam T}}

\def\ZZ{{\fam=\Bbbfam Z}}
\let\Ss=\Ss

\font\sansserif=cmss10 at 12 pt 
\font\scriptsansserif=cmss10 at 8.4 pt
\font\scriptscriptsansserif=cmss10 at 6 pt
\newfam\ssfam
\textfont\ssfam=\sansserif
\scriptfont\ssfam=\scriptsansserif
\scriptscriptfont\ssfam=\scriptscriptsansserif

\let\usuDelta=\Delta  
\let\usuGamma=\Gamma  
\let\usuLambda=\Lambda  
\let\usuOmega=\Omega  
\let\usuPhi=\Phi  
\let\usuPi=\Pi  
\let\usuPsi=\Psi  
\let\usuSigma=\Sigma  
\let\usuTheta=\Theta  
\let\usuUpsilon=\Upsilon  
\let\usuXi=\Xi  

\def\Alpha{{\fam=\ssfam A}}

\def\Chi{{\fam=\ssfam \usuChi}}
\def\Delta{{\fam=\ssfam \usuDelta}}
\def\Epsilon{{\fam=\ssfam E}}
\def\Eta{{\fam=\ssfam H}}
\def\Gamma{{\fam=\ssfam \usuGamma}}
\def\Chi{{\fam=\ssfam X}}

\def\Kappa{{\fam=\ssfam K}}
\def\Lambda{{\fam=\ssfam \usuLambda}}
\def\Mu{{\fam=\ssfam M}}
\def\Nu{{\fam=\ssfam N}}
\def\Omega{{\fam=\ssfam \usuOmega}}
\def\Phi{{\fam=\ssfam \usuPhi}}
\def\Pi{{\fam=\ssfam \usuPi}}
\def\Psi{{\fam=\ssfam \usuPsi}}

\def\Rho{{\fam=\ssfam P}}
\def\Sigma{{\fam=\ssfam \usuSigma}}
\def\Tau{{\fam=\ssfam T}}
\def\Theta{{\fam=\ssfam \usuTheta}}
\def\Upsilon{{\fam=\ssfam \usuUpsilon}}
\def\Ypsilon{{\fam=\ssfam \usuUpsilon}}
\def\Xi{{\fam=\ssfam \usuXi}}
\def\Zeta{{\fam=\ssfam Z}}

\let\epsilon=\varepsilon
\let\theta=\vartheta
\let\phi=\varphi
\let\rho=\varrho

\font\teneufm=eufm10 at 12 pt
\font\seveneufm=eufm7 at 12 pt
\font\fiveeufm=eufm5 at 12 pt
\newfam\eufmfam
\textfont\eufmfam=\teneufm
\scriptfont\eufmfam=\seveneufm
\scriptscriptfont\eufmfam=\fiveeufm
\def\frak#1{{\fam\eufmfam\relax#1}}

\def\frA{{\frak A}}

\def\frC{{\frak C}}
\def\frD{{\frak D}}
\def\frE{{\frak E}}

\def\frL{{\frak L}}

\def\frO{{\frak O}}

\def\frb{{\frak b}}
\def\frc{{\frak c}}

\def\fre{{\frak e}}

\def\frl{{\frak l}}

\def\frr{{\frak r}}
\def\frs{{\frak s}}

\def\frv{{\frak v}}

\def\frx{{\frak x}}

\begin{document}

\abovedisplayskip=3pt
\belowdisplayskip=3pt

\overfullrule=.5pt
\font\bf=cmbx10 scaled 1240
\font\bbf=cmbx10 scaled 1500
\font\Bf=cmbx10 scaled 1400
\def\bullett{\raise1pt\hbox{$\scriptscriptstyle\bullet$}} 
\def\ssm{\smallsetminus}
\def\rk{\rm rk\,}
\let\hat=\widehat
\let\tilde=\widetilde
\let\bold=\bf
\let\ß=\ss{}
\def\bib{\bibitem{} }
\font\teneufm=eufm10 scaled 1200 
\font\seveneufm=eufm10 
\newfam\eufmfam
\textfont\eufmfam=\teneufm
\scriptfont\eufmfam=\seveneufm

\def\frak#1{{\fam\eufmfam\relax#1}}
\def\frL{{\frak L}}
\def\fre{{\frak e}}

\def\3ast{$(\lower1pt\hbox{$\scriptstyle**$})\hskip-13pt\raise3.2pt\hbox{$\scriptstyle*$}\hskip6pt$}
\def\Rtimes{\times\hskip-3pt\raise0pt\hbox{$\vrule height 6pt width .8pt$}}
\def\smotimes{\begin{scriptsize}$\otimes$\end{scriptsize}}
\def\smoplus{\begin{scriptsize}$\oplus$\end{scriptsize}} 
\def\tsmotimes{\begin{scriptsize}$\tilde\otimes$\end{scriptsize}}
\def\tsmoplus{\begin{scriptsize}$\tilde\oplus$\end{scriptsize}}

\vskip24pt

\centerline {\bbf  Compact planes, mostly 8-dimensional. A  retrospect} 

\par\bigskip
\centerline{by Helmut R. \sc Salzmann}
\par\bigskip
\begin{abstract}
 Results on $8$-dimensional topological planes are scattered in the literature. It is the aim of the present paper to give a survey of  these geometries, in particular of information obtained after the appearance of the treatise 
{\it Compact Projective Planes\/} \cite{cp} or not included in \cite{cp}. For some theorems new proofs are given and a few related results concerning planes of other dimensions are presented.
\end{abstract}
\par\smallskip
As the word {\it plane\/} is associated with the notion of a $2$-dimensional geometry, at  first glance the title might seem to be self-contradictory. It is a fact, however, that the point set $P$ of a projective plane  $\cP{\,=\,}(P,\frL)$ may carry a (locally) compact topology of the {\it topological\/} (covering) dimension  $\dim P{\,=\,}8$ such that the geometric operations of joining and intersecting are continuous with respect to a suitable topology on the set $\frL$ of lines. In this case $\cP$ will be called a (compact) $8$-dimensional (topological)  plane.  The classical example, of course, is the Desarguesian projective plane $\cH{\,=\,}\cP_{\HH}$ over the locally compact skew field $\HH$ of the (real) quaternions (see \cite{eb} Chapt.\,7 and \cite{cp} \S\hskip2pt14, cf. also \cite{cs}).
\par\smallskip
Other models abound; most of them have been constructed by suitably modifying the algebraic operations of the skew field $\HH$. There is a vast variety of possibilities; as in the case of finite planes, only the more homogeneous ones  can be classified or can be treated in a reasonable way. The degree of homogeneity can be expressed by the size of the automorphism group $\Sigma{\,=\,}\Aut\cP$ of all continuous collineations. Taken with the compact-open topology, $\Sigma$ is always a locally compact transformation group of $P$ (as well as of $\frL$) with a countable basis 
(\cite{cp} 44.3). For a locally compact group $\Sigma$, the covering dimension $\dim\Sigma$ coincides with the inductive dimension ${\rm ind\,}\Sigma$ and with the maximal dimension of a euclidean ball contained in $\Sigma$, see \cite{cp} 93.\hskip1pt5,\hskip1pt6. This dimension has turned out to be a measure for the size of $\Sigma$ best adapted to our purposes.
\par\smallskip
By a fundamental theorem due to L\"owen \cite{Lw} or \cite{cp} 54.11, compact projective planes with a point space $P$ of positive finite covering dimension exist only for $\dim P{\,=\,}2^m$ and $1{\,\le\,}m{\,\le\,}4$. For reasons to be mentioned later, the case $\dim P{\,=\,}8$ is the most complicated and difficult one, and we focus on this case. For $16$-dimensional planes see 
\cite{hs}, \cite{sz9}.
\par\smallskip
Homogeneity in the sense that $\Sigma$ acts transitively on $P$ is too strong a condition, it is satisfied only by the $4$ classical planes over a locally compact connected (skew) field or the (real) octonion algebra $\OO$, see  \cite{lw1} or \cite{cp} 63.8. In  fact, each point transitive group contains 
the compact elliptic motion group of the corresponding classical plane (\cite{cp} 63.8), hence it is even transitive on the set of {\it flags\/} (=\;incident point-line pairs). 
\par\bigskip
\newpage
{\Bf 1. Review of basic facts}
\par\medskip
We collect some results necessary for understanding the subsequent parts. The symbol $\cP$ will denote  a compact $8$-dimensional projective plane, if not stated otherwise. Note that line pencils 
$\frL_p{\,=\,}\{L{\,\in\,}\frL \mid p{\,\in\,}L\}$ are homeomorphic to lines.  Hence the dual of $\cP$ is also a compact $8$-dimensional plane.
In fact, the spaces $\cP$ and $\frL$ are homeomorphic (Kramer \cite{kr2}). {\tt Notation} is standard, a quaternion will be written in the  forms 
$c{\,=\,}c_0{\,+\,}c_1i{\,+\,}c_2j{\,+\,}c_3k{\,=\,}c_0{\,+\,}\frc{\,=\,}c'{\,+\,}c''j$ with $c_\nu{\,\in\,}\RR$ and $c',c''{\,\in\,}\CC$. For a locally compact group $\Gamma$ and a closed subgroup $\Delta$ the 
coset space $\{\Delta\gamma\mid\gamma{\,\in\,}\Gamma\}$ will be denoted by $\Gamma/\Delta$, its dimension $\dim\Gamma{-}\dim\Delta$ by $\Gamma{:\hskip2pt}\Delta$.
\par\medskip
{\bf 1.1 Lines} (L\"owen). {\it Each line of $\cP$  is homotopy equivalent to a $4$-sphere $\Ss_4$\/}, see \cite{cp} 54.11. In all known examples, the lines of $\cP$ are even topological manifolds and 
then they are  homeomorphic  $(\approx)$ to $\Ss_4\,$ (\cite{cp} 52.3). 
It seems to be very difficult to decide whether or not this is true in general.
 (In contrast,  lines of a compact plane of dimension $2$ or $4$ are known to be $m$-spheres, 
 cf. \cite{sz1} 2.0 and \cite{cp} 53.15. Each line $L$ of a $16$-dimensional plane is  homotopy 
 equivalent  to an $8$-sphere $\Ss_8$, in all known cases $L{\,\approx\,}\Ss_8$, see again
 \cite{cp} 54.11 and 52.3.\,)
\par\medskip
{\bf 1.2 Baer subplanes.}  {\it Each $4$-dimensional closed subplane $\cB$ of a  compact $8$-dimensional plane  $\cP$ is a Baer subplane\/}, i.e., each point of $\cP$ is incident with a line of $\cB$ (and dually, each line of $\cP$ contains a point of $\cB$), see \cite{sz2} \S~\hskip-3pt3 or 
\cite{cp} 55.5 for details. If $\cP$ contains a closed Baer subplane $\cB$, it follows easily that the pencil of lines through a point outside~$\cB$ is a manifold, and hence
the lines of $\cP$ are homeomorphic to $\Ss_4$, see \cite{cp} 53.10 or \cite{sz2}~3.7. By a result of L\"owen \cite{lw2}, any two closed Baer subplanes of $\cP$ have a point and a line in common. This is remarkable because a finite Pappian plane of order $q^2$ is a  union of  $q^2{-}q{+}1$ disjoint Baer subplanes, see \cite{sz2} 2.5. Generally, 
$\langle \cM\rangle$ will denote the smallest closed subplane of $\cP$ containing the  set $\cM$ of points and lines. We write $\cB{\,\,\ldot}\cP$ if $\cB$ is a Baer subplane.
\par\medskip
{\bold 1.3 Large groups.} It is convenient to restrict attention to a {\it connected\/}  locally compact (and therefore closed)   subgroup $\Delta{\,\le\,}\Sigma{\,=\,}\Aut\cP$, usually to the connected component $\Sigma^1$. The letter  $\Delta$ will always be used in this sense; $\Delta'$ denotes the commutator subgroup. \\  {\it $\Delta$ is a Lie group whenever $\dim\Sigma{\,\ge\,}12$\/}, see \cite{pr}. {\it If $\dim\Sigma{\,>\,}18$, then $\cP$ is isomorphic to the classical plane $\cH$, and $\Sigma{\,\cong\,}\PSL3\HH$ is a simple group of dimension $35$\/}. For $\dim\Sigma{\,\ge\,}23$, a proof is given in \cite{cp} 84.27. All planes with $\dim\Sigma{\,\ge\,}17$ have been described explicitly; in fact, they are either translation planes or Hughes planes, see  \cite{sz3}. In a series of papers, H\"ahl has determined all $8$-dimensional translation planes admitting a $17$-dimensional group, cf. \cite{hl1} or \cite{cp} 82.25; more details are given in 1.10 below.
The  Hughes planes depend on a real parameter, they  can be characterized as  follows: $\cP$ has a $\Sigma$-invariant Baer subplane  $\cC$ such that $\Sigma$ induces on $\cC$~the group $\PSL3\CC$, see \cite{cp} \S\hskip2pt86. If $\cP$ is a proper Hughes plane, then $\dim\Sigma{\,=\,}17$ and  $\Sigma$ is transitive on the set of flags of the {\it outer\/} subgeometry consisting of the points and lines not belonging to $\cC$, see \cite{cp} 86.5. 
Sometimes it is convenient to consider~the~classical plane together with the  stabilizer of a Baer subplane also as a Hughes  plane, cf. \cite{sz2}~3.19. 
\par\smallskip
{\tt Remark.} In the case of $16$-dimensional planes, $\Delta$ is known to be a Lie group, if 
$\dim\Delta{\,\ge\,}27$, see \cite{psz}; the full automorphism group $\Sigma$ is a Lie group, if 
$\dim\Sigma{\,\ge\,}29$ \ (\cite{sz13}).

\begin{picture}(0,0)(0,-0)  
\put(0,0){}                        
\put(300,-120){\includegraphics[width=6cm]{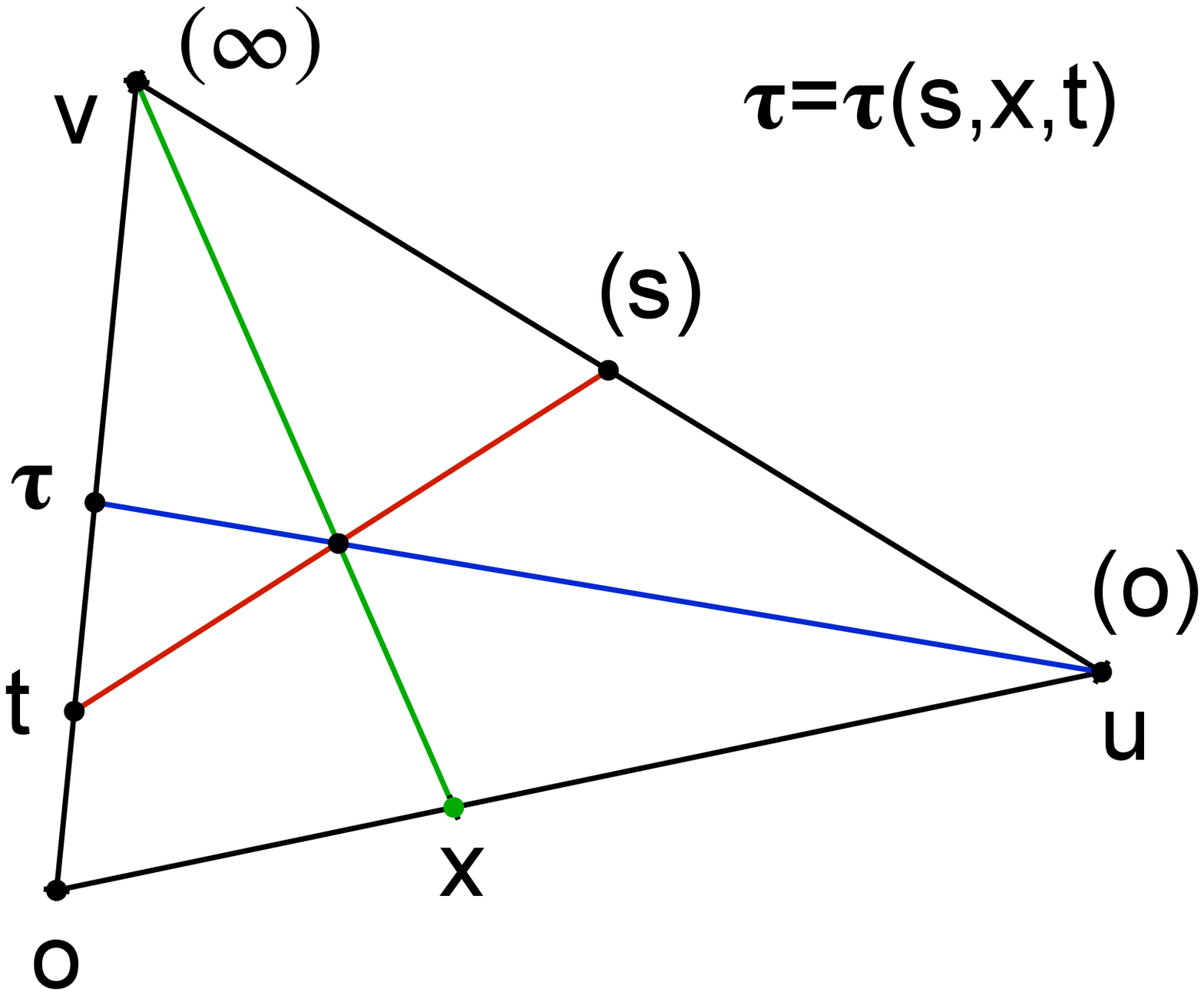}}
\end{picture}

\parbox{300pt}{
{\bold 1.4 Coordinates.} Let ${\fre}{\,=\,}(o,e,u,v)$ be a (non-degenerate) quadrangle in $\cP$. Then the affine subplane  $\cP\sm uv$ can be coordinatized with respect to  the {\it frame\/} $\fre$ by a so-called  {\it ternary field\/} $H_\tau$, where $H$ is homeomorphic to an affine line, the affine point set is written as $H{\times}H$, lines  are given by an equation $y{\,=\,}\tau(s,x,t)$ or $x{\,=\,}c$ (verticals),
they are denoted by $[\hskip1pt s,t\hskip1pt]$ or  $[\hskip1pt c\hskip1pt]$ respectively. 
Parallels are} \\ [6pt] 
 either vertical or both have the same {\it``slope''\/} $\hskip-2pt s$. The axioms of an affine plane can easily~be translated into properties of $\tau$,  see \cite{cp} \S\hskip2pt22. 
Write $e{\,=\,}(1,1)$ and put $x{\,+\,}t{\,=\,}\tau(1,x,t)$. In the Desarguesian case  
$\tau(s,x,t){\,=\,}s\hskip1pt x+t$. A~ternary field $H_\tau$ is called a {\it Cartesian field\/} if  $(H,+)$ 
is a group and if identically $\tau(s,x,t){\,=\,}s{\bullett}x+t$ with  multiplication 
$s{\bullett}x{\,=\,}\tau(s,x,0)$; equivalently, the maps $(x,y)\mapsto(x,y{+}t)$ form a transitive group of translations.
Any ternary field $H_\tau$ such that $H$ is homeomorphic to $\RR^4$ and $\tau{\,:\,} H^3{\,\to\,} H$ is continuous yields a compact projective plane $\cP$, 
 cf. \cite{cp} 43.6.
\par\medskip

{\bold 1.5 Stiffness.} In the classical plane $\cH$, the stabilizer $\Lambda{\,=\,}\Sigma_{\fre}$ of any frame $\fre$ is isomorphic to $\SO3\RR$; in particular,  $\Lambda$ is compact and 
$\dim\Lambda{\,=\,}3$. In any plane, $\Lambda$ can be identified with the automorphism group of the ternary field $H_\tau$ defined with respect to $\fre$. The fixed elements  of $\Lambda$ form a closed subplane 
$\cE{\,=\,}\cF_\Lambda$.  It is not known if $\cE$ is always connected or if $\Lambda$ is  compact in general. Therefore the following 
{\it stiffness\/} results play an important r\^ole:
 (1) $\dim\Lambda{\,\le\,}4$ \ (B\"odi \cite{B}), \par
\qquad\ (2) {\it If $\cE$ is connected or if $\Lambda$ is compact,  then $\dim\Lambda{\,\le\,}3$\/} \ 
(\cite{cp} 83.\,12,\,13),\par
\qquad\ (3) {\it if $\cE$ is contained in a Baer subplane $\cB$,  then $\cE$ is connected and the connected \par \hskip45pt component $\Lambda^1$ of  $\Lambda$ is compact\/} \     (\cite{cp} 55.4 and \cite{cp} 83.9 or   \cite{sz5} (*)\,), \par
\qquad\ (4) {\it if, moreover, $\cB$ is $\Lambda$-invariant, then $\dim\Lambda{\,\le\,}1$\/} \ (\cite{cp} 83.11), hence \par
\qquad\ (5) {\it if $\Lambda$ is compact and $\dim\Lambda{\,>\,}1$, then $\Lambda^1{\,\cong\,}\SO3\RR$, because 
$\Lambda$ does not contain a \par\hskip45pt central  involution\/}, \par
\qquad\ (6) {\it the stabilizer $\Omega$ of a degenerate quadrangle has dimension at most $7$\/} (\cite{cp} 83.17), \par
\qquad\ (7) {\it if $\dim\Omega{\,=\,}7$, then $\Omega^1{\,\cong\,}e^\RR{\cdot\,}\SO4\RR$\/}\ 
{\it and lines are $4$-spheres\/}\ (\cite{sz5} (**)\,), \par
\qquad\ (8) {\it If a group $\Phi{\,\cong\,}\SO3\RR$ fixes a line $W$\hskip-3pt, then each involution 
in $\Phi$ is planar. \par \hskip45pt Either $\Phi$ has no fixed point on $W$ or $\cF_\Phi$ is a\/} flat 
(=\,{\it $2$-dimensional{\rm)} subplane\/} \par\hskip45pt (\cite{sz4} Observation).
\par\medskip
{\bf 1.6 Reflections and translations.} {\it Let $\sigma$ be a reflection with axis $W$ and center~$c$ in the connected group $\Delta$, and let $\Tau$ denote the  group of translations in $\Delta$ with axis~$W$. If  $W^\Delta{\,=\,}W$  and   $\dim c^\Delta{\,=\,}k{\,>\,}0$, then  
$\,\sigma^\Delta\sigma{\,=\,}\Tau{\,\cong\,}\RR^k$, 
$\,\tau^\sigma{\,=\,}\tau^{-1}$ for each $\tau{\,\in\,}\Tau$, and $k$ is even\/}.
This improves \cite{cp} 61.19b; a detailed {\tt proof} can be found in \cite{sz7} Lemma 2. 
\par\medskip
\break
{\bold 1.7 Stabilizer of a triangle.} {\it Assume that $\Delta$ fixes a triangle  and that
$\dim\Delta{\,\ge\,}10$. Then $\Delta$ is a Lie group and the lines of $\cP$ are homeomorphic to 
$\Ss_4$. If $\dim\Delta{\,=\,}11$, then $\Delta$ is equivalent to its classical  counterpart   in 
$\cH$\/}, see \cite{gs} 10.9. {\it In any case, $\dim\Delta{\,\le\,}11$ and $\Delta$ has a compact subgroup isomorphic to $\SO4\RR$ or to $\Spin4\RR$\/}.
\par\smallskip
{\tt Proof.} Denote the fixed points of $\Delta$ by $o,u,v$ and let $S_0{\,=\,}uv\sm\{u,v\}$,
$\,S_1{\,=\,}ou\sm \{o,u\}$, and $\,S_2{\,=\,}ov\sm \{o,v\}$ be the sides of the fixed triangle. 
For each point $z$ on one of these sides, the stabilizer $\Delta_z$ satisfies $\dim\Delta_z{\,\le\,}7$
by 1.5\,(6);  in case of equality, lines are manifolds and $\Delta_z$ is a Lie group with a subgroup  $\SO4\RR$, see 1.5\,(7). Richardson's classification  of compact groups on $\Ss_4$ (to be referred to by $(\dagger)$, for details see \cite{cp} 96.34) applies to the action of a maximal compact subgroup $\Phi$ on the sides $S_\nu$. Note that  all maximal compact subgroups of $\Delta$ are conjugate by the Mal'lcev-Iwasawa theorem  \cite{cp} 93.10. There are two cases: \\
(a) $\dim\Delta_x{\,=\,}\dim\Delta_y{\,=\,}7$ for points $x,y$ on distinct sides, say $x{\,\in\,}S_1$ and 
 $y{\,\in\,}S_2$. Up to conjugacy, $\Phi_x{\,\cong\,}\Phi_y{\,\cong\,}\SO4\RR$. By~$(\dagger)$, 
 $\,\Phi_y$ induces on $ov$ the group  $\Phi_y/\Kappa{\,\cong\,}\SO3\RR$ with 
 $\Kappa{\,\cong\,}\Spin3\RR$, and  $x^{\Kappa}{\,\approx\,}\Ss_3$. Hence 
 $\Phi{\,=\,}\Phi_x\Kappa$,  $\,\dim\Phi{\,=\,}9$, each factor of $\Phi$ is isomorphic to $\Spin3\RR$ and acts trivially on one side $S_\nu$. Therefore  $\Phi$ is equivalent to its classical counterpart. 
The fixed elements of a subgroup $\Lambda{\,\cong\,}\SO3\RR$ of $\Phi$ form a $2$-dimensional subplane 
$\cF_\Lambda$ (see 1.5(8)\,), the radical $\Rho{\,=\,}\sqrt\Delta$ acts effectively on 
$\cF_\Lambda$, and $\Rho$ is commutative by \cite{cp} 33.10. This completes the proof for 
$\dim\Delta{\,=\,}11$.  \\
(b) $\dim\Delta{\,=\,}10$ and there are least two sides $S_\nu$ such that each  orbit of $\Delta$ on  $S_\nu$ is open. Then 
$\Delta$ is transitive, say on $S_1$ and $S_2$. By \cite{cp} 96.14 or \cite{HK}, lines are manifolds 
and $\Delta$ induces  Lie groups on these two sides. Because the kernels of the two actions intersect trivially, 
$\Delta$ itself is a Lie group. We may assume that a maximal compact subgroup of $\Delta_x$ with $x{\,\in\,}S_1$ is contained in $\Phi$. As $S_1$ is homotopy equivalent to $\Ss_3$, there is an exact homotopy sequence
 $$\dots \ZZ_2{\,=\,}\pi_4S_1\to \pi_3\Phi_x\to \pi_3\Phi\to \pi_3 S_1{\,=\,}\ZZ\to 0 \,,\leqno{(*)}$$
 see \cite{cp} 96.12.
 In particular, $\Phi_x{\,<\,}\Phi$ and $\Phi{\,\ne\,}\1$. By $(\dagger)$ the group $\SU3\CC$ does not act on 
 $\Ss_4$. As  the fixed elements of $\SO3\RR$ form a $2$-dimensional subplane, $\Phi$ is a product of factors  
 $\Spin3\RR$.  Either $\dim\Delta_z{\,=\,}7$ for some $z{\,\in\,}S_0$ and $\Phi$ contains $\SO4\RR$ by   
1.5\,(7), or $\Delta$ is also transitive on~$S_0$. If $\dim\Phi{\,=\,}3$, the central involution $\sigma{\,\in\,}\Phi$ is a reflection (or else $\cF_\sigma$ is a Baer subplane and $\Phi|_{\cF_\sigma}{\,\cong\,}\SO3\RR$ by Stiffness, but  
 $\Phi$ is transitive or trivial on each fixed line of $\cF_\sigma)$. Let $S_\nu$ be the axis of $\sigma$. Then 
$(\dagger)$ implies  $\Phi_x{\,=\,}\Phi$ for some $x{\,\in\,}S_\nu$, but this is in conflict with $(\ast)$. Hence $\dim\Phi{\,\ge\,}6$. 
\par\medskip
{\bf 1.8 Fixed configurations.} 
It is the ultimate goal to describe all sufficiently homogeneous planes $\cP$, more precisely, to find
all pairs $(\cP,\Delta)$, where $\dim\Delta{\,\ge\,}b$ for a suitable bound $b$. According to 1.3, 
this programme has been completed for $b{\,=\,}17$. With few exceptions, the bound $b$ will be chosen in the range $12{\,\le\,}b{\,<\,}17$, and then $\Delta$ is always a Lie group \cite{pr}. In particular,  $\Delta$ is semi-simple, or $\Delta$ contains a central torus subgroup or a minimal normal vector subgroup, cf. \cite{cp} 94.26.  The results to be discussed vary with the structure of $\Delta$ and the configuration $\cF_\Delta$ of the fixed elements of 
$\Delta$. In some cases, 
\newpage
only the group  $\Delta$ and its action on the plane can be determined, but the corresponding planes cannot be described explicitly. As shown in 1.7, it follows from 
$\dim\Delta{\,\ge\,}12$ that $\cF_\Delta$ does not contain a triangle. Up to duality, there remain the following possibilities: \par 
\qquad (0) $\cF_\Delta{\,=\,}\emptyset$, \par
\qquad (1) $\cF_\Delta{\,=\,}\{W\}$ is a unique line, \par
\qquad (2) $\cF_\Delta$ consists of exactly $2$ elements, necessarily a point and a line,  \par
\hskip45pt and these may be incident or not, \par
\qquad (3) $\cF_\Delta{\,=\,}\{u,v,uv\}$ consists of exactly $3$ elements, \par
\qquad (4) $\cF_\Delta$ consists of more than $2$ points on a line, \par
\qquad (5)  $\cF_\Delta$ consists of $2$ points and $2$ lines (a {\it double flag\/}).
\par
Least can be said in case (2) or if $\Delta$ has a normal vector subgroup. In some cases even bounds $b{\,\ge\,}9$ may be considered.
\par\medskip
{\bf 1.9 Fixed elements.} The Lefschetz fixed point theorem implies that \par 
\hskip20pt {\it each homeomorphism $\phi:P\to P$ has a fixed point\/}. \\
(a) {\it By duality, each automorphism of $\cP$ fixes a point and a line\/}, 
see \cite{cp} 55.19,\hskip1pt45. \\ (b) {\it The solvable radical $\Rho{\,=\,}\sqrt\Delta$ of $\Delta$ fixes some element of $\cP$\/}. \\ (c) {\it If $\cF_\Delta{\,=\,}\emptyset$, then $\Delta$ is semi-simple with trivial center, or $\Delta$ induces a simple group on \par\hskip20pt some connected  closed $\Delta$-invariant  subplane\/}.
\par\smallskip
{\tt Proof.} The following argument or  slight variations thereof (referred to by$(A)\,$)  will be used repeatedly: if $\Theta$ is a commutative connected normal subgroup of $\Delta$ and if 
$\1{\,\ne\,}\zeta{\,\in\,}\Cs{}\Theta$,  then $p^\zeta{\,=\,}p$ for some point $p$, either 
$p^\Theta{\,=\,}p$, 
or  $p^\Theta$ is contained in a fixed line of  $\Theta$, or  $p^\Theta$ 
generates a connected (closed) subplane $\cS{\,=\,}\langle p^\Theta\rangle$ and 
$\zeta|_\cS{\,=\,}\1$. In the latter case, $\overline\Theta{\,=\,}\Theta|_\cS{\,\ne\,}\1$, 
and $\cS$ is a proper subplane of $\cP$.\\
(b) If $\dim\cS{\,=\,}2$, then $\cS$ has no proper closed subplane (cf. \cite{cp} 32.7), and $\Theta$ 
has a fixed element in $\cS$. If $\cS$ is a Baer subplane, then  $(A)$  can be applied to 
$\overline\Theta$; again  $\cF_\Theta{\,\ne\,}\emptyset$. Induction over the solvable length of 
$\Rho$ and repeated use of (A)  implies that also $\cF_\Rho{\,\ne\,}\emptyset$. \\
(c)  will be proved successively for planes $\cR$ of dimension $2,\,4,$ and $8$. 
If $\Delta$ is not semi-simple, then  $\Rho{\,=\,}\sqrt\Delta{\,\ne\,}\1$ by definition, and $\Rho$ fixes some element by step (b), say $p^\Rho{\,=\,}p$.  Assume also that  $\cF_\Delta{\,=\,}\emptyset$. 
Then $p^\Delta$ is not contained in a line and  $\langle p^\Delta\rangle{\,=\,}\cS{\,\le\,}\cR$ is a 
closed subplane; normality of $\Rho$ implies $\Rho|_\cS{\,=\,}\1$. If $\zeta{\,\ne\,}\1$ is a central element of  $\Delta$, then (A) yields a common fixed element $p$ of $\zeta$ and $\Rho$, 
 and  $\zeta|_\cS{\,=\,}\Rho|_\cS{\,=\,}\1$. \\ 
If $\dim\cR{\,=\,}2$, there is no proper closed subplane, $\Rho|_\cR{\,=\,}\1{\,=\,}\zeta|_\cR$,  and $\Delta$ is semi-simple with trivial center, hence $\Delta$ is strictly simple, cf. \cite{cp} 33.7 or \cite{sz1} 5.2. 
If $\dim\cR{\,=\,}4$,  then $\Rho{\,\ne\,}\1$ or $\zeta{\,\ne\,}\1$ implies 
$\cS{\, \ne\,}\cR$,  $\,\dim\cS{\,=\,}2$, and  $\overline\Delta{\,=\,}\Delta|_\cS{\,\ne\,}\1$ is simple. Finally, let  
$\dim\cR{\,=\,}8$. Then $\cS{\,=\,}\langle p^\Delta\rangle{\,<\,}\cR$,   $\,\dim\cS{\,\le\,}4$, and 
$\cF_{\overline\Delta}{\,=\,}\emptyset$. Either $\dim\cS{\,=\,}2$  and $\Delta|_\cS$ is simple by what has just 
been proved, or   $\dim\cS{\,=\,}4$ and $\overline\Delta$ is semi-simple with trivial center. In the latter case 
$\overline\Delta$ is simple by \cite{cp} 71.8. \Qed
\par\medskip
{\bold 1.10 Nearly classical translation planes.} There are several families of non-classical planes 
$\cP$ admitting a group $\Delta$ of dimension at least $17$. They can be distinguished by the structure of a maximal compact subgroup $\Phi$ of $\Delta$. In each case $\dim\Phi{\,\ge\,}3$ \ 
(see \cite{hl2} Hauptsatz or \cite{cp} 82.15). 
The Hughes planes are characterized by $\Phi{\,\cong\,}\U3\CC$. As stated in 1.3, in all other cases $\cP$ is a translation plane. With the exception of the families (3) and (5) in the following list, $\cP$ is also a dual translation plane (a plane of Lenz type V) and can be coordinatized by a semi-field (= not necessarily associative division ring) $\,(H,+,\circ)$.   \par
\qquad (1) $\,\Phi{\,\cong\,}\SO3\RR$, $\,c{\circ}z{\,=\,}t{\cdot}cz{\,+\,}(1{-}t){\cdot}zc$ with 
$\frac{1}{2}<t{\,\in\,}\RR$ ({\it mutations\/} $\HH_{(t)}$), $\,\dim\Delta{\,=\,}17$, \par
\qquad (2) $\,\Phi{\,\cong\,}(\SO2\RR)^3$, $\,(H,+,\circ)$ is a Rees algebra\footnote{\ multiplication in 
$H{\,=\,}{\CC}^2$ is defined as $(a,b){\circ}(x,y){\,=\,}(ax{+}e^{i\vartheta}y\overline b,xb{+}\overline ay)$, where 
$0{\,<\,}\vartheta{\,<\,}\pi.$} (\cite{cp} 82.16),  $\,\dim\Delta{\,=\,}17$, \par
\qquad (3) $\,\dim\Phi{\,=\,}9$,  $\,(H,+,\circ)$ is a near-field (\cite{cp} 64.19\hskip1.5pt-22),  $\,\dim\Delta{\,=\,}17$, \par

\qquad (4) $\,\Phi{\,\cong\,}\U2\CC$, $\dim\Delta{\,=\,}18$:  there exist $3$ closely related  one-parameter families of 
\par \hskip45pt   planes\footnote{\ the planes of one family are  coordinatized by  semi-fields $H_a{\,= \,}(\HH,+,\circ)$ defined as follows:  let $a{\,=\,}e^{i\alpha}$ with  
$0{\,<\,}\alpha{\,<\,}\frac{\pi}{2}$ and $V{\,=\,}\{z{\,\in\,}\HH \mid z{+}\overline z{\,=\,}0\}$, and 
note that each quaternion can be written  in the form  $\sigma{\,+\,}\frs\hskip1pt a$ with unique elements $\sigma{\,\in\,}\RR$ and $\frs{\,\in\,}V$. Put  $(\sigma{\,+\,}\frs\hskip1pt a){\,\circ\,}x{\,=\,}\sigma\hskip1pt x{\,+\,}\frs\hskip1pt x\hskip1pt a$.}  
  of Lenz type V, see  \cite{cp} 82.2,  \par 

\qquad (5)  $\,\Phi{\,\cong\,}\U2\CC$,  $\dim\Delta{\,=\,}17$: there is a $2$-parameter family of  translation planes such
\par \hskip45pt  that $\Delta$ contains a $3$-dimensional group of elations.  These planes  are  
described \par \hskip45pt in \cite{cp} 82.20. \par
 Obviously, $\cF_\Delta$ is a  flag for all planes of Lenz type V and for the planes of the family (5); in the near-field planes, $\Delta$ has exactly $2$ fixed points.
In any case,  $\Delta$ fixes a  flag, see also \cite{cp} p.\hskip1pt500, Theorem.
\par\medskip
{\bf 1.11 Dimension formula.} By \cite{hal} or  \cite{cp} 96.10, the following holds for the action of $\Delta$ on 
$P$ or on any closed $\Delta$-invariant subset $M$ of $\cP$, and for any point $a{\,\in\,}M$:
$$\dim\Delta{\,=\,}\dim\Delta_a{\,+\,}\dim a^\Delta \quad {\rm or}  \quad 
\dim a^\Delta{\,=\,}\Delta{\,:\,}\Delta_a\,.$$
\par\medskip

{\bf 1.12 Planes of other dimensions.} All $2$-dimensional (=\,{\it flat\/}) planes admitting a group of dimension ${\ge\,}3$ are known explicitly, see \cite{sz1} and \cite{cp} Chapter 3, cf. also 9.15. \\ 
The classification \cite{cp} 74.27 of all $4$-dimensional planes with a group of dimension at least 
$7$ has been extended, mainly by Betten, to {\it flexible\/} planes
 defined by the property that the automorphism group has an open orbit in the flag space, see \cite{bt}. If $\cP$ is not classical, then $\dim\Sigma{\,\le\,}8$. With the  exception of a single  shift plane found by Knarr \cite{kn0}, all planes with a $7$-dimensional group are (dual) translation planes. (The lines of a {\it shift plane\/} are the shifts (translates) of the graph of a suitable function  $f$;  a class of examples is given by 
the maps $f{\,:\,}\CC\to\CC:z\mapsto z^2|z|^c$ ($c{\,\in\,}\CC$), the parameter  $c{\,=\,}0$~yields the classical complex plane $\cP_\CC$. Knarr's most homogeneous shift plane is given by \hskip6pt 
$x{\,+\,}iy\mapsto xy{-}\frac{1}{3}s^3{\,+\,}i(\frac{1}{2}y^2{-}\frac{1}{12}x^4)$.)  \\ 
{\tt Theorem.} {\it A semi-simple group $\Delta$ of automorphisms of a $4$-dimensional plane is in fact almost simple, and  $\dim\Delta{\,\le\,}8$ or $\Delta{\,\cong\,}\PSL3\CC$\/}, see \cite{cp} 71.8; 
this will be used repeatedly in the following.
\par\smallskip 
All $16$-dimensional planes with a group of dimension ${\ge\,}35$ have been described, provided the group does not fix exactly one point and one line (\cite{hs2}, \cite{hs}). By \cite{cp} 83.26, a group of dimension ${>\,}30$ fixes at most some collinear points and their join or two points and two lines. Only for the Hughes planes (including the classical octonion plane $\cO$), the group $\Sigma$ has no fixed element.
The most homogeneous planes other than $\cO$ have a group 
$\Sigma$ with $\dim\Sigma{\,=\,}40$, they are coordinatized by a {\it mutation\/} of the octonions 
(a semi-field with a multiplication as in 1.10\hskip1pt(1),  see \cite{cp} 87.7). If $\dim\Sigma{\,=\,}39$, then $\cF_\Sigma$ is a flag, but the corresponding planes are not known explicitly. 
Because of some distinctive features, $8$-dimensional planes cannot be treated in a completely analogous way. Main differences between the two cases  are the following: 
\par\smallskip
\halign{#\ &#\hfil&\qquad#\hfil\cr
&{\tt 8-dimensional planes} & {\tt 16-dimensional planes} \cr
\noalign{\par\smallskip\hrule\par\medskip}
1)& transitive homology groups may exist & homology groups are never transitive, \cr
&& in the octonion plane they are isom. to $\RR^{\times}$ \cr
2)& compact Lie groups on $\Ss_4$ are known & there is no  list of all compact groups on $\Ss_8$ \cr
3)& $\SO3\RR$ has representations in all odd  & $\Gtwo$ has representations only in  \cr
& dimensions ${>\,}1$ & dimensions 7, 14, $>16$  \cr} 
\par\smallskip
The groups $\SO3\RR$ and $\Gtwo$ are the largest possible stabilizers of a quadrangle. The significance of 3) is due to the fact that many proofs use the complete reducibility of the representation of a compact or semi-simple subgroup of $\Delta$  on the vector space underlying the Lie algebra 
$\frl\hskip1pt\Delta$. The irreducible representation of  $\Delta$ on a minimal normal vector subgroup plays a particular r\^ole. \\
The automorphism groups of the four classical planes have the dimensions
$C_m{\,=\,}8,\,16,\,35,$ and $78$, respectively.
For each $m{\,\le\,}4$ there is a largest number $c_m$ such that there exist non-classical planes 
$\cP$ of dimension  $2^m$ with $\dim \Aut\cP{\,=\,}c_m$.   In fact, $c_m{\,=\,}4,\,8,\,18,$ and~$40$, respectively, and $\frac{1}{2}C_m{\,\le\,}c_m{\,\le\,}\frac{1}{2}C_m{\,+\,}1$.
\par\medskip
{\bf 1.13 Construction of models.} Families of curves in $\RR^2$ which form the lines of an affine part of a flat plane can easily be found.  Early examples of non-Desarguesian planes have been obtained in this way, see Hilbert's {\it Grundlagen der Geometrie\/} \S\hskip2pt23, \cite{mm}, \cite{ts}; such models may have many automorphisms, or they may be {\it rigid\/} (the automorphism group is trivial),  cf. \cite{sst}. 
The most homogeneous flat planes with $\dim\Sigma{\,=\,}4$ are due to Moulton \cite{Mo}, see  \cite{sz1} 4.8 and
\cite{cp}  \S\,34. They can be coordinatized by a Cartesian field $(\RR,+,\circ)$, where 
{$s{\circ}x{\,=\,}sx$ for $s{\,\ge\,}0$ or $x{\,\ge\,}0$,  and  $s{\circ}x{\,=\,}skx$ for $s,x{\,<\,}0$ with a fixed real number  $k{\,>\,}1$.} \\
 It is difficult, however, to find in an  intuitive way 
 $\ell$-dimensional  surfaces in $\RR^{2\ell}$ which satisfy the axioms of an affine plane. Two other procedures can be applied in all dimensions: \\
(A)   As mentioned above, suitable modifications of the algebraic operations of one of the $4$ classical (semi-)fields (sometimes called {\it mutations\/} or {\it perturbations\/}) yield topological ternary fields. 
Most of these constructions have been determined so as to retain a fairly large group of automorphisms, they do not yield rigid planes. \\ 
A noteworthy example has been given by  Tschetweruchin \cite{ts}: on the real field, 
let $\tau(s,x,t){\,=\,}sx{+}t$ for $s{\,\ge\,}0$ and $\tau(s,x,t)^3{\,=\,}s^3x^3{+}t^3$ for $s{\,<\,}0$. 
Then $\tau$ defines the ordinary addition and multiplication (see 1.4) and 
$(\RR,+,\ )$ is the real field, but $\RR_\tau$ coordinatizes a non-Desarguesian plane; cf. also 
\S\hskip1pt13 below.  \\
(B) In some situations, a geometry can be reconstructed from a suitable group of automorphisms:  let
$\cD{\,=\,}(D,\frD)$  be a  subgeometry of $\cP$ (e.g.,  a {\it stable\/} plane as defined in 
\S\hskip2pt8)   and  suppose that $\Delta$ acts transitively on the set of flags of $\cD$.  Let $(p,L)$ be any flag in  $\cD$, and put 
$\Pi{\,=\,}\Delta_p$ and $\Lambda{\,=\,}\Delta_L$.  
Then $\cD$ is isomorphic to $(\Delta/\Pi,\Delta/\Lambda)$ with incidence given by $\Pi\alpha{\smcap}\Lambda\beta{\,\ne\,}\emptyset$.  The fact that any two distinct points are joined by a line  is expressed by  $\Pi\Lambda\Pi{\,=\,}\Delta$, uniqueness of the joining line by $\Pi\Lambda{\smcap}\Lambda\Pi{\,=\,}\Lambda{\smcup}\Pi$.  This construction  is due to Freudenthal \cite{fr} \S\,6, see also \cite{stf}; it can be used in particular to describe the {\it outer\/} part of a Hughes plane.
\par\medskip
{\bf 1.14 Smooth planes.} In 1955 H. Kneser asked whether or not it is possible to introduce 
coordinates for the points and lines of a Moulton plane such that the geometric operations of 
joining two points  by a line and intersecting two lines are differentiable, say of class 
$C^\infty$\,(=\,{\it smooth\/}). In 1971 this question was answered in the negative by Betten \cite{Bt},
but a systematic treatment of smooth planes began only in 1996 with B\"odi's Habilitations\-schrift~\cite{B1}.  Because smoothness is a local property, B\"odi developped first the foundations for a theory of smooth stable planes and then applied his results to the projective case; for details see 
\S\,9. We mention the following basic \\
{\tt Theorem:} {\it Each continuous collineation of a smooth stable plane is differentiable\/}. \\
In general, smoothness strengthens other homogeneity conditions: if $\tilde c_m$ is the largest value of $\dim\Aut\cP$ for a non-classical smooth plane $\cP$ of dimension $2^m$, B\"odi's  main result is  
\begin{center}
 \begin{tabular}{|c|} 
\hline
$\tilde c_m{\,\le\,}c_m{\,-\,}2\,.$ \\  \hline
\end{tabular}
\end{center}
Examples of non-classical smooth planes have been constructed by Otte \cite{ot} by slightly disturbing the algebraic operations of the classical (semi-)fields near the element $0$. \\
Characterizations of smooth stable or projective planes in terms of submersion and trans\-versality are given in \cite{Bi}; for the projective case see also \cite{ilp}.
\par\smallskip
The assumption that the geometric operations are even {\it analytic\/} or {\it algebraic\/} are very restrictive:
{\it The only complex analytic projective plane is the classical plane over $\CC$\/}, see \cite{br1},
\cite{kr1}, and \cite{cp} 75.1.  For the algebraic case cf. \cite{stb}. There are, however, 
non-Desarguesian real analytic projective planes of dimension $2^m{\,\le\,}8$, see \cite{Iv1}. If 
points and lines are described in the usual way by homogeneous coordinates over a classical (skew) field, incidence is given by 
$(|a|^2{+}|b|^2{+}|c|^2)(|x|^2{+}|y|^2{+}|z|^2)(ax{+}by{+}cz){\,+\,}t\hskip1pt|c|^2|z|^2cz{\,=\,}0$, where 
$t$ is a real parameter with $|t|{\,<\,}\frac{1}{9}$. These planes admit compact groups of automorphisms  of dimension $1$, $4$, or $13$, respectively. Obviously, the planes are self-dual;  in fact, they admit a polarity as defined next.
\par\medskip
{\bf 1.15 Polarities.} A {\it polarity\/} is a bijective map $\rho$ of $P{\smcup}\frL$ which interchanges 
$P$ and $\frL$ such that $x{\,\in\,}y^\rho{\,\Leftrightarrow\,}y{\,\in\,}x^\rho$; this implies that $\rho$ is an involution (if the pencil $\frL_x$ is identified with the point $x$).  Besides the familiar elliptic and the hyperbolic polarity of $\cH$ there is a third kind: the {\it planar\/} polarity is a product of the elliptic polarity with a planar involution, see \cite{cp} 13.18 and 18.28 ff. 
The absolute elements of the polarity $\rho$ are defined by $x{\,\in\,}x^\rho$; if they exist, they form a {\it unital\/}. The point set of the unital belonging to the planar~polar\-ity $\pi$ of $\cH$ is a  $5$-sphere $U$; the motion group  $\Psi{\,\cong\,}\PSU4(\CC,1)$ corresponding to~$\pi$, i.e., the~central\-izer of 
$\pi$ in  $\Aut\cH$, is doubly transitive on $U$  and transitive on the complement of~$U$. Recall that the unital of the hyperbolic polarity is a $7$-sphere, its motion group $\PU3(\HH,1)$ is 
flag-transitive on the {\it interior\/} hyperbolic plane, see \cite{cp} 13.17. The polarities of the 
classical $16$-dimensional octonion plane $\cO$ are discussed in detail in \cite{cp} \S\,18.
\par\bigskip
\break
{\Bf 2. No fixed elements}
\par\medskip
In consequence of 1.9, it is a rare phenomenon that a group $\Delta$ acts without fixed elements on a plane  $\cP$. Familiar examples are the classical motion groups, i.e., centralizers of polarities of 
$\cH$. If there exists a $\Delta$-invariant subplane of dimension $4$ or $2$, then Stiffness implies that  $\dim\Delta{\,\le\,}17$ or ${\,\le\,}11$, respectively.
\par\medskip
{\bf 2.1 Theorem.} {\it Assume that $\cF_\Delta{\,=\,}\emptyset$. If $\dim\Delta{\,\ge\,}11$, then 
$\Delta$ is a Lie group\/} by \cite{sz4}~Th.~1.1.\par
{\it If $11{\,<\,}\dim\Delta{\,<\,}35$, then $\Delta$ is one of the $3$ motion groups of $\cH$, or 
$\Delta'{\,\cong\,}\SL3\CC$ and $\cP$ is a Hughes plane\/}.
\par\smallskip
{\tt Proof.}   From 1.9(c) it follows that $\Delta$ is a semi-simple group with trivial center, i.e., a direct product of strictly simple Lie groups, or there is a $\Delta$-invariant subplane $\cC$ and \cite{cp} 71.8 implies $\Delta|_\cC{\,\cong\,}\PSL3\CC$. \\
(a) In the first case, suppose that $\Delta{\,=\,}\Gamma{\times}\Psi$, where $\Gamma$ is a factor of minimal dimension and $\dim\Psi{\,\ge\,}6$. Discussion of the action of $\Psi$ on the set of fixed elements of a conjugacy class of involutions in $\Gamma$ will lead to a contradiction. \\
(a$_{\hskip-.7pt1}$)  If $\alpha{\,\in\,}\Gamma_{\hskip-1.3pt[c,A]}$ is a reflection, then $\Gamma$ is not a group of elations, and $\Gamma$ does not consist of homologies (otherwise  
$\Gamma{\,\cong\,}\SO3\RR$ by \cite{cp} 61.2, but commuting reflections have different axes by
\cite{cp} 55.35). Consequently $A^\Gamma{\,\ne\,}A$ and, dually, $c^\Gamma{\,\ne\,}c$. As
$\Psi|_{c^\Gamma}{\,=\,}\1$, the set $c^\Gamma$ does not generate a subplane, and  $c^\Gamma$ 
is contained in a fixed line $L$ of $\Gamma$. Dually $A^\Gamma$ belongs to a pencil $\frL_z$ and 
$z^\Gamma{\,=\,}z$. There is a point $x{\,\in\,}A$ such that 
$\cF{\,=\,}\langle c^\Gamma,A^\Gamma,x\rangle$ is a subplane. We have 
$\dim\Psi_{\hskip-1.5pt x}{\,>1\,}$ and $\Psi_{\hskip-1.5pt x}|_\cF{\,=\,}\1$. Therefore $\cF$ is flat, 
$\dim\Psi_{\hskip-1.5pt x}{\,\le\,}3$ by Stiffness,  $\dim\Psi{\,<\,}8$, $\,\cF{\,=\,}\cF_{\Psi_x}$,
 $\,\cF^\Gamma{=\,}\cF$ and  $\dim\Gamma{\,=\,}6$, which is impossible. \\
(a$_{\hskip-.5pt2}$) Now let $\beta$ be a planar involution in $\Gamma$. Different planes 
$\cF_{\hskip-2pt\beta}$ and $\cF_{\hskip-2.5pt\beta}^\gamma$ with $\gamma{\,\in\,}\Gamma$ are not disjoint  (\cite{cp} 55.38 or  \cite{lw2}), and \cite{cp} 71.8 shows that $\dim\Psi{\,\le\,}8$.
It follows that $\cF_{\hskip-2pt\beta}$ is the classical complex plane and that $\Psi$ acts on this plane in the standard way without fixed elements, see \cite{cp} 72.\hskip1pt1,\hskip1pt3, and 4. \ 
As $\dim\Psi{\,\le\,}8$, we have $\dim\Gamma{\,\ge\,}6$, and there is an involution 
$\alpha{\,\in\,}\Cs\Gamma\beta$. The action of $\Psi$ implies that $\alpha|_{\cF_{\hskip-2pt\beta}}$  
is not a reflection. Consequently $\cE{\,=\,}\break\cF_{\hskip-2pt\alpha,\beta}$ is a flat subplane and 
$\Psi{\,\cong\,}\SL3\RR$. Suppose that $\cE^\gamma{\,\ne\,}\cE$ for some 
$\gamma{\,\in\,}\Gamma$. Then the $\Psi$-orbits $\cE$ and $\cE^\gamma$ are disjoint.  
A point $x$ in $\cE$ is the only fixed point of $\Psi_{\hskip-1.5pt x}$ in $\cE$, and $x^\gamma$ is a  fixed point of  $\Psi_{\hskip-1.5pt x}$ in $\cE^\gamma$. Let $z$ be another point in $\cE$.
Analogously,  $\Psi_{\hskip-1.5pt z}$ fixes each point of  the orbit $z^\Gamma$, and the fixed points 
of the $4$-dimensional stabilizer   $\Psi_{\hskip-1.5pt x,z}$ form a connected subplane. This contradicts Stiffness.  Consequently $\cE^\Gamma{\,=\,}\cE$ and $\Gamma|_\cE{\,=\,}\1$, but this is impossible.  Hence a product of two or more factors has dimension at most $11$, and 
$\Delta$ is simple in the first case. \\
(b) In the second case,  $\Delta|_\cC$ is covered by a subgroup $\SL3\CC$ of $\Delta$, 
see \cite{sz6} (3.8). The asser\-tion is now a consequence of the following far more general 
theorem of Stroppel  \cite{str}~4.5:
\par\medskip
{\bf 2.2 Theorem.} {\it If $\Delta$ is almost simple and if $\dim\Delta{\,>\,}10$, then $\Delta$ is isomorphic to one of the groups $\PSL3\HH,\;\PU3(\HH,r)$ with $r{\,\in\,}\{0,1\},\;\PSU4(\CC,1),\;\SL3\CC$, or 
$\SL2\HH$. Either $\cP$ is classical, or $\Delta{\,\cong\,}\SL3\CC$ and $\cP$ is a Hughes plane\/}. \\
If $\Delta{\,\cong\,}\SL2\HH$, the center of $\Delta$ contains a reflection $\zeta$, and $\Delta$ fixes center and axis of $\zeta$. This completes the proof of 2.1.
\par\medskip
{\bold 2.3 Theorem.} {\it Assume that  $\cF_\Delta{\,=\,}\emptyset$ and that $\Delta$ has a 
minimal normal vector subgroup. If $\dim\Delta{\,>\,}6$, then $\Delta$ 
is a Lie group  and $\dim\Delta{\,\le\,}10$\/}. {\it In the case of equality, 
$\Delta{\,=\,}\Lambda{\times}\Psi$, where 
$\Lambda{\,\cong\,}{\rm L}_2{\,=\,}\{t{\,\mapsto\,}at{+}b{\,:\,}\RR{\,\to\,}\RR\mid 
a{\,>\,}0,\,b{\,\in\,}\RR\}$ 
fixes all elements of a real subplane $\cE$ and  $\Psi{\,\cong\,}\SL3\RR$ induces the full 
collineation group on~$\cE$\/}.
\par\smallskip
{\tt Remarks.} No example with $\dim\Delta{\,=\,}10$ is known, and 
probably none does exist. For a {\tt proof} of 2.3 and further properties of the hypothetical case 
$\dim\Delta{\,=\,}10$ see \cite{sz7}.
\par\bigskip
{\Bf 3. Exactly one fixed element}
\par\medskip
{\bf 3.1 Theorem.} {\it Suppose that $\cF_\Delta$ consists of a unique line and that $\Delta$ is a semi-simple group of dimension $>\!3$.  Then $\Delta$ is a Lie group and $\dim\Delta{\,\le\,}10$. In the case of equality, $\Delta$ is isomorphic to  $\Opr5(\RR,1)$ or to some covering group of 
 $\Opr5(\RR,2)$\/}.
\par\smallskip
{\tt Proof.} Except for the last part,  this has been proved in \cite{sz4} 1.3 and 2.1; for almost simple groups the assertion follows from Stroppel's theorem \cite{str}. If $\dim\Delta{\,=\,}10$, then $\Delta$ is locally isomorphic to a group $\Opr5(\RR,r),\ r{\,\le\,}2$. By \cite{cp}  55.29, a central involution 
$\zeta$ of $\Delta$ is either planar or a reflection. In the first case, $\dim\Delta|_{\cF_\zeta}{\,=\,}10$ 
contrary to \cite{cp} 71.8; in the second case, the center of $\zeta$ would be a fixed point of 
$\Delta$. By Richardson's theorem $(\dagger)$, \,(see \cite{cp} 96.34), the compact form cannot act on the fixed line of $\Delta$. \Qed
\par\medskip
{\bf 3.2 Theorem.} {\it If $\cF_\Delta{\,=\,}\{W\}$ and if $\Delta$ has a normal torus subgroup, then 
$\dim\Delta{\,\le\,}13$ and this bound is sharp\/}.
\par\smallskip
{\tt Proof.} The torus is even central (cf. \cite{cp} 93.19). Hence there exists a central 
involution~$\zeta$, and 
$\zeta$ is planar, because the center of a reflection would be a fixed point of $\Delta$.
Let $\Lambda$ denote the kernel of the action of $\Delta$ on $\cF_\zeta$. 
Then $\dim\Lambda{\,\le\,}1$ by Stiffness, and 
$\cF_\zeta$ is the classical complex plane whenever $\dim\Delta{\,\ge\,}10$, see \cite{cp} 72.8  or  1.12 above. Consequently, $\Delta{\,:\,}\Lambda{\,\le\,}12$. In a proper Hughes plane, the stabilizer of
a line of the invariant Baer subplane ({\it an affine Hughes group\/}) is $13$-dimensional. \Qed
\par
{\tt Remark.}  There is little chance that all planes admitting an affine Hughes group can be described explicitly. In the case of $4$-dimensional planes, an  analogous situation has been considered.
In this case all planes admitting a ($6$-dimensional) affine Hughes group have been determined 
in the papers \cite{kkl} and \cite{lw3}.
\par\medskip
{\bf 3.3 Groups with a minimal normal vector subgroup $\Theta$.} If $\cP{\,\not\cong\,}\cH$, then $\Delta$ has a fixed element by 2.1,   and  the last statement of 1.10 implies that $\dim\Delta{\,\le\,}16$. Let
$\cF_\Delta{\,=\,}\{W\}$ and suppose that $\dim\Delta{\,\ge\,}15$. A long and rather involved proof shows that 
$\Theta$ is contained in the group $\Tau{\,=\,}\Delta_{[W,W]}$ of translations with axis $W$, see \cite{sz7} Theorem~2. \par
{\it If $\dim\Delta{\,=\,}16$, then $\Tau$ is transitive and $\Delta$ has a subgroup $\Gamma{\,\cong\,}\SL2\CC$\/}.
\par
Transitivity of $\Tau$  has been announced with a short sketch of  proof in  \cite{sz3}~(5); a detailed proof is given in Boekholt's dissertation \cite{bh} Satz 2.2. According to Satz 
7.11 in \cite{bh}, a Levi complement of the radical  $\sqrt\Delta\,$ (a maximal semi-simple subgroup of $\Delta$) is isomorphic to  $\SL2\CC$. All translation planes admitting the group $\Gamma$ have been described explicitly by H\"ahl \cite{hl3}. The planes with $\dim\Delta{\,=\,}16$ 
depend on an integer and a real parameter in $(0,1)$. There is a $2$-sphere $S{\,\subset\,}W$ such that $\Delta$ is transitive on $W\sm S$ and $\Gamma$  acts doubly transitively on $S$.
\par\bigskip
\break
{\Bf 4. Two fixed elements} 
\par\medskip
{\bf 4.1 Theorem.} {\it Assume that $\Delta$ is semi-simple and that $\cF_\Delta{\,=\,}\{v,W\}$ is a flag. If $\dim\Delta{\,>\,}3$, then $\Delta$ is a Lie group and $\dim\Delta{\,\le\,}11$; in the unlikely case of equality,  
$\,\Delta$ is simply connected and the smaller factor consists of $(v,W)$-translations\/}.
\par\smallskip
{\tt Remarks.} By \cite{sz4} Theorem 1.3, a semi-simple group of dimension ${\!\!>\,}3$ is a Lie group if it fixes a unique line $W$ and no point outside $W$. In the classical plane, a semi-simple sub\-group of the stabilizer of a 
flag has dimension at most $9$; the same is true for all translation planes by 1.10. For almost simple groups the assertion follows from Stroppel's theorem \cite{str} 4.5. No example is known with $\dim\Delta{\,>\,}9$. For a proof of the theorem see \cite{sz4} 3.2.
\par\medskip
{\bf 4.2 Proposition.} {\it If $\Delta$ has a normal torus subgroup and if $\cF_\Delta$ is a flag, then 
$\dim\Delta{\,\le\,}11$\/}.
\par\smallskip
{\tt Proof.} The involution $\iota$ of the normal torus is in the center of $\Delta$, cf. \cite{cp} 93.19, hence $\iota$ is not a reflection and the fixed elements of $\iota$ form a Baer subplane $\cF_\iota$ containing $v$ and~$W$.  
By  \cite{cp} 71.7 and the dimension formula, we have $\dim\,\Delta|_{\cF_\iota}{\,\le\,}10$,  and Stiffness implies that the kernel of the action of $\Delta$ on $\cF_\iota$ has dimension ${\,\le\,}1$. Equality holds for the stabilizer of an {\it inner\/} flag in a Hughes plane. \Qed 
\par\medskip
{\bf 4.3 Remark.} Each of the  planes of Lenz type V described in  1.10 provides an example of a group $\Delta$ with a normal vector subgroup such that $\cF_\Delta$ is a flag. Conversely, if  $\cF_\Delta$ is a flag and if 
$\dim\Delta{\,\ge\,}17$, then $\cP$ is a translation plane.
\par\medskip

{\bf 4.4 Theorem.} {\it Let $\cF_\Delta{\,=\,}\{o,W\}$ with $o{\,\notin\,}W$. For the sake of simplicity, assume that lines are manifolds.  If $\Delta$ is semi-simple and if $\dim\Delta{\,>\,}13$, then $\cP$ is  classical. If 
$\dim\Delta{\,=\,}13$, then  $\Delta{\,\cong\,}\U2(\HH,r){\,\cdot\,}\SU2\CC$  with $r{\,\in\,}\{0,1\}$, the smaller factor consisting of homologies\/}.
\par\smallskip
{\tt Proof.} Because of 1.10, it suffices to consider only groups $\Delta$ with 
$13{\,\le\,}\dim\Delta{\,<\,}17$.\\
(a) For almost simple groups, the assertion is an immediate consequence of Stroppel's theorem  \cite{str}.  If there exists a $\Delta$-invariant proper subplane $\cS$, then $\dim\Delta|_\cS{\,\le\,}8$, or $\cS$ is a Baer subplane, $\Delta$ acts almost effectively on $\cS$ by Stiffness, and  \cite{cp} 71.8 implies again  $\dim\Delta{\,\le\,}8$. Hence each central involution is a reflection with axis $W$ and center~$o$. In fact, each element in the center of $\Delta$ acts trivially on $W$, or the stabilizer $\Delta_x$ of some point $x{\,\notin\,}W$  
would fix a quadrangle, and then $\dim\Delta{\,\le\,}12$ by Stiffness and the dimension formula. \\
(b) Suppose first that the center of $\Delta$ is trivial. Then $\Delta{\,=\,}\Gamma{\times}\Psi$, where $\Gamma$ is a factor of minimal dimension. If some involution $\iota{\,\in\,}\Gamma$ is a reflection, then $\Psi$ fixes center and axis of  $\iota$ and of each of its conjugates. Either  
$\iota|_W{\,\ne\,}\1$, $\,\iota$ and some conjugate fix more than $2$ points on $W$, and the stiffness result 1.5(7) implies $\dim\Psi{\,<\,}7$  and $\dim\Delta{\,<\,}13$, 
 or $\Gamma|_W{\,=\,}\1$,   $\,\Gamma{\,\cong\,}\Spin3\RR$ by \cite{cp} 61.2 and 1.5(8), and $\iota$ would be central. Hence $\iota$ is planar,  $W{\,\approx\,}\Ss_4$ holds automatically, and $6{\,\le\,}\dim\Gamma{\,\le\,}\dim\Psi{\,\le\,}8$ by part (a) of the proof. Richardson's  theorem $(\dagger)$ implies that $\Delta$ has torus rank $2$; therefore  $\rk\Gamma{\,=\,}\rk\Psi{\,=\,}1$. In particular, $\Gamma$ is isomorphic to $\SO3\CC$ or 
$\SL3\RR$, and  $\Gamma$ has a subgroup $\Phi{\,\cong\,}\SO3\RR$. By 1.5(8), the fixed 
elements of $\Phi$ form a flat  subplane $\cE$ containing $o$ and $W$, and $\dim\Psi|_\cE{\,\le\,}4$, a contradiction.\par
(c)  Assume next that $\Delta$ is transitive on $W$. Then  $\Delta|_W$ contains a transitive subgroup $\SO5\RR$ by  \cite{cp} 96.19\hskip1pt--22.  The group  $\SO5\RR$ cannot act on $\cP$ because there are too many involutions 
(\cite{cp} 55.40),  hence it is covered by a subgroup $\U2\HH{\,\cong\,}\Spin5\RR$  of~$\Delta$, denoted for short by the Cyrillic letter  {\cyss Yu}.
The rank argument shows that the smaller factor $\Gamma$ of $\Delta$ consists of homologies with axis~$W$, and then $\Gamma{\,\cong\,}\Spin3\RR{\,\cong\,}\SU2\CC$. If $\dim\Delta{\,>\,}13$, then {\cyss Yu} is properly contained in $\Psi$, and the list \cite{cp} 94.33 of almost simple Lie groups shows that $\Psi{\,\cong\,}\SL2\HH$ or $\dim\Psi{\,\ge\,}20$. In both cases $\cP$ is classical. \\
(d) The last part of step (a) implies that $\overline\Delta{\,:=\,}\Delta|_W{\,=\,}\overline\Gamma{\times}\overline\Psi$ is a direct product of strictly simple Lie groups. 
In this step, the assumption $W{\,\approx\,}\Ss_4$ will be used by applying $(\dagger)$.
There are two possibilities: either $\dim\Gamma{\,=\,}3$, 
$\hskip1pt10{\,\le}\dim\Psi{\,<\hskip1pt}14$, $\hskip1pt\rk\overline\Psi{\,=\,}2$,  $\,\overline\Gamma{\,=\,}\1$, and then
$\,\Gamma{\,\cong\,}\Spin3\RR$ is a group of homologies as in step (b), and $\dim\Psi{\,=\,}10$, or ,
 secondly, $\dim\Gamma{\,\ge\,}6$, $\,\dim\Psi{\,\ge\,}8$, 
$\,\rk\overline\Gamma{\,=\,}\rk\overline\Psi{\,=\,}1$,  $\,\overline\Psi{\,\cong\,}\SL3\RR$, and 
$\overline\Delta$ has a subgroup $(\SO3\RR)^2$, but this contradicts  $(\dagger)$. Therefore 
$\Gamma{\,\cong\,}\Spin3\RR$. If $\Psi$ has a subgroup $\Phi{\,\cong\,}\SO3\RR$, then the fixed elements of $\Phi$ form a flat subplane 
$\cE$, and $\Gamma$ acts as a group of homologies on $\cE$, but this is impossible. By  
$(\dagger)$, a maximal compact subgroup $\Kappa$ of $\Psi$ fixes some point $z{\,\in\,}W$, and 
$\Kappa{\,\le\,}\Psi_{\hskip-2pt z}{\,<\,}\Psi$. 
According to \cite{cp} 96.34,  $\,\Kappa$ is even a maximal subgroup of $\Psi$. Consequently, 
$\Kappa{\,=\,}\Psi_{\hskip-2pt z}$, 
$\,\dim\Kappa{\,=\,}6$ and $\Psi$ is isomorphic to $\Opr5(\RR,1)$ or to its universal covering group $\U2(\HH,1)$, but in the first case $\,\Psi$ would have a subgroup $\SO3\RR$. If $\dim\Psi{\,>\,}10$, then $\dim\Psi{\,\ge\,}15$, and $\cP$ is classical.
\par\smallskip
{\tt Remark.} The theorem has been proved in a different way without assuming  lines to be manifolds in \cite{sz6} (3.1--4) and \cite{sz4} 3.3.\Qed
\par\medskip
{\bold 4.5 Theorem.} {\it If  $\cF_\Delta{\,=\,}\{o,W\}$ with $o{\,\notin\,}W$,  if $\Delta$ has a normal torus subgroup $\Theta{\,\cong\,}\SO2\RR$, and if $\dim\Delta{\,\ge\,}13$, then the plane $\cP$ is classical\/}.
\par\smallskip
{\tt Remark.} The proof  given in \cite{sz4} 3.4 can be  simplified as follows. \\
(a) The group $\Theta$ is in the center of $\Delta$, see \cite{cp} 93.19. If $z^\Theta{\,\ne\,}z$ for some 
$z{\,\in\,}W$, then $\Delta_z$ fixes a connected subset of $W$, and $9{\,\le\,}\dim\Delta_z{\,\le\,}7$ by the dimension formula and Stiffness. Hence $\Theta$ consists  of homologies with axis $W$ and center $o$. \\
(b) {\it The stabilizer $\Lambda{\,=\,}\Delta_\fre$ of a frame $\fre{\,=\,}(o,e,u,v)$ with $u,v{\,\in\,}W$ satisfies  $\dim\Lambda{\,\le\,}1$\/}. In fact, $\Theta$ acts as a group of homologies on the plane 
$\cE{\,=\,}\cF_\Lambda$ of the fixed elements of $\Lambda$,  and $\dim\cE{\,>\,}2$, see \cite{cp} 32.17. Stiffness implies  $\dim\Lambda{\,\le\,}1$. \\
(c)  Consequently, $\dim\Delta_{u,v,w}{\,\le\,}5$ for any $3$ points $u,v,w{\,\in\,}W$. Let  
$V{\,=\,}v^\Delta,\ u,w{\,\in\,}V$,  and $\nabla{\,=\,}\Delta_{u,v}$. By the dimension  formula, 
either $\Delta$ is transitive on $W$, or $\Delta$ is doubly transitive on some $3$-dimensional  
orbit $V$ and $\dim w^\nabla{\,>\,}1$.  In the second case, the classification of $2$-transitive 
groups as summarized in \cite{cp} 96.\hskip1pt15\hskip1pt--18 shows that $V$ is a sphere 
$\Ss_3$ and that  $\Delta|_V{\,=\,}\Delta/\Kappa{\,\cong\,}\Opr5(\RR,1)$. \\
(d) If $\Delta$ is transitive on $W$, then $\Delta$ has a compact subgroup  
{\cyss Yu}${\cong\,}\U2\HH$ as 
in 4.4(c)~above. The connected component $\Gamma$ of the centralizer 
${\rm Cs}_\Delta$\hskip-2pt{\cyss Yu} consists of homologies with axis~$W$, because each 
$z{\,\in\,}W$ is an isolated fixed point of {\cyss Yu}$\hskip-2pt_z$.
 Let $\Eta$ be a maximal compact subgroup of $\Gamma$. Then $\Theta{\,\trianglelefteq\,}\Eta$ and 
$\rk\Eta{\,\le\,}1$. Therefore $\Eta{\,=\,}\Theta$ and $\dim\Eta{\,=\,}1$, see \cite{cp} 94.31(c),  moreover, $\dim\Gamma/\Eta{\,\le\,}1$ by \cite{cp} 61.2, and $\dim\Gamma{\,\le\,}2$. The representation of {\cyss Yu} on the Lie algebra $\lie\hskip-1pt\Delta$ is completely reducible.  
Hence $\dim\Delta/\Theta\!${\cyss Yu}${\,\ge}5$ and  $\dim\Delta{\,\ge\,}16$. 
The group $\Theta\!${\cyss Yu}$\hskip-1pt_v$ fixes a unique second point $u{\,\in\,}W$;  
because of $(\dagger)$ it induces on $ov$ a group $\U2\CC$. Consequently there exists a 
group $\Spin3\RR$ of homologies with axis   $ov$, in particular, $\Xi{\,=\,}\Delta_v$ contains a reflection $\sigma$ with center $u$.  By 1.7 we have $u^\Xi{\,\ne\,}u$ and then $\dim u^\Xi{\,\ge\,}3$. 
Now 1.6 implies that $\sigma^\Xi\sigma{\,=\,}\Epsilon{\,\cong\,}\RR^4$ is a transitive group of 
$(v,ov)$-elations. Each group  $\Epsilon^\delta$ is also transitive, in other words, $\cP$ is a plane of Lenz tye III, and then $\cP$ is classical by  \cite{sz8} or \cite{cp} 64.18. The group $\Delta$ is doubly transitive on $W$, and $\Delta'{\,\cong\,}\SL2\HH$. \\
(e) If $\Delta$ is doubly transitive on $V{\,\subset\,}W$,  then $\Delta$ contains a covering group 
$\Omega$ of $\Delta|_V{\,=\,}\Delta/\Kappa$ by  \cite{cp} 94.27.  Obviously,  $\Kappa$ acts freely on $P\sm (o{\smcup}W)$ and $3{\,\le\,}\dim\Kappa{\,\le\,}4$.  From the action of $\Omega$ on the Lie algebra $\lie\hskip-1pt\Delta$ we conclude that $\Kappa{\,\le\,}\Gamma{\,:=\,}\Cs\Delta\Omega$. On the other hand, the arguments in step (d) (with $\Omega$ instead of {\cyss Yu}) show that 
$\dim\Gamma{\,\le\,}2$, a contradiction. \Qed
 \par\medskip
{\bold 4.6 Theorem.} {\it If $\Delta$ has a minimal normal vector subgroup, if 
$\cF_\Delta{\,=\,}\{o,W\}$ with $o{\,\notin\,}W$, and if $\dim\Delta{\,\ge\,}15$, then the plane $\cP$ is classical\/}.
\par\smallskip
{\tt Remark.} The proof is similar to the previous one; however, fact (b) above is no longer true. For this reason, 
the larger bound for $\dim\Delta$ is required. Details can be found in \cite{sz7}.
\par\bigskip

{\Bf 5. Three fixed elements}
\par\medskip
Asssume in this chapter that $\Delta$ fixes two points $u,v$, the line  $W{\,=\,}uv$, and possibly further points on $W$, but no other line. If $\Delta$ is semi-simple and if $\dim\Delta{\,>\,}3$, then  $\Delta$ is a Lie group by \cite{sz4} 1.3.
 
\par\medskip
{\bf 5.1 Theorem.} {\it If $\Delta$ is semi-simple, then $\dim\Delta{\,\le\,}10$\/}.
\par\smallskip
{\tt Proof.} (a)  {\it $\Delta$ does not contain a reflection\/}. In fact, if $\sigma$ is a reflection, then
1.6 shows that $\1{\,\ne\,}\sigma^\Delta\sigma{\,=\,}\Tau{\,\cong\,}\RR^k$ is the normal subgroup of
translations in $\Delta$, and $\Delta$ is not semi-simple. \\
(b) {\it $\Delta$ has no subgroup $\Omega{\,\cong\,}\Spin3\RR$\/}: the involution $\beta$ of such a group would be planar, and $\Omega|_{\cF_\beta}{\,\cong\,}\SO3\RR$ would act without fixed point on $\cF_\beta$ by Stiffness and  \cite{cp} 71.10.
\par
(c) Because of (a), a maximal compact subgroup $\Phi$ of $\Delta$ acts almost effectively on $W$. 
If $\Phi{\,\ne\,}\1$, then $\Phi$ contains a planar involution, $W{\,\approx\,}\Ss_4$, and Richardson's classification $(\dagger)$ shows that $\Phi$ acts on $W$ as a subgroup of $\SO4\RR$. 
Step (b) implies that  $\Phi{\,\cong\,}\SO3\RR$ or $\Phi$ is a torus of dimension at most $2$, in particular, $\dim\Phi{\,\le\,}3$. \par
(d) {\it $\Delta$ is a product of $3$-dimensional factors\/}: if $\Delta$ has a factor $\Psi$ with 
$\dim\Psi{\,>\,}3$,  then $\Psi$ has a subgroup $\Phi{\,\cong\,}\SO3\RR$ and $\dim\Psi{\,\le\,}10$ 
by step (c) and and the list  \cite{cp} 94.33 of almost simple Lie groups; all other factors have trivial compact subgroups. 
Suppose from now on that $\dim\Delta{\,>\,}10$. It follows that  $\Delta$  has some $3$-dimensional factor $\Gamma$, but then $\Gamma$ induces a solvable group on the flat plane  $\cF_\Phi$, see 1.5(8) and \cite{cp} 33.8. \\
(e) Let $\Gamma$ denote the simply connected covering group of $\SL2\RR$. Then {\it $\Delta$ is a direct product of $4$\/} (or more) {\it factors $\Gamma_{\hskip-1.5pt\nu}{\,\cong\,}\Gamma$\/}: if some  $\Gamma_{\hskip-1.5pt\nu}$ contains an involution 
$\beta$, the product of the other factors would induce an almost simple group on the Baer subplane $\cF_\beta$,  see \cite{cp} 71.8. \\
(f) {\it $x^\Gamma{\,\ne\,}x$ for each point $x{\,\notin\,}W$ and each factor $\Gamma$ of $\Delta$\/}: if  $x^\Gamma=x$, then $\Gamma|_{x^\Delta}{\,=\,}\1$, $\,x^\Delta$ is not contained in a line, 
$\langle x^\Delta\rangle$ is a $2$-dimensional subplane, and $\dim\Delta$ would be too small. \\
 (g) It is now easy to derive a contradiction.  At most one factor of $\Delta$ consists of translations 
 (or the translation group would be commutative). If $\Gamma$ is another factor and if 
 $\Delta{\,=\,}\Gamma{\times}\Psi$, not all orbits $x^\Gamma$ are contained in a line (because of step (f)\hskip1pt). Hence, for a suitable point $x{\,\notin\,}W$, the orbit $x^\Gamma$ generates a subplane $\cC$. As $\dim\Psi_{\hskip-1.5pt x}{\,>\,}0$ and 
 $\Psi_{\hskip-1.5pt x}|_{x^\Gamma}{\,=\,}\1$, it follows that $\cC$ has dimension $2$ or $4$. 
 In the first case, $\Gamma$ would induce a solvable group on $\cC$; in the second case, 
 $\Psi_{\hskip-1.5pt x}$ would be compact.  A different proof has been given in \cite{sz4}.  \Qed  

\par\medskip
{\bf 5.2 Corollary.} {\it If $\Delta$ fixes a unique line $W$ and more than $2$ points on $W$, and if $\Delta$ is semi-simple, then $\dim\Delta{\,\le\,}9$, or $\Delta$ is an infinite covering group of $\Opr5(\RR,2)$\/}.  
\par\smallskip
{\tt Proof.} The fact that $\Delta$ has exactly $2$ fixed points  has not been used in the proof of 5.1. Therefore 
it suffices to consider a $10$-dimensional (and hence almost simple) group $\Delta$. Let $\Zeta$ denote the center of $\Delta$. Then $\Delta/\Zeta{\,\cong\,}\Opr5(\RR,r)$, $\,0{\,<\,}r{\,\le\,}2$, and $\Zeta$ does not contain  a reflection or a planar involution. Consequently, $r{\,=\,}2$, $\,\Delta$ has a maximal compact subgroup 
$\Phi{\,\cong\,}\SO3\RR$, and $\Zeta{\,\cong\,}\ZZ$.  A direct proof can be found in \cite{sz4}.\Qed 
\par\medskip
{\bf 5.3 Proposition.} {\it If $\Delta$ has a normal torus subgroup, then $\dim\Delta{\,\le\,}9$\/}. 
\par\smallskip
{\tt Proof.} $\Delta$ contains a central planar involution $\beta$, and $\dim\Delta|_{\cF_\beta}{\,\le\,}8$. \Qed
\par\medskip
{\bf 5.4 Theorem.} {\it If $\cF_\Delta{\,=\,}\{u,v,uv\}$, if $\Delta$ has a minimal normal vector subgroup $\Theta$, and if  $\dim\Delta{\,\ge\,}14$,  then $\Theta$ is a group of translations having one of the  fixed points of 
$\Delta$ as center; either $\dim\Delta{\,=\,}14$, or $\Theta$ is transitive\/}.
\par\smallskip
{\tt Remark.} For the long and complicated proof see \cite{sz7}; it shall not even be sketched here. 
Under the additional assumption that  $\Delta$ has a $6$-dimensional compact  subgroup $\Phi$ acting almost effectively on the fixed line, H\"ahl \cite{hl4},\cite{hl5} has determined all translation planes satisfying the conditions of Theorem 5.4. There are two families of such planes.\\
(1) In the case $\Phi{\,\cong\,}\Spin4\RR$ they are coordinatized  by generalized Andr\'e systems 
$(\HH,+,{\bullett\,})$ defined as follows: 
let $\phi{\,:\,}e^\RR{\,\to\,} \HH$  be any continuous map,  and put $s{\bullett}x{\,=\,}s\,x^{\phi(|s|)}$,~  
$\,0{\bullett}x{\,=\,}0$. Restricted to the affine plane, $\Phi$ consists of the maps 
$(x,y)\mapsto(ax,by)$ of $\HH^2$ with $a\overline a{\,=\,}b\overline b{\,=\,}1$.
If $\phi$ is a non-trivial homomorphism into $\{c{\,\in\,}\CC\mid c\overline c{\,=\,}1\}$, then 
 $(\HH,+,{\bullett\,})$ is a near-field and $\dim\Sigma{\,=\,}17$.\\
(2)  Let $\Phi{\,\cong\,}\SO4\RR$.  Write each quaternion in the form  $x{\,=\,}x_0{\,+\,}\frx$ with $x_0{\,=\,}{\rm Re}\;x{\,=\,}\frac{1}{2}(x{+}\overline x)$. 
Choose an arbitrary homeomorphism $\rho$ of $[0,\infty)$ with $\rho(1){\,=\,}1$ and put 
$\phi_s{\,=\,}\rho(|s|)|s|^{-1}$.  Define a multiplication on $\HH$ by 
$s{\circ}(x_0{\,+\,}\frx){\,=\,}s(x_0{\,+\,}\phi_s\hskip1pt\frx)$. The maps of $\Phi$ are given by 
$(x,y)\mapsto(x^c,ayc)$ with $|a|{\,=\,}|c|{\,=\,}1$. Without additional hypotheses there is no chance to describe all planes admitting a group as in the theorem.
\par\medskip
{\bf 5.5 Theorem.}  {\it If $\Delta$ fixes a unique line $W$ and $3$ distinct points $u,v,w$, if $\Delta$ has a minimal normal vector subgroup $\Theta$ and  if $\dim\Delta{\,>\,}12$,  then $\Delta$ contains a transitive translation group 
$\Tau$. Either  $\dim\Delta{\,=\,}13$ and a complement of  $\Tau$ is isomorphic to 
$e^\RR{\cdot}\U2\CC$,  or  $\cP$ is classical\/}. 
\par\smallskip
A {\tt proof} can be accomplished in the following steps: (a) Let $a{\,\notin\,}W$ and 
$\nabla{\,=\,}(\Delta_a)^1$. Then $5{\,\le\,}\dim\nabla{\,\le\,}7$ and $\dim a^\Delta{\,\ge\,}6$, 
$\langle a^\Delta\rangle{\,=\,}\cP$ and  $a^\Theta{\,\ne\,}a$. \quad
(b) $\langle a^\Theta,u,v,w\rangle{\,=\,}\cP$ and $\nabla$ acts faithfully and irreducibly on 
$\Theta{\,\cong\,}\RR^4$. Now $\nabla'$ contains a central reflection $\alpha$ with axis $W$ and 
$\alpha^\Delta\alpha{\,=\,}\Tau$ is a translation group   of dimension $\dim a^\Delta{\,\ge\,}6$, see 1.6. \\
(c) Each translation group $\Tau_{\hskip-1.5pt z}$ with center $z{\,\in\,}\{u,v,w\}$ contains a minimal normal vector subgroup of $\Delta$, and $\dim\Tau_{\hskip-1.5pt z}{\,=\,}4$ by step (b). Hence $\Tau{\,\cong\,}\RR^8$ is transitive, and there exists a $1$-dimensional group $\Rho$ of homologies. We may and will assume that 
$\Rho{\,<\,}\nabla$. \\
(d) Either $\nabla{\,\cong\,}\CC^{\times}\cdot\hskip1pt\SU2\CC$ and $\dim\Delta{\,=\,}13$, or 
$\dim\nabla'{\,=\,}6$, $\,\dim\nabla{\,=\,}7$, and $\nabla{\,\cong\,}e^\RR{\,\cdot\,}\SO4\RR$. As 
$\nabla$ fixes $u,v,w$, it follows that $\nabla|_W{\,\cong\,}\SO3\RR$, and $\nabla$ 
contains a transitive homology group with axis $W$. Hence   $\cP$ is  coordinatized by an associative quasi-field satisfying both distributive laws, i.e., by a (skew) field. Details of the proof can be found in \cite{sz7} Theorem~5.\Qed 
\par\smallskip
{\tt Remark.} The planes described in 1.10\hskip2pt(4) admit a $13$-dimensional group $\Delta$ as in Theorem 5.5.  The complement in $\Delta$ of the translation group consists of the maps 
$(x,y)\mapsto(hxc,hyc)$ with  $h\overline h{\,=\,}1$ and $c{\,\in\,}\CC^{\times}$, mapping a line of slope $s$ to a line of slope $hs\overline h$.
\par\bigskip
{\Bf 6. Two fixed points and two fixed lines}
\par\medskip
Assume throughout this chapter that $\Delta$ fixes a {\it double flag\/} $\;\{u,v,uv,ov\}$. Stiffness 
implies that $\dim\Delta{\,\le\,}15$.
\par\medskip
{\bold 6.1 Theorem.} {\it If $\Delta$ is semi-simple, then $\dim\Delta{\,\le\,}10$\/}.
\par\smallskip
In proving Theorem 5.1, essential use has been made of the assumption  that $\Delta$ fixes only one line. Other arguments   will be needed. In the following, we provide a different {\tt proof}  than the one given in \cite{sz7}.  Suppose  that 
$\dim\Delta{\,>\,}10$. \\
(a) If $\Delta$ is almost simple, Stroppel's theorem (2.2 above) implies that  $\cP{\,\cong\,}\cH$ is the classical quaternion plane, but in $\cH$ a maximal semi-simple subgroup $\Ypsilon$ of the stabilizer of a double flag  fixes a triangle, and $\Ypsilon$ is a compact group of dimension at most $9$. \\
(b) {\it Each reflection $\sigma{\,\in\,}\Delta$ has the axis $ov$\/}: if $\sigma$ has axis $uv$ or
center $v$, then $\sigma^\Delta\sigma{\,=\,}\Tau{\,\cong\,}\RR^k$ is a group  of translations with center $v$, see 1.6;  on the other hand, $\Tau$ is an almost simple factor of  $\Delta$, an obvious contradiction. By \cite{cp} 55.32(ii), no other reflection commutes with~$\sigma$. \\
(c) Write $\Delta{\,=\,}\Gamma\Psi$, where $\Gamma$ is an almost simple factor of minimal dimension and $\Psi$ is the product of the other factors. Suppose that some element 
$\zeta{\,\ne\,}\1$ in the center $\Zeta$ of $\Delta$ acts non-trivially on $ov$, and let 
$o^\zeta{\,\ne\,}o$. Then $\Delta_o{\,\cong\,}e^\RR\,\SO4\RR$ by 1.5\,(7), 
$\dim\Delta{\,=\,}11$ and $\dim\Psi{\,=\,}8$, moreover, $o^\Delta$ is open in $ov$, lines are manifolds and Richardson's theorem $(\dagger)$ can be applied. The compact part of $\Delta_o$ induces on $ov$ a group $\SO3\RR$, and we may assume that 
$\Gamma{\,\cong\,}\Spin3\RR$ is a group of homologies with axis $ov$. 
Now $\Psi_{\hskip-1.5pt o,z}|_{z^\Gamma}{\,=\,}\1$ for each $z{\,\in\,}S{\,:=\,}uv\sm\{u,v\}$ and 
$\Psi_{\hskip-1.5pt o}$ is sharply transitive on $S$. Consider the transitive action of $\Psi$ on the 
$1$-dimensional orbit space $S/\Gamma$. Brouwer's Theorem \cite{cp} 96.30 implies that
$\Psi|_{S/\Gamma}{\,=\,}\Psi/\Nu$ with $5{\,\le\,}\dim\Nu{\,<\,}8$, but $\Psi$ is almost simple.
Therefore $\Zeta|_{ov}{\,=\,}\1$. \\
(d) We will prove that $\dim\Psi{\,\ge\,}8$. Suppose that 
$\dim\Psi{\,\le\,}6$. Then both $\Gamma$ and $\Psi$ are locally isomorphic to $\SL2\CC$. Let 
$\sigma$ be a central involution in $\Delta$. If $\sigma$ is planar, $\Delta$ acts almost effectively 
on $\cF_\sigma$, but $\Delta|_{\cF_\sigma}$ is almost simple by \cite{cp} 71.8, see 1.12. Hence $\sigma$ is a reflection, and $\sigma$ has axis $ov$ by step (b). Moreover, the center of 
$\Delta$ has order at most~$2$, and $\Delta$ has a subgroup $\SO3\RR$ containing planar involutions. This implies that lines are manifolds. As homology groups have dimension at most $4$, the group $\Delta$ acts almost effectively on $ov$. By $(\dagger)$, a maximal compact subgroup $\Phi$ of $\Delta$ induces on $ov$ the group $\SO4\RR$. On the other hand, $\Phi{\,\cong\,}\SO4\RR$ 
and the central involution is a reflection with axis $ov$. This contradiction shows that  
$\dim\Psi{\,\ge\,}8$.  \\
(e)  Assume  that  $o^\Gamma{\,\ne\,}o$ for some $o$. Then the subplane 
$\cF{\,=\,}\langle o^\Gamma,z^\Gamma,u\rangle$ is $\Gamma$-invariant and 
$\Psi_{\hskip-1pt o,z}|_\cF{\,=\,}\1$ for each $z{\,\in\,}S{\,=\,}uv\sm\{u,v\}$. If  $\cF$ is flat, $\Gamma$ would be solvable (\cite{cp} 33.8). If $\cF$ is a Baer subplane, lines are manifolds by 1.2. 
In the case $\cF{\,=\,}\cP$, the stabilizer  $\Psi_{\hskip-1pt o}$ is transitive on $S$. 
Hence~$uv{\,\approx\,}\Ss_4$. As $\Zeta|_{ov}$  is trivial,  $\overline\Delta{\,:=\,}\Delta|_{ov}$ 
is a direct product  $\overline\Gamma{\times}\overline\Psi$ of (strictly) simple groups, and 
$(\dagger)$ implies that both factors have torus rank $1$. Therefore $\Psi{\,\cong\,}\SL3\RR$, 
and a maximal compact subgroup  $\Phi$ of $\overline\Delta$ is  isomorphic to $\SO2\RR{\times}\SO3\RR$, but according to  $(\dagger)$  there is no faithful action of  $\Phi$ on  $\Ss_4$ with a 
fixed point. This  contradiction proves that 
$\Gamma|_{ov}{\,=\,}\1$, and then $\Gamma{\,\cong\,}\Spin3\RR$. 
The same is true for each $3$-dimensional factor of $\Delta$.  Consequently $\Gamma$ is unique,
and each other factor of $\Delta$ has dimension at least $6$. In fact, $\dim\Psi{\,\in\,}\{8,12\}$. \\
(f) If $z^{\Gamma\Psi}{\,=\,}z^\Gamma{\,=:\,}Z{\,\approx\,}\Ss_3$ for some $z{\,\in\,}S$, then 
$\Psi_{\hskip-1.5pt z}|_Z{\,=\,}\1$ and  $\Psi_{\hskip-1.5pt z}$ is normal in $\Psi$.
Since $\,\dim z^{\hskip-1pt\Psi}{\,=\,}\Psi{:}\Psi_{\hskip-1.5pt z}{\,\le\,}3$,  the last 
part of step (e) implies  $\Psi|_Z{\,=\,}\1$. It follows that $\Psi$ is sharply $2$-transitive on
$ov\sm\{v\}{\,\approx\,}\RR^4$, but then $\Psi$ cannot be semi-simple, cf. \cite{cp} 96.16. 
Now \cite{cp} 96.11a shows that $\dim z^\Delta{\,=\,}4$ for each $z{\,\in\,}S$, and then 
$uv{\,\approx\,}\Ss_4$ and $\Delta$ is transitive on $S$, see \cite{cp} 53.2. Therefore 
$\Psi$ is transitive on the orbit space $S/\Gamma{\,\approx\,}\RR$ in contradiction to Brouwer's 
theorem \cite{cp} 96.30. \Qed
\par\medskip

{\bf 6.2 Theorem.} {\it If $\Delta$ has a normal subgroup $\Theta{\,\cong\,}\SO2\RR$, then 
$\dim\Delta{\,\le\,}13$.  In the case of equality, $\cP$ is the classical quaternion plane\/}.
\par\smallskip
{\tt Proof.} (a)  Assume that $\dim\Delta{\,\ge\,}13$ and recall that $\Theta$ is contained in the center of $\Delta$ \  (\cite{cp} 93.19). It follows that $\Theta$ is a group of homologies with axis $ov\,$: 
if $o^\Theta{\,\ne\,}o$, then  $\Delta_o|_{o^\Theta}{\,=\,}\1$ and $9{\,\le\,}\dim\Delta_o{\,\le\,}7$ by 1.5\,(7), which is absurd. Moreover, $\Delta$ is transitive on  $ov\sm \{v\}\,$: if not, then 
$\dim o^\Delta{\,\le\,}3$ for a suitable point $o$, and Stiffness implies 
$\Delta_{o,c}{\,\cong\,}e^\RR\;\SO4\RR$ for $o{\,\ne\,}c{\,\in\,}o^\Delta$, but then $\Theta$ is not a subgroup of  $\Delta_{o,c}$. \\
(b) If $\Lambda$ fixes a (non-degenerate) quadrangle, then $\Theta$ acts on $\cF_\Lambda$ and  
$\cF_\Lambda$ is not flat. Hence $\dim\cF_\Lambda{\,\ge\,}4$ and Stiffness shows that 
$\dim\Lambda{\,\le\,}1$. Consequently,  $\dim\Delta_o{\,\le\,}2{\cdot}4{\,+\,}1$ and 
$\dim\Delta{\,=\,}13$. Let $w{\,\in\,}S{\,:=\,}uv\sm \{u,v\}$ and put 
$\Omega{\,=\,}\Delta_w$. Then $\dim\Omega{\,=\,}9$ and  $\Omega$ is doub\-ly transitive on 
$K{\,=\,}ov\sm \{v\}$. Therefore  $\Omega$ is an  extension of $\CC^2$ by 
$\CC^{\times}{\cdot\,}\SU2\CC$, cf. \cite{cp} 96.16,19 and 95.6. Moreover, $\nabla{\,=\,}\Delta_o$ is transitive on the set $M$ of all points not on any of the lines $uv, ou, ov$. \\
(c) Consider a maximal compact subgroup $\Phi$ of $\nabla$ and note that $M$ is homotopy 
equi\-valent to  $\Ss_3{\times}\Ss_3$. The exact homotopy sequence (\cite{cp} 96.12) shows that 
$\pi_3\Phi{\,\cong\,}\ZZ^2$, and we conclude that $\Phi$ has two compact $3$-dimensional factors, cf. \cite{cp} 94.36. Because there is no  $\Theta$-invariant flat subplane, $\Phi$ has no subgroup 
$\SO3\RR$, and $\Phi'{\,\cong\,}\Spin4\RR$. Recall from (a) that $\Theta|_{ov}{\,=\,}\1$. 
By Richardson's theorem $(\dagger)$, the two factors of $\Phi'$ consist of homologies with axis $uv$ or $ou$, respectively, say $\Phi_0|_{uv}{\,=\,}\1$ and $\Phi_{\hskip-1pt\infty}|_{ou}{\,=\,}\1$. 
By 1.6 and the last fact, $\Delta$ contains a transitive group  $\Tau{\,\cong\,}\RR^4$ of 
$(v,uv)$-translations. \\
(d) As $\Theta$ is contained in the radical $\Gamma{\,=\,}\sqrt\nabla$, the group $\Phi'$ is a maximal semi-simple subgroup of $\nabla$, $\,\dim\Gamma{\,=\,}3$, and $\Gamma$ is the connected component of  $\Cs\nabla\Phi'$. It follows that $\nabla$ induces on $\Tau$ a group 
$\nabla/\Nu{\,\cong\,}e^\RR{\cdot\,}\SO4\RR$. The kernel $\Nu$ is a group of homologies with axis $ov$, and \cite{cp} 61.2 implies $\Nu{\,\cong\,}\CC^{\times}$. Therefore 
$\Gamma{\,=\,}\Theta{\times}\Rho$ with $\Rho{\,\approx\,}\RR^2$. In the next step, we will show that $\Rho$ is commutative. \\
(e) Suppose that $\Rho'{\,\ne\,}\1$. Because $\nabla$ is transitive on $S{\,=\,}S_o$ as well as on
$S_z{\,=\,}oz\sm \{o,z\}$, $\,z{\,\in\,}\{u,v\}$, the group $\Rho$ is transitive on the $3$ orbit spaces $R_x{\,=\,}S_x/\Phi'{\,\approx\,}\RR$. Either $\Rho'|_{R_x}{\,=\,}\1$ or $\Rho$ acts faithfully on 
$R_x$, see \cite{cp} 96.30. We know already that $\dim\Rho|_{S_v}{\,=\,}\1$. Hence 
$\Rho'|_{ov}{\,=\,}\1$ and $\Rho|_{R_o}{\,\cong\,}\Rho|_{R_u}{\,\cong\,}\Rho$. Therefore 
any complement $\Pi$ of $\Rho'$ in $\Rho$ acts trans\-itively on $R_v$ and fixes orbits 
$w^{\Phi'}{=\,}w^{\Phi_\infty}{\,\subset\,} S_o$ and $b^{\Phi'}{=\,}b^{\Phi_0}{\,\subset\,} S_u$. 
Put $p{\,=\,}bv{\smcap}ow$ and $c{\,=\,}pu{\smcap}ov$. Then 
$p^{\Phi_0}{\,=\,}b^{\Phi_0}v{\smcap}ow$ and
$p^\Pi{\,\subseteq\,}p^{\Phi'}{=\,}p^{\Phi_0\Phi_\infty}{=\,}
b^{\Phi_0}v{\,\smcap\,}ow^{\Phi_\infty}$. Consequently 
$c^\Pi{\,\subseteq\,}p^\Pi u{\,\smcap\,}ov{\,\subseteq\,}\break p^{\Phi'}\hskip-1.5pt u{\,\smcap\,}ov
{\,=\,}c^{\Phi'}$.  This is impossible  because $\Phi'$ is compact. \\
(f) Now $\Rho{\,\cong\,}\RR^2$ and $\Gamma$ is in the center of $\nabla$. Let $\Pi$ be a 
one-parameter subgroup of $\nabla_{\hskip-1pt w}$ which induces on $\Tau$ the group of real dilatations (see steps (b) and (d)\,). Then $\Pi$ centralizes $\Phi'$ and hence $\Pi{\,\le\,}\Gamma$. Recall that 
$\Nu$ consists of homologies with axis $ov$. The orbit $w^\Nu$ is a connected set which joins the 
points $u$ and $v$, and $\Pi|_{w^\Nu}{\,=\,}\1$. According to  Richardson's theorem, the action of 
$\Phi'$ on $S$ is equivalent to the standard linear action of $\Spin3\RR$ on $\RR^4\sm \{0\}$.
This implies that $w^\Nu$ meets each orbit $s^{\Phi'}$ with $s{\,\in\,}S$. Therefore 
$\Pi|_{s^{\Phi'}}{\,=\,}\1$, and $\Pi\Phi_0$ is a transitive group of homologies with center $o$. 
Analogously, the homology group $\Delta_{[v,ou]}$ is transitive. \\
(g) By the last statement in step (c) and 1.4,  the affine plane $\cP\sm uv$ can be coordinatized by a Cartesian field  $(\HH,+,\bullett)$. The transitive homology groups have the form
$$\{(x,y)\mapsto(x, b\bullett y) \mid  0{\,\ne\,}b{\,\in\,}\HH\}\quad {\rm and}\quad 
\{(x,y)\mapsto(x\bullett\hskip.5pt c, y\bullett c) \mid  0{\,\ne\,}c{\,\in\,}\HH\}.$$
Consequently, the multiplication $\bullett$ is associative and satisfies both distributive laws; hence 
$(\HH,+,\bullett)$ is a skew field, and then $(\HH,+,\bullett)$ is isomorphic to the quaternion skew field. \Qed \\
 This proof is slightly simpler than the one given in \cite{sz4} 5.2.
 \par\medskip
{\bf 6.3 Theorem.} {\it If $\Delta$ has a normal vector subgroup $\Theta$ and if  
$\dim\Delta{\,\ge\,}14$, then the translation group $\Tau$ with center $v$ and axis $uv$ is isomorphic to $\RR^4$ and $(\Delta_o)'$ is  a $9$-dimensional compact group. In the case 
$\dim\Delta{\,=\,}15$ the plane is classical\/}.
\par\smallskip 
{\tt Remarks.} The last statement follows easily from 1.7 and the arguments of step (g) in the 
proof of 6.2.  A proof of the other assertions is given in \cite{sz7} Th.~7.   
\par\medskip
{\bf 6.4 Theorem.} {\it A plane $\cP$ has a group $\Tau{\hskip1pt\cdot\,}(\Delta_o)'$ as in\/} 6.3 
{\it if, and only if, it can be coordinatized 
by a topological Cartesian field  $(\HH,+,\circ)$ of so-called\/} distorted quaternions {\it defined as follows\/}: 
{\it  let $(\RR,+,\ast,1)$ be any topological Cartesian field  such that identically
$(-r){\,*\,} s{\,=\,}-(r{\,*\,} s){\,=\,}r{\,*\,}(-s)$. Put  
$c\circ z{\,=\,}(|c|{\ast}|z|) |cz|^{-1}cz$ for  $z{\,\in\,}\HH^{\times}$, and  $\,c\circ 0{\,=\,}0\circ z{\,=\,}0$\/}.
\par\smallskip
{\tt Remarks.} (1) An analogous theorem for $16$-dimensional planes has been proved in \cite{hs2}. 
In particular, distorted octonions have been defined in the same way. The proof that these 
form a topological Cartesian field can easily be adapted to the quaternion case; in fact, the distorted quaternions are a Cartesian subfield of the distorted octonions over the same real Cartesian field. \\
(2) Suppose that $\cP$ is coordinatized by a Cartesian field $(\HH,+,\circ)$ as in the theorem.
If $c{\,\in\,}\HH'{\,=\,}\{h{\,\in\,}\HH \mid h\overline h{\,=\,}1\}$, then $c{\,\circ\,}z{\,=\,}cz$, 
$\,\{(x,y)\mapsto(axc,byc)\mid a,b,c{\,\in\,}\HH'\}$ is a compact group of dimension $9$, and the maps $(x,y)\mapsto(x,y{+}t)$ are translations.\\
(3) See \cite{sz7} for a proof that $\cP$ can be coordinatized by distorted quaternions if 
$\Tau{\,\cong\,}\RR^4$ and $\Delta$ has a $9$-dimensional compact subgroup. \par
(4) Examples of Cartesian fields $(\RR,+,\ast)$ with unit element $1$  and associative multiplication satisfying the identity above can be obtained as follows: let $\rho$  be a homeomorphism of 
$[\hskip1pt0,\infty)$ with  $1^\rho{\,=\,}1$, and put 
$r{\,*\,}s{\,=\,}{\rm sgn}\,r{\cdot\,}{\rm sgn}\,s{\cdot}(|r|^\rho{\cdot}|s|^\rho)^{\rho^{-1}}\!$. In this case, $\circ$ is also associative,   $\{(x,y){\,\mapsto\,}(ax,y)\mid a{\,\in\,}\HH^{\times}\}$ is a transitive group of homologies, 
and $\dim\Delta{\,\ge\,}14$. 
\par\bigskip

{\Bf 7. Characterizations}
\par\medskip
{\bf 7.1 Corollary.} {\it  Assume that $\cP$ is not classical. If $\Delta$ is semi-simple  and if 
$\dim\Delta{\,>\,}11$, then $\cP$ is  a Hughes plane or $\cF_\Delta{\,=\,}\{o,W\}$ with 
$o{\,\notin\,}W$ and $\dim\Delta{\,\le\,}13$\/}.
\par\smallskip
This is a consequence of theorems 2.1, \ 3.1, \ 4.1, \ 4.4, \ 5.1,\  and 6.1. \Qed
\par\medskip
{\bf 7.2 Corollary.} {\it  Assume that $\cP$ is not classical. If $\Delta$ has a normal torus subgroup  and if  $\dim\Delta{\,\ge\,}13$, then $\cP$ is  a Hughes plane or $\cF_\Delta$ is a single line and 
$\Delta$ is an\/} affine {\it Hughes group\/} (up to duality).
\par\smallskip
This follows from 2.1, \ 3.2, \ 4.2, \ 4.5, \ 5.3,\ and 6.2. \Qed
\par\medskip

{\bf 7.3 Summary.} The next table is to be read as follows: \ $s$-$s$\,=\,semi-simple 
\par\smallskip 
\begin{tabular}{l l}
$\dim\Delta\ge b \Rightarrow \cP$ is known, &
$\dim\Delta\ge b' \Rightarrow \cP$ is a translation plane, \\
$\dim\Delta\ge b''\Rightarrow\cP$ is a  Cartesian plane,\hspace{30pt} &
$\dim\Delta\ge b^*\Rightarrow \cP$ is a Hughes plane, \\
$\dim\Delta\ge c \Rightarrow \cP$ is classical, &
$\dim\Delta\le d$, \   $\dim\Delta\ge g \Rightarrow \Delta$ known.
\end{tabular}
\par\smallskip
\begin{center}
 \begin{tabular}{|c||c|c|c|c||l|} 
\hline
$\cF_\Delta$ & $\Delta$ $s$-$s$  & $\TT\triangleleft \Delta$ & $\RR^t\triangleleft\Delta$   & 
$\Delta$ arbitrary & References   \\  \hline 
$\emptyset$ & $b^*{=}12\;^{1)}$ & $b^*{=}10$ & $d{\,=}10\;^{2)}$ & $b{\,=}12$\quad$c{\,=}18$ & 
2.1,\enspace\cite{sz4},\enspace2.3    \\
$\{W\}$ & $d{\,=}10\;^{3)}$ & $d{\,=}13$ & $b{\,=}16$ & $b'{=}16$\quad$c{\,=}17$  & 
3.1,2,3,\enspace1.10   \\ 
flag & $d{\,=}11\;^{4)}$ & $d{\,=}11$ & $b'{=}17$& $b{\,=}17$\quad$c{\,=}19$ & 
4.1,2,\enspace\cite{sz3},\enspace1.10 \\ \hline
$\{o,W\}$ & $c{\,=}14\;^{5)}$ & $c{\,=}13$ & $c{\,=}15$ &  \hspace{34pt} $c{\,=}15$ & 4.4,5,6 \\  
$\langle u,v\rangle$ & $d{\,=}10$ & $d{\,=\;}9$ & $b''{=}15$ & $b{\,=}17$\quad$c{\,=}18$ & 
5.1,3,4,\enspace1.10(3)   \\
$\langle u,v,w\rangle$ & $d{\,=}10\;^{6)}$ & $d{\,=\;}7$ & $b'{=}13\;^{7)}$ & \hspace{34pt} $c{\,=}14$ & 5.2,3,5 \\  \hline
$\langle u,v,ov\rangle$ & $d{\,=}10$ & $c{\,=}13$ & $b{\,=}14$ &\hspace{34pt} $c{\,=}15$  & 6.1,2,4 \\  
$\langle o,u,v\rangle$ & $d{\,=\;}9$ & $d{\,=\;}9$ & $d{\,=}11$ & $d{\,=}11$\quad$g{\,=}11$ & 
\cite{sz4}\,6.1,\enspace1.7 \\
arbitrary & $b^*{=}14$ & $b^*{=}14$ & $b'{=}17 $ & $b{\,=}17$\quad$c{\,=}19$ & 7.1,2, \quad \cite{sz3} \\  \hline
\end{tabular}

\par\medskip
\begin{tabular}{l l l}
$^{1)}$ & $\dim\Delta{\,=\,}11\Rightarrow\Delta{\,\cong\,}\SL3\RR{\times}\SO3\RR$ & 
\cite{sz4} Cor.\,2.3\hspace{9pt} \\
$^{2)}$ & $\dim\Delta{\,=\,}10\Rightarrow\Delta{\,\cong\,}\SL3\RR{\times}{\rm L}_2$ & 
\cite{sz7} Th.\,1 \\
$^{3)}$ & $\dim\Delta{\,=\,}10\Rightarrow\Delta{\,\cong\,}\Opr5(\RR,1){\,\smor\,}\Delta/\Zeta{\,\cong\,}\Opr5(\RR,2)$ & 3.1  \\

$^{4)}$ & $\dim\Delta{\,=\,}11\Rightarrow\Delta/\Zeta{\,\cong\,}\PSU3(\CC,1){\times}\PSL2\RR{\,\and\,}
\pi\Delta{\,=\1\,}$\hspace{9pt} & \cite{sz4} Th.\,3.2 \\
$^{5)}$ & $\dim\Delta{\,=\,}13\Rightarrow\Delta{\,\cong\,}\U2(\HH,r){\,\cdot\;}\SU2\CC$ & 
\cite{sz4} Th. 3.3 \\
$^{6)}$ & $\dim\Delta{\,=\,}10\Rightarrow\Delta/\Zeta{\,\cong\,}\Opr5(\RR,2){\,\and\,}\pi\Delta{\,=\1\,}$ &
\cite{sz4} Th.\,4.1 \\
$^{7)}$ & $\dim\Delta{\,=\,}13\Rightarrow\Delta{\,\cong\,}\RR^4{\rtimes\,}e^\RR\U2\CC$ & 
\cite{sz7} Th.\,5 \\
\end{tabular} 
\end{center}
\par\bigskip
\break
{\Bf 8. Stable planes}
\par\medskip 
The classical hyperbolic planes or, more generally, open (convex) parts of the classical af\-fine planes are typical stable planes.  A survey of these geometries can be found  in \cite{gl} \S\S~3 and 5.
Precisely, a {\it stable plane\/} $\cS$ is a pair $(P,\frL)$ such that \par
\quad (1) $P$ is a locally compact space of positive (covering) dimension $\dim P{\,<\,}\infty$, \par
\quad (2) each line $L{\,\in\,}\frL$ is a closed subset of $P$, and $\frL$ is a Hausdorff space, \par
\quad  (3) any two distinct points $x,y{\,\in\,}P$ are joined by a unique line $xy{\,\in\,}\frL$,  \par
\quad  (4) the set $\{(K,L) \mid K{\,\ne\,}L,\;K{\smcap}L{\,\ne\,}\emptyset\}$  is an open 
subspace $\frO{\,\subset\,}\frL^2$ ({\it stability\/}), \par
\quad  (5) the maps $(x,y)\mapsto xy$ and $(K,L)\mapsto K{\smcap}L$ are continuous where defined.
\par
Note that the stability axiom  excludes space-like geometries. A stable plane may not admit an open embedding into a compact  projective plane. An easy example is $(\RR^2,\frE)$, where $\frE$ consists  of all translates of the exponential line  $y{\,=\,}e^x$ and of all ordinary straight lines of the real plane except those of positive slope. Suppose that  $(\RR^2,\frE)$ is embedded into a flat {\it projective\/} plane.  Then equally many lines through two distinct points $p,q$ are disjoint from a given line $H$, but this is false, if $H$ has the equation $y{\,=\,}0$ and $p,q$ are on different sides of $H$; for details and further examples see \cite{st9}  \S\S\,2--4.
\par\smallskip
{\bold 8.0.} {\it The point space $P$ and the line space $\frL$ of a stable plane are metrizable and have a countable basis, every pencil is arcwise connected and locally connected\/};   cf. \cite{gl} 3.5 and 3.6(ii). Hence 
{\it the restriction to a connected component of $P$ is again a stable plane\/}.
\par\smallskip
Several other results on compact projective planes can be extended to stable planes. The most prominent one is due to  L\"owen \cite{Lw}, see also \cite{gl} 3.\hskip1pt28,\hskip1pt29,\hskip1pt25(c):
\par\smallskip
{\bf 8.1 Topology.} {\it Each line pencil $\frL_p$ in a stable plane $(P,\frL)$ is compact and  homotopy equivalent to a sphere 
$\Ss_\ell$, each line $L$ satisfies $\dim L{\,=\,}\ell\,|\,8$, and $\dim P{\,=\,}2\ell$. A closed subset $C$ 
of $L$ contains interior points if, and only if, $\dim C{\,=\,}\ell$\/}.
\par\medskip
{\bold 8.2 Compactness.} {\it A stable plane $\cS{\,=\,}(P,\frL)$ is projective if, and only if, $P$ is compact.
The set $\frC$ of all compact lines is open in $\frL$, and 
$(\frC,\raise1.5pt\hbox{$\scriptstyle\bigcup$}\frC)$ is a stable plane, a sort of ``dual'' 
of~$\cS$\/}. ({\it Note that $\frC$ may be empty\/}.)
\par\smallskip
{\tt Proof.}  A compact line $K$ meets each line in the pencil $\frL_p,\, p{\,\notin\,}K$, because the perspectivity  $K{\,\to\,}\frL_p$ has a compact open image and $\frL_p$ is connected.  
See also \cite{gl} 3.11,12.\Qed
\par\medskip
Let $\Gamma{\,=\,}\Aut\cS$ be the group of all continuous collineations of the stable plane 
$\cS{\,=\,}(P,\frL)$. As in the projective case, the following holds by \cite{lw4} (2.9):
\par\medskip
{\bf 8.3 Groups.} {\it The compact-open topologies of $\,\Gamma$ on $P$ and on $\frL$ coincide. 
In this topology, the group $\Gamma$ is locally compact, and $\Gamma$ acts as
a topological transformation group on $P$ as well as on $\frL$\/}.
\par\smallskip 
A collineation $\gamma$ of a stable plane $(P,\frL)$ is said to be {\it central\/} with center 
$c{\,\in\,}P$, if $\gamma|_{\frL_c}{\,=\,}\1$; dually $\gamma$ is $axial$ with axis $A{\,\in\,}\frL$, if 
$\gamma|_A{\,=\,}\1$.  Other than in the projective case, an axial collineation of a stable plane need not have a center, and the axis of a central collineation  may be removed. 
As usual, the group of collineations with center $c$ will be denoted by $\Gamma_{\hskip-1pt[c]}$, 
dually, $\Gamma_{\hskip-1pt[A]}$ is the group of collineations with axis $A$, and
$\Gamma_{\hskip-1pt[c,A]}{\,=\,}\Gamma_{\hskip-1pt[c]}{\smcap}\Gamma_{\hskip-1pt[A]}$. 
{\it Reflections\/}  are central or axial collineations of order $2$. 
If $c{\,\in\,}A$, the elements of $\Gamma_{\hskip-1pt[c,A]}$ are called {\it elations\/}.
\par\smallskip
In sharp contrast to the projective case, transitivity of $\Gamma$ on the point set of a stable plane is a rather weak condition; it is satisfied in each translation plane and, more generally, in each subplane 
$(O,\frL_O)$ of a stable plane $(P,\frL)$, where $O$ is an open orbit of $\Gamma$ in $P$ and  
$\frL_O{\,=\,}\{L{\smcap}O\mid L{\,\in\,}\frL,\, L{\smcap}O{\,\ne\,}\emptyset\}$. Line-transitivity of 
$\Gamma$ seems to be  somewhat stronger;  in the {\it affine\/} case it implies that the 
plane is classical (\cite{cp} 63.1 or \cite{sz12}).
Flag-transitivity of $\Gamma$ on a {\it stable\/}  plane $\cS$ is possible only in classical situations 
(cf.  \cite{lw5} 1.4\,,\ \cite{gl} 5.8). In his proof, L\"owen \cite{lw5} shows that the stabilizer  
$\Gamma_{\hskip-2pt p}$ of a flag-transitive group 
$\Gamma$ of $\cS$ contains a central reflection $\sigma_p$, where  
$\sigma_{\hskip-1pt p}^{\hskip.8pt\gamma}{\,=\,}\sigma_{p^\gamma}$. 
Via $\Gamma_{\hskip-2pt p}\delta{\,\mapsto\,}p^\delta$, the point set $P$ inherits a differentiable structure from the coset space 
$\Gamma/\Gamma_{\hskip-2pt p}$ of the Lie group $\Gamma$. The product 
$p^\gamma{\cdot\hskip1pt}p^\delta{\,=\,}
p^{\hskip-.3pt\gamma\hskip.5pt\sigma_{\hskip-1pt p}^{\hskip.5pt\delta}}$ turns $P$ into a symmetric space;  the reflections or {\it symmetries\/}  $\sigma_x$  are diffeomorphic automorphisms of this space, they generate the 
{\it motion group\/}. The tangent translation plane 
${\rm T}_{\hskip-1pt p}P$ can now be used to prove that $\cS$ is classical. The  proof  given below in 
8.\,11 and 12 is quite different.
\par\medskip
{\bf 8.4 Continuity.} {\it Each central or axial collineation $\gamma$ of a stable plane is continuous\/}.
\par\smallskip
A {\tt proof} follows from the fact that $\gamma$ can be reconstructed locally from center or axis, a suitable triangle and its image; cf. also \cite{lw4} (3.2).
\par\medskip

{\bf 8.5 Baer groups.} {\it If $\dim\cS{\,=\,}2\dim\cF_\Delta$, then $\Delta$ is compact\/}. 
\par\smallskip
As in the projective case, the {\tt proof} uses  the Arzela-Ascoli Theorem, see  \cite{st6} (6.1).
\par\medskip
{\bf 8.6 Semi-simple groups of $4$-dimensional planes.}  {\it If $\Delta$ is a semi-simple group of automorphisms 
of a  $4$-dimensional stable plane, then $\Delta $ is almost simple. If $\dim\Delta{\,\ge\,} 6$, then 
$\Delta$ is a Lie group and $\dim\Delta{\,\in\,}\{6,8,16\}$\/}, see \cite{lw7}, \cite{st8} 16.1, cf. also 1.12 above.
\par\medskip
{\bf 8.7 Compact groups.} {\it Let $\Delta$ be a compact group of automorphisms of
a stable plane $\cS{\,=\,}(P,\frL)$. If  $\dim x^\Delta {\,=\,}\dim P$ for some point $x$, then $\cS$ 
is a classical pro\-jective plane and $\Delta$ is isomorphic to its elliptic motion group\/}.
\par\smallskip
This a consequence of the following observations: (a)  $x^\Delta$ is compact and  open in $P$ by \cite{Lw} 11.c) or \cite{cp} 96.11a. \ (b) The subgeometry induced on $x^\Delta$ is a stable  plane, which is projective by 8.2.  \ (c)  A homogeneous projective plane is classical (\cite{lw1} or \cite{cp} 63.8).
\par\medskip
{\bold 8.8 Large  groups.} {\it Let $d_m$ denote the dimension of the automorphism group of the classical $2^m$-dimensional\/} affine {\it plane \,$(d_m{\,=\,}6,12,27,$ or $62$, respectively\/}), {\it and let $\cS$ be a $2^m$-dimensional stable plane. Up to duality,   $\cS$ is  isomorphic to the  classical affine or projective plane, or 
$\dim\Aut\cS{\,<\,}d_m$\/},  see \cite{st8} Main theorem.
\par\medskip
In addition to the notation in 8.8,  the elliptic motion group of the $2^m$-dimensional classical projective plane $\cP_m$ will be called 
$\Epsilon_m$.
\par\smallskip
{\bold 8.9 Large compact groups.} {\it If $\Delta$ is a compact connected group of automorphisms of a $2^m\hskip-2pt$-dimensional stable plane $\cS$, then $\Delta{\,\cong\,}\Epsilon_m$ or $\dim\Delta{\,<\,}d_m$. If $\Delta$ is 
locally isomorphic to $\Epsilon_m$, then $\cS{\,\cong\,}\cP_m$\/}, see \cite{st7} Main Theorem and \cite{lw6}.
\par\medskip
\break
{\bf 8.10 Semi-simple groups of $8$-dimensional planes.} {\it Let $\Delta$ be a group of automorphisms of an
$8$-dimensional stable plane. If $\Delta$ is compact, connected and almost simple, then 
$\dim\Delta{\,\le\,}10$ or $\Delta{\,\cong\,}\Epsilon_3{\,\cong\,}\PU3\HH$. If $\Delta$ is semi-simple and if 
$\dim\Delta{\,>\,}16$, then $\Delta$ is a Lie group of dimension $18, 21$, or $35$, and $\Delta$ is equivalent to a subgroup of $\PSL3\HH$\/}; see \cite{st7} (4.3) and  \cite{st8} 16.1, cf. also \cite{st2} and Theorem 2.2 above.
\par\medskip

{\bf 8.11 Flag-homogeneous  planes.} {\it Let $\cS{\hskip1pt=\hskip1pt}(P,\frL)$ be an $8$\hskip-1pt-dimensional stable plane. If there exists a flag-transitive group 
$\Delta{\,\le\,}\Aut\cS$, then $\cS$ is a classical projective or affine plane, or $\cS$ is isomorphic to the interior of the absolute sphere of the hyperbolic polarity of the classical plane\/}.
\par\smallskip

{\tt Proof.} (a) $\Delta{\hskip1pt:\hskip1pt}\Delta_p{\,=\,}\dim p^\Delta{\,=\,}8$,  and $\Delta_p$ 
is trans\-itive on the pencil $\frL_p$. As $\frL_p$ is homotopy equivalent to $\Ss_4$, it follows from \cite{cp} 96.19 that a maximal compact connected subgroup $\Phi$ of $\Delta_p$ is transitive on 
$\frL_p$, and \cite{cp} 96.23 implies that $\Phi/\Phi_{[p]}{\,\cong\,}\SO5\RR$. Hence 
$\dim\Delta{\,\ge\,}18$. The structure theorem \cite{cp} 93.11 for compact groups shows that $\Phi$ has a transitive normal subgroup $\Omega$ isomorphic to $\SO5\RR$ or to its covering group $\U2\HH$.
\par
(b) If $\Nu$ is a normal subgroup of $\Delta$ and if $\Omega{\,\le\,}\Cs{}\Nu$, then $p^\Nu{\,=\,}p$ and $\Nu{\,=\,}\1$ because $p$ is the only fixed point of $\Omega$ and $\Delta$ is transitive on $P$.
In particular, the center of $\Delta$ is trivial.
\par
(c) Suppose that $\Delta$ is semi-simple. Then $\Delta$ is a direct product of simple factors, and 
$\Omega$ is contained in one of these factors, say $\Omega{\,\le\,}\Psi$. By the previous step, the other factors are trivial, and $\Delta{\,=\,}\Psi$ is simple. According to Stroppel \cite{st1} Theorem B, the stable plane $\cS$ can be embedded into the classical projective plane $\cH$, and  
$\Delta{\,=\,}\Aut\cH$ or  $\Delta$  is one of the motion groups $\PU3(\HH,r)$ of $\cH$. Stroppel's proof rests on the behavior of involutions and excludes all other groups. In the flag-transitive case there are easier arguments: by  8.7 and 8.8, it suffices to discuss non-compact strictly simple groups $\Delta$ containing  $\Omega$ such that $18{\,\le\,}\dim\Delta{\,<\,}27$, that is, the groups 
$\SO5\CC$, $\Opr7(\RR,r),\,r{\,\in\,}\{1,2\}$, $\PSU5(\CC,1)$, $\SL5\RR$, and  $\PU3(\HH,1)$.   
In the first case and the last case but one,  $\Omega$ is a maximal subgroup  of~$\Delta$ by  \cite{cp} 96.34, 
but $\Omega{\,<\,}\Delta_p{\,<\,}\Delta$, a contradiction. If $\Delta$ is a group $\Opr7(\RR,r)$ of 
type ${\rm B}_3$, then the representation of $\Omega$ on $\frl\hskip1pt\Delta_p$ shows that 
$\dim\Cs{\Delta_p}\Omega{\,=\,}3$, but $\dim\Cs\Delta\Omega{\,=\,}2$. 
Finally, let $\Omega{\,<\,}\Upsilon{\,\cong\,}\U4\CC{\,<\,}\PSU5(\CC,1)$ and let 
$\Gamma{\,\cong\,}\SO2\RR$ be the center of $\Upsilon$. Then $p^\Gamma{\,=\,}p$, since $p$ is the only fixed point of $\Omega$, and $p^\Upsilon{\,\ne\,}p$ because $\Upsilon$ does not act on
$\frL_p$ (see \cite{cp} 96.13). Thus $\Gamma|_{p^\Upsilon}{\,=\,}\1$ and $\Gamma$ fixes each line in $\frL_p$. This implies $\Gamma{\,=\,}\Gamma_{\hskip-1pt[p]}{\,=\,}\Gamma_{\hskip-1pt[x]}$ 
for each $x{\,\in\,}p^\Upsilon$ and then $\Gamma{\,=\,}\1$. Therefore 
$\Delta{\,\not\cong\,}\PSU5(\CC,1)$, and only the hyperbolic motion group remains. \par  
(d)  If $\Delta{\,\cong\,}\PU3(\HH,1)$, then  $\Delta_p{\,=\,}\Omega\,\Cs{\hskip-1.5pt}\Omega$ 
is a maximal compact subgroup of $\Delta$; up to conjugacy, it is the same for each plane 
by the Mal'cev-Iwasawa theorem \cite{cp} 93.10.  In fact, even the pair $(\Delta_p,\Delta_L)$ 
is unique up to conjugation: the group $\Omega_L$ is conjugate to 
$\bigl\{\hskip-1pt\bigl({a \atop }{ \atop b}\bigr)\big| a,b{\,\in\,}\HH'\bigr\}{\,\cong\,}\Spin4\RR$. 
Moreover, $\Delta_{[p]}{\,\le\,}\Delta_L$ and $\dim\Delta_L{\,=\,}13$. By complete reducibility, 
$\Omega$ acts in the standard way on a complement of $\frl\hskip1pt\Omega$  in the vector space underlying the Lie algebra $\frl\hskip1pt\Delta$. Restriction to $\Omega_L$ shows that only two 
subgroups satisfy the conditions for $\Delta_L$; they are conjugate and correspond to $L$ and 
the orthogonal line $L^{\perp}$. See also step~(i) for a similar argument.  Freudenthal's 
construction (see 1.13) proves that $\cS$ is uniquely determined by $\Delta$.
\par
(e) Now let $\Theta$ be a minimal normal subgroup of $\Delta$. Then $\Omega$ 
acts non-trivially on~$\Theta$ by step (b). Hence $\Theta{\,\cong\,}\RR^t$ is a vector group and 
$\dim\Theta{\,\ge\,}5$. Minimality of $\Theta$ implies that  $\Delta$ induces an irreducible group 
$\Delta|_\Theta{\,=\,}\Delta/\Nu$ on $\Theta$, where 
$\Nu{\,=\,}\Cs{}\Theta$. Therefore $\Delta|_\Theta$ has a semi-simple subgroup of codimension at most $2$, see \cite{cp} 95.6.  Note that  $p^\Theta{\,\ne\,}p$, because $\Theta$ is normal and  
$\Delta$ is transitive on $P$.  We will show that $\Theta$ and even $\Nu$ acts freely on the point 
set~$P$.  If $\Nu_p{\,\ne\,}\1$, in particular, if $\dim\Theta{\,>\,}8$, then  
$\Nu_p|_{p^\Theta}{\,=\,}\1$, and $p^{\Theta\Omega}{\,=\,}p^{\Omega\Theta}{\,=\,}p^\Theta$ contains points of each  line in  $\frL_p$. Consequently,  $\Nu_p{\,=\,}\Nu_{[p]}{\,=\,}\Nu_{[x]}$ for
each $x{\,\in\,}p^{\Theta}$. This implies $\Nu_p{\,=\,}\1$ and then 
$\Nu_{p^\delta}{\,=\,}(\Nu_p)^\delta{\,=\,}\1$ for each $\delta{\,\in\,}\Delta$.
\par
(f) Assume next that $\dim\Theta{\,=\,}5$. The group $\Omega$ is contained in an almost simple factor $\Psi$ of the Levi complement of the radical $\Rho{\,=\,}\sqrt\Delta$, and 
$\,\Psi{\,=\,}\Omega\,$ or $\,\Psi{\,\cong\,}\SL5\RR$, see \cite{cp} 94.34. In the last case, 
$\dim\Delta{\,\ge\,}29$ and $\cS{\,\cong\,}\cH$ by 8.8. In the remaining case  
$10{\,\le\,}\dim\Delta|_\Theta{\,\le\,}11$ and $\dim\Nu{\,\le\,}8$ by step (e), so that 
$\dim\Delta{\,<\,}20$. Now   the representation of $\Omega$  on the Lie algebra 
$\frl\hskip.5pt\Nu$ shows that $\Gamma{\,=\,}\Cs{}{\Theta\Omega}$ satisfies
$2{\,\le}\dim\Gamma{\le\,}3$. As $p$ is the only fixed point of $\Omega$, it follows that 
$p^\Gamma{\,=\,}p$ and $\Gamma|_{p^{\Theta\Omega}}{\,=\,}\1$, but then 
$\Gamma{\,=\,}\Gamma_{[x]}$ for $x{\,\in\,}p^{\Theta}$, and $\Gamma{\,=\,}\1$, a contradiction.
\par
(g) If $\Theta{\,\cong\,}\RR^6$, then Clifford's Lemma \cite{cp} 95.5 implies  $\Omega{\,<\,}\Psi$, 
the list \cite{cp} 95.10 of irreducible representations shows that $\dim\Psi{\,=\,}15$, and $\Psi$ is locally isomorphic to $\Opr6(\RR,r)$ with $r{\,\le\,}1$. By 8.8, we may assume that 
$\dim\Delta{\,<\,}27$.  Put $\Gamma{\,=\,}\Cs{}\Psi$, and note that $p^\Gamma{\,=\,}p$.
From the action of $\Psi$ on the Lie algebra $\frl\hskip1pt\Delta$ we conclude 
$\dim\Delta{\,=\,}\dim\Psi{\hskip1pt+\hskip1pt}\dim\Theta{\hskip1pt+\hskip1pt}\dim\Gamma$.
Consider the natural projection $\pi{\,:\,}\Delta{\,\to\,}\Delta/\Theta$. We have
$\Psi{\smcap}\Theta{\,=\,}\1$ and $\pi{\,:\,}\Psi{\,\cong \,}\Psi^\pi$,
since $\Psi$ acts irreducibly on $\Theta$; moreover $\pi{\,:\,}\Gamma{\,\cong\,}\Gamma^\pi$,
and $\Delta^{\hskip-1.5pt\pi}{\,=\,}\Psi^\pi\Gamma^\pi$ is an almost direct product.
The image of the canonical projection 
$\rho{\,:\,}\Delta^{\hskip-1.5pt\pi}{\,\to\,}\Delta^{\hskip-1.5pt\pi}/\Gamma^\pi$
is locally isomorphic to $\Psi$.
 The restriction $\pi|_{\Delta_p}$ is injective because $\Theta_p{\,=\,}\1$.  Write
$\Chi{\,=\,}(\Delta_p)^{\pi\rho}$. Then
$$\dim\Chi{\,=\,}\dim\Delta_p{\,-\,}\dim\Gamma{\,=\,}\dim\Psi{\,+\,}\dim\Theta{\,-\,}8{\,=\,}13$$ 
 and   $\Omega^{\pi\rho}{\,<\,}\Chi{\,<\,}\Psi^{\pi\rho}$, but such a group $\Chi$ does not exist. 
\par
(h)  In the case $\Theta{\,\cong\,}\RR^7$ it follows from \cite{cp} 95.10 that $\Psi$ is the compact group $\Gtwo$ or $\dim\Psi{\,\ge}21$. If $\dim\Psi{\,=\,}14$, then $\Psi$ would act on the 
$4$-dimensional manifold $\Psi/\Omega$, but this is impossible by \cite{cp} 96.13. Hence 
$\dim\Delta{\,\ge\,}28$ and $\cS{\,\cong\,}\cH$ by 8.8.
\par
(i) By steps (a--h), we may assume that $\Delta{\,=\,}\Theta\Omega$ with $\Theta{\,\cong\,}\RR^8$ 
acting freely on $P$. The orbit $p^\Theta$ is open in $P$, cf. \cite{cp} 96.11, hence $\Theta$ is sharply transitive on $P$, and $\Omega{\,=\,}\Delta_p{\,\cong\,}\U2\HH$ is a maximal compact subgroup of $\Delta{\,\cong\,}\HH^2{\rtimes}\U2\HH$. The structure of $\Delta$ determines the
pair $(\Delta_p,\Delta_L)$  of stabilizers up to conjugation:  $\Omega_L$ is conjugate
to $\bigl\{\hskip-1pt\bigl({a \atop }{ \atop b}\bigr)\big| a,b{\,\in\,}\HH'\bigr\}{\,\cong\,}\Spin4\RR$, and 
$\Delta_L$ is one of the two $10$-dimensional non-compact subgroups containig $\Omega_L$.
Freudenthal's  construction (1.13) shows that $\cS$ is isomorphic to an affine subplane of $\cH$. \Qed
\par\medskip
An analogous result for $16$-dimensional planes can be proved in a similar way:
\par\smallskip
{\bf 8.12 Theorem.}  {\it If a $16\hskip-1pt$-dimensional stable plane $\cS{\,=\,}(P,\frL)$ admits a 
flag-transitive~group $\Delta{\,\le\,}\Aut\cS$, then $\cS$ is a classical projective or affine plane, or 
$\cS$ is isomorphic to the interior of the absolute sphere of the hyperbolic polarity of the classical octonion plane $\cO$\/}.
\par\smallskip
{\tt Proof.} (a) From \cite{cp} 96.\,19 and 23 it follows that $\Delta_p$ has a compact connected subgroup $\Phi$ such that $\Phi/\Phi_{[p]}{\,\cong\,}\SO9\RR$ is transitive on the pencil 
$\frL_p$. Hence   $\dim\Delta{\,\ge\,}52$. Note that $p$ is the only fixed point of $\Phi$. \par
\break
(b) According to \cite{cp} 94.27, there is a subgroup $\Omega{\,\le\,}\Phi$ locally  isomorphic to 
$\Spin9\RR$, and $\Omega$ is contained in a maximal almost simple subgroup $\Psi$ of $\Delta$. 
As in 8.11, each normal subgroup of $\Delta$ in $\Cs{}\Omega$ is trivial, and so is the center of 
$\Delta$. \par
(c) If $\Delta$ is semi-simple, then $\Delta$ is even simple. By 8.8 we may assume that 
$52{\,\le}\dim\Delta{\hskip1pt<\hskip1pt}62$. It follows that $\Delta$ is the classical elliptic or hyperbolic motion group of type ${\rm F\hskip-1.3pt}_4$ and dimension $52$ of the projective plane 
$\cO$, or $\dim\Delta{\,=\,}55$. The fact that $\Delta$ contains the group $\Omega$ implies
in the second case  that   $\Delta{\,\cong\,}\Opr{11}(r)$ 
with $r{\,\le\,}2$ or $\Delta{\,\cong\,}\PU5\HH$. If $\Delta$ is compact, then $P$ is compact, and  
 8.2 together with the last part of the introduction shows that $\cS{\,\cong\,}\cO$. If 
$\Delta$ is a group $\Opr{11}(r)$, the centralizer $\Gamma$ of $\Omega$ has dimension at most 
$2$, but $\dim\Delta_p{\,=\,}39$ and $\Delta_p{\,\cong\,}\Omega{\times}\Gamma$. Therefore 
$\dim\Delta{\,=\,}52$. Only the possibility $\Delta{\,\cong\,}{\rm F\hskip-1.3pt}_4(-20)$ remains. 
In this case, $\Delta_p{\,=\,}\Omega{\,\cong\,}\Spin9\RR$ is a maximal compact subgroup of 
$\Delta$, by the  Mal'cev-Iwasawa theorem \cite{cp} 93.10, it is unique up to conjugacy. 
Under the action of $\Omega_L{\,\cong\,}\Spin8\RR$ on the vector space underlying the
Lie algebra $\frl\hskip1pt\Delta$, a complement of $\frl\hskip1pt\Omega$ splits into two subspaces $\RR^8$, 
and $\Omega_L$ is contained in two conjugate $36$-dimensional subgroups $\Delta_L$ and
$\Delta_{L^\perp}$.  The claim follows from Freudenthal's construction 1.13.
\par
(d) Now let $\Theta{\,\cong\,}\RR^t$ be a minimal normal subgroup of $\Delta$. As in step (e) above, 
$\Nu{\,=\,}\Cs{}\Theta$ acts freely on the point set $P$, and $\Omega{\smcap}\Nu{\,=\,}\1$. 
In particular, $\Omega|_\Theta{\,\cong\,}\Omega$ and $t{\,\ge\,}9$.
\par
(e) Suppose that $\Theta{\,\cong\,}\RR^9$. Then $\dim\Delta|_\Theta{\,=\,}\Delta{\,:\,}\Nu{\,\le\,}37$, because $\SO9\RR$ is a maximal subgroup of $\SL9\RR$ by \cite{cp} 94.34 and the latter group is too large. Consequently $\dim\Delta{\,\le\,}53$ and $15{\,\le\,}\dim\Nu{\,\le\,}16$. The action of 
$\Omega$ on the vector space of the Lie algebra $\frl\hskip.5pt\Nu$ is completely reducible 
(\cite{cp} 95.3). Hence $\Gamma{\,=\,}\Nu{\smcap}\Cs{}\Omega$ satisfies 
$6{\,\le\,}\dim\Gamma{\,<\,}8$, see \cite{cp} 95.10. On the other hand, $p^\Gamma{\,=\,}p$ is the only fixed point of  $\Omega$, and $\Gamma{\,\le\,}\Cs{}\Theta$ implies $\Gamma|_{p^\Theta}{\,=\,}\1$. Again 
$p^{\Theta\Omega}{\,=\,}p^\Theta$ contains points of each line in the pencil $\frL_p$, and each point $x{\,\in\,}p^\Theta$ is a center of $\Gamma$. Therefore $\Gamma{\,=\,}\1$, a contradiction. \par
(f) If $9{\,<\,}t{\,<\,}16$, then $\Gamma{:=\,}\Theta{\smcap}\Cs{}\Omega{\,\ne\,}\1$, since the action of $\Omega$ on $\Theta$ is completely reducible.   As $p$ is the only fixed point of $\Omega$,  we have $p\hskip1pt^\Gamma{=\,}p$, but  $\Gamma{\,\le\,}\Theta_p{\,=\,}\1$ by step (d). \par
(g) Consequently $\dim\Theta{\,=\,}16$, and $\Theta$  is sharply transitive on $P$. We may 
assume  that $\Delta{\,\cong\,}\RR^{16}{\rtimes\,}\Spin9\RR$. As in step (c),  the pair 
$(\Delta_p,\Delta_L)$ of stabilizers of a flag is unique up to conjugation, and $\cS$ is the classical affine plane. \Qed
\par\bigskip
{\Bf 9. Smooth planes}
\par\medskip
This section is based mainly on B\"odi's thesis \cite{B1}. See also McKay \cite{mck} for further aspects.
\par\medskip
{\bf 9.0 Definition.} {\it A stable plane $\cS{\,=\,}(P,\frL)$ is said to be\/} smooth {\it if\/} \par
\quad (1) {\it $\,P$ and $\frL$ are smooth\/} (=\,C$^\infty$) {\it manifolds, and\/} \par
\quad (2) {\it\hskip2pt joining two points and intersecting two lines are smooth maps on their respective \par \hskip30pt domains of definition\/}. 
\par\medskip

{\bf 9.1 Theorem.} {\it Each line of a smooth stable plane $(P,\frL)$ is a closed submanifold of~$P$. Each pencil $\frL_p$ is a compact submanifold of $\frL$ and $\frL_p{\,\approx\,}\Ss_\ell$. 
Moreover, the flag space {\bf F}=$\{(p,L) \mid p{\,\in\,}L\}$ of  $(P,\frL)$ is a closed submanifold of
the product manifold $P{\times}\frL$\/}.
\par\medskip
The tangent space ${\rm T}_{\hskip-1pt p}P$ at the point $p$ of the manifold $P$ is a  $2\ell$-dimensional real vector space, each tangent space ${\rm T}_{\hskip-1pt p}L$ is an $\ell$-dimensional subspace, and 
$\bigcup_{L{\in}\frL_p}\hskip-3pt{\rm T}_{\hskip-1pt p}L{\,=\,}{\rm T}_{\hskip-1pt p}P$. 
\par\smallskip
{\bf 9.2 Derived planes.} {\it If $(P,\frL)$ is a smooth stable plane, the translates of all subspaces 
${\rm T}_{\hskip-1pt p}L$ with $L{\,\in\,}\frL_p$ in the vector space ${\rm T}_{\hskip-1pt p}P$ form
a locally compact  topological affine translation plane~$\cA_p$. The projective closure  
$\overline\cA_p$ of  $\cA_p$ is a compact $2\ell$-dimensional topological projective translation~plane\/}.
\par\smallskip
The translation planes $\cA_p$ are  among the most important tools in the investigation of smooth planes. The proof that  $\cA_p$ is a {\it topological\/} plane uses the following
\par\medskip
{\bf 9.3 Linearization.} {\it There exists an open neighbourhood $U$ of $p$ in $P$ and a homeomorphism $\lambda{\,:\,}U \to {\rm T}_{\hskip-1pt p}P$ such that 
$\lambda(L{\smcap}U){\,=\,}{\rm T}_{\hskip-1pt p}L$ for every line $L{\,\in\,}\frL_p$\/} \ (cf. \cite{B2} 3.12).
\par\medskip
{\bf 9.4 Regularity.} {\it The map $L\mapsto{\rm T}_{\hskip-1pt p}L$ is a smooth homeomorphism 
of the pencil $\frL_p$ onto its image; if this map is even a diffeomorphism, $(P,\frL)$ is said to be regular at $p$. Equivalently,  $\cA_p$ is a smooth  affine  plane; if the projective plane  $(P,\frL)$
is regular everywhere, then the dual plane  $(\frL,P)$ is regular\/}, see \cite{mck} Theorem 9 and Corollary 11.

\par\medskip
{\bf 9.5 Automorphisms.} {\it Every continuous collineation $\gamma{\,:\,}\cS\to\cS'$ 
 between two smooth stable planes is a smooth collineation\/} \ (see \cite{B3}). 
 \par\smallskip
 This is {\tt proved\/} in the following steps: (a) By  \cite{BK}, 
{\it every continuous local homomorphism between smooth local left or right loops  is smooth in a neighborhood of the neutral element\/}. (b) For each flag $(o,K)$, a local ternary field can be introduced similar to the global construction  in 1.4 such that  addition defines a local loop on K with neutral element $o$. \\ (c) For each line $K$, the restriction $\gamma|_K$ is locally and globally smooth. \\ (d) As $P$ is locally the product of two lines, smoothness 
of $\gamma$ follows. \Qed
 \par\smallskip
 We write $\Gamma{=}\Aut\cS$ for the group of all smooth collineations of the smooth stable plane~$\cS$. 
The group $\Gamma$, taken with the compact-open topology, is a locally compact transformation group of $P$ as well as of $\frL$ or {\bf F} \ (see 8.3).  Now \cite{mz} \S\,5.2,\;Th.\,2 together with 8.0  implies
\par\medskip
{\bf 9.6 Lie groups.}  {The group \it $\Gamma{=}\Aut\cS$ is a smooth Lie transformation group of $P$, $\frL$, and {\bf F};  moreover,
the stabilizer of a triangle is a linear Lie group, the stabilizer of a non-degenerate quadrangle is compact\/}.
\par\medskip
{\bf 9.7 Derivations.} {\it The derivations 
${\rm D}_p{\,:\,}\Gamma_{\hskip-1pt p}\to \GL{}({\rm T}_{\hskip-1pt p}P)$ and its dual ${\rm D}_L$ are continuous homomorphisms\/}, cf.   \cite{mz} \S\,5.2,\;p.208,\;Th.\,3.
\par\medskip
{\tt Lemma.} {\it  If $\gamma{\,\in\,}\Gamma$ fixes more than $2$ lines of the pencil $\frL_p$, then the fixed lines 
of $\gamma$ in $\frL_p$ form a homology $m$-sphere with $m{\,>\,}0$\/}.
\par\smallskip
This Lemma is a crucial step in a long  discussion leading to the following results:
\par\smallskip
{\bf 9.8 Elations and reflections.} {\it The elations with center $c$ form a torsion free closed subgroup  
${\rm ker\,D}_c{\,=\,}\Gamma_{\hskip-1pt[c,c]}{\,=\,}\bigcup_{A{\in}\frL_c}\Gamma_{\hskip-1pt[c,A]}$ of 
$\Gamma$.  Dually, ${\rm ker\,D}_L{\,=\,}\Gamma_{\hskip-1pt[L,L]}$, and each compact subgroup of  
$\Gamma_{\hskip-1pt[L,L]}$ is trivial. The connected components of $\Gamma_{\hskip-1pt[c,c]}$ and 
$\Gamma_{\hskip-1pt[L,L]}$ are simply connected and solvable. If $\alpha,\beta$ are reflections in 
$\Gamma_{\hskip-1pt[c]}$, then $\alpha\beta$ is an elation\/}. 
\par\medskip
{\bf 9.9 Stabilizers.} {\it If  $(P,\frL)$ is a smooth stable plane, and if  $\Delta$ is a {\it connected\/} closed subgroup of  $\Gamma{=}\Aut(P,\frL)$, then each  point stabilizer $\Delta_p$ is a 
\emph{linear} Lie group, and each semi-simple subgroup of the stabilizer of a \emph{flag} is compact\/}, \ see \cite{B4}, in particular Th. 3.7. 
\par\medskip
Let $\Xi$ denote the stabilizer 
$\Delta_{H,K,L}$ of $3$ distinct lines $H,K,L$ in some line pencil $\frL_p$. In the classical $2\ell$-dimensional projective plane we have $\dim\Xi{\,=\,}\dim\Sigma{\,-\,}5\ell{\,:=\,} k_\ell{\,=\,}3,6,15,38$.
 We say that  $(P,\frL)$ is an almost projective translation plane, if there is a compact connected projective translation plane 
$(\overline P,\overline\frL)$ such that $\overline P\sm P$ is a compact connected subset of the translation axis $W$ and $\frL{\,\subseteq\,}\overline\frL{\,\subseteq\,}\frL{\smcup\{W\}}$. 
\par\medskip
{\bf 9.10 Theorem.} {\it If $\cS$ is a $2\ell$-dimensional smooth  stable plane, then either 
$\dim\Xi{\,\le\,}k_\ell{\,-\,}\ell$ or $\cS$ is an almost projective translation plane\/} \ (see \cite{B4}).
\par\bigskip
The remainder of this section is concerned with smooth {\it projective\/} planes $\cP{\,=\,}(P,\frL)$. 
These can be characterized as compact, connected smooth local planes, see \cite{ilp}.

\par\medskip
{\bf 9.11 Topology.} {\it The point manifold of a smooth projective plane is homeomorphic to its classical counterpart, and so is the line manifold\/}, see Kramer \cite{kr1}.
\par\medskip
For the sake of simplicity, the next result is stated only for $8$-dimensional planes:
\par\medskip
{\bf 9.12 Solvable groups.} {\it If $\Delta$ is a connected\/} solvable {\it group of automorphisms of a\/} non-classical {\it $8$-dimensional projective plane, then $\dim\Delta{\,\le\,}9$, or $\Delta$ fixes a 
flag $(p,L)$ and $\dim\Delta{\,\le\,}13{\,+\,}\min\{3,\dim\Delta_{[p,L]}\}$\/}, see \cite{B7} p.\hskip2pt131.
\par\medskip
{\bf 9.13 Subplanes.} {\it  If $\cE$ is a closed subplane of the smooth projective plane $\cP$, then the  point space $E$ of $\cE$ is smoothly embedded into $P$\/} \ (\cite{B1} 8.3).
\par\medskip
{\bf 9.14 Translation planes.}  {\it Every smooth projective translation plane $\cP$ is isomorphic to one of the four classical planes\/} \ (Otte \cite{ot}, \cite{Ot}).
\par\smallskip
We may assume that $\cP$ is not flat (cf. \cite{sz1} 5.4 or \cite{cp} 33.5).  A {\tt proof} of 9.14  is achieved in the following way: 
(a) If  ${\rm Diff\,}\Ss_n$ denotes the group of diffeomorphisms of $\Ss_n{\,=\,}\RR^n{\smcup}\{\infty\}$ with respect to the standard smooth structure, then 
$\GL n\RR{\smcap}{\rm Diff\,}\Ss_n{\,=\,}{\rm GO}_n\RR$ is the group of similitudes, hence it is rather small. This observation suggests \\
(b) Suppose that the topological quasi-field $Q{\,=\,}(\RR^n,+,\circ\hskip2pt)$ coordinatizes a smooth projective translation plane $\cP$ and that $\hat Q{\,=\,}Q{\smcup}\{\infty\}$ carries the smooth structure induced by a projective line of $\cP$. Then each map
$\mu_{a,b}{\,:\,}x{\,\mapsto\,}a{\circ}x{\,-\,}b{\circ}x$ with $a{\,\ne\,}b$ extends to a diffeomorphism 
$\hat\mu_{a,b}$ of  $\hat Q$, in fact, $\hat\mu_{a,b}$ is contained in the group $\Alpha$ of affine projectivities. In \cite{Ot} it is shown that the group
$\langle \hat\mu_{a,b}\mid a,b{\,\in\,}Q, a{\ne}b\rangle{\,=\,}\Alpha_0{\,\le\,}\GL n\RR$ does not contain   $\SL n\RR$, for $n{\,=\,}4$ not even $\SL2\CC$. According to 
Grundh\"ofer--Strambach \cite{GSt}, the latter conditions characterize the classical Moufang planes. \Qed
 
\par\medskip
We write again $\Sigma{\,=\,}\Aut\cP$ for the full automorphism group.
\par\medskip
{\bf 9.15 Flat planes.} {\it If $\cP$ is a smooth flat plane and if $\dim\Sigma{\,\ge\,}3$, then $\cP$ is the classical real plane\/}.
\par\smallskip
{\tt Remarks.} (1) Moulton planes are the most homogeneous non-Desarguesian flat projective planes. They have a $4$-dimensional automorphism group $\Sigma$, and this characterizes these planes. 
$\Sigma$ fixes a non-incident point-line pair $(o,W)$, the commutator subgroup $\Sigma'$ is isomorphic to the simply connected covering group $\tilde\Omega$ of $\Omega{\,=\,}\SL2\RR$.  
In the usual description of the Moulton planes as indicated in 1.13 the line $W$ is a vertical affine line; for a discussion of affine models such that $W$ is the line at infinity see \cite{cp} \S\,34 or 
\cite{sz0} p.\,14. 
\par\smallskip
{\tt Theorem.} {\it The affine Moulton planes having the fixed line $W$ as line at infinity cannot be turned into a smooth plane\/}. \ 
(This follows from 9.9 and the fact that $\tilde\Omega$ has no faithful linear representation, 
cf. \cite{cp} 95.10). 
\par\smallskip
(2) There are 3 families of non-classical flat planes with $\dim\Sigma{\,=\,}3$, viz., the skew 
hyperbolic planes with  $\Sigma{\,\cong\,}\PSL2\RR$ and $\cF_\Sigma{\,=\,}\emptyset$, the skew parabola planes (or flexible flat shift planes, a $2$-parameter family such that $\cF_\Sigma$ is a flag), and a $3$-parameter family of Cartesian planes ($\cF_\Sigma$ is a double flag), see \cite{sz1} or 
\cite{cp} \S\S~35--37. Lengthy calculations show that the proper skew hyperbolic planes are not even C$^1$-differentiable and the skew parabola planes are not even C$^2$-differentiable.
\par\medskip
{\bf 9.16 Four-dimensional planes.} {\it If $\cP$ is a $4$-dimensional smooth projective plane, then 
$\dim\Sigma{\,\le\,}6$ or $\cP$ is classical. Moreover, the only $4$-dimensional  flexible smooth shift plane is the complex projective plane\/} \ (\cite{B1}\,Chap. 10). 
\par\smallskip
{\tt Remarks.} As mentioned in 1.12, with the exception of Knarr's single shift plane,  all $4$-dimensional planes with $\dim\Sigma{\,\ge\,}7$ are translation planes (up to duality); if smooth, the latter are classical by Otte's theorem 9.14. The generating function for a flexible shift plane is  
always differentiable. For affine points $p$, the tangent translation planes $\cA_p$ are pairwise isomorphic as topological affine planes (\cite{B1} 10.5).
\par\medskip
{\bf 9.17 Eight-dimensional planes.} {\it If $\cP$ is a smooth projective plane such that 
$\dim P{\,=\,}8$ and $\dim\Sigma{\,\ge\,}17$, then $\cP$ is the classical quaternion plane\/} \ 
(\cite{B6} Theorem B).
\par\smallskip
{\tt Proof.} According to 1.10, all $8$-dimensional planes with  $\dim\Sigma{\,\ge\,}17$ are translation planes or Hughes planes (\cite{sz3}). The smooth translation planes are classical by 9.14. Now let $\cP$ be a Hughes plane. Then $\dim\Sigma{\,=\,}17$, and $\Sigma$  induces on the invariant complex subplane $\cC{\,=\,}(C,\frC)$ the group $\Sigma/\Kappa{\,\cong\,}\PSL3\CC$. Moreover, 
$\Kappa{\,\cong\,}\SO2\RR$, and $\PSL3\CC$ is covered by a subgroup 
$\Gamma{\,\cong\,}\SL3\CC$ of $\Sigma$, cf. \cite{sz6} (3.\,6--8) or \cite{cp} \S\,86. All involutions of $\Gamma$ are reflections, their centers and axes are the points and lines of $\cC$, and $\Gamma$ is flag-transitive on the outer subgeometry of $\cP$, see \cite{cp} 86.5 for the last assertion..  \\
(a) Choose an outer flag $(p,L)$, put $\Pi{\,=\,}\Gamma_{\hskip-1.5pt p}$, 
$\,\Lambda{\,=\,}\Gamma_{\hskip-2pt L}$, and $\Delta{\,=\,}\Gamma_{\hskip-1.5pt p,L}$, and note that $\Pi$ fixes a unique line $W$ of $\cC$ and $\Lambda$ fixes a unique point $o{\,\notin\,}W$ of 
$\cC$. Let $\Omega{\,=\,}\Sigma_{o,W}$; then $\dim\Omega{\,=\,}9$, and 
$\Delta{\,<\,}\Gamma{\smcap}\Omega{\,\cong\,}\GL2\CC$. 
If $\cP{\,\cong\,}\cH$, an easy calculation with coordinates shows that 
$\Delta{\,\cong\,}e^\RR{\cdot\,}\SU2\CC$. In the case of proper Hughes planes, 
$\Delta{\,\cong\,}\Rho{\times}\SU2\CC$ with $\Rho{\,=\,}\{e^{(1{+}\alpha i)t} \mid t{\,\in\,}\RR\}$ and
$\alpha{\,\ne\,}0$, see \cite{cp} 86.\,2,12. As the involutions in $\Gamma$ are reflections, the group 
$\Tau{\,=\,}\Gamma_{\hskip-1.5pt[W,W]}$ is sharply transitive on $C\sm W$, dually, 
$\Epsilon{\,=\,}\Gamma_{\hskip-1.5pt[o,o]}{\,\cong\,}\RR^4$. Thus $\Pi{\,=\,}\Delta{\ltimes}\Tau$ 
and $\Lambda{\,=\,}\Delta{\ltimes}\Epsilon$. Freudenthal's construction (1.13) shows that the 
(everywhere dense) outer geometry is uniquely determined by $\alpha$, and so is the plane 
$\cP{\,=\,}\cH_\alpha$. For distinct positive  values of $\alpha$,  the re\-pre\-sentations of $\Delta$ on $\Tau$ are inequivalent; hence $\cH_\alpha{\,\not\cong\,}\cH_{\alpha'}$ for 
$0{\,<\,}\alpha{\,<\,} \alpha'$, cf. \cite{cp} 86.37.  
For later use note that as in \cite{cp} 86.\,2 and 35 one may write \\[-10pt] \par
\hskip9pt$\Sigma{\,=\,}\{S{\,\in\,}\GL3\CC\mid |\det S|{=}1\}$ \  and \ 
$\Sigma_p{\,=\,} \left \{ \left ( \begin{array}{cc}  C|h|&  \\ \frv&h^{-2} \end{array} \right )  
\;  \bigg | \;  C{\,\in\,}\U2\CC, \ h{\,\in\,}\Rho, \ \frv{\,\in\,}\CC^2 \right \}$.   \\[2pt]
Then $\Pi{\,=\,}\Sigma_p{\smcap}\Gamma$, and $\Lambda{\,=\,}\Pi^t$ is the transpose of~$\Pi$. \\
(b) Consider the tangent translation plane $\cA_o$ of $\cP$  at $o$ and the stabilizer $\Alpha$ of 
$\Aut\cA_o$ at~$o$.  The restriction ${\rm D}_o{\,:\,}\Omega\to\Alpha$  of the derivation ${\rm D}_o$ is injective because the kernel consists of elations only (see 9.8).  Therefore  
$\dim\Aut\cA_o{\,=\,}\dim\Alpha{\,+\,}8{\,\ge\,}17$, and  $\cA_o$ is one of the translation planes 
listed in 1.10, or $\cA_o$ is classical. If $\Phi$ is a maximal compact subgroup of $\Omega$, then
$\Phi/\Kappa{\,\cong\,}\U2\CC$ and $\dim\Phi{\,=\,}5$. Hence ${\rm D}_o\Phi$ is a maximal compact subgroup of  $\Alpha$. This excludes all proper translation planes of 1.10, and  
$\overline\cA_o{\,\cong\,}\cH$. 
The action of $\Omega$ on $\frL_o$ is equivalent to the action of ${\rm D}_o\Omega$ on the pencil 
$\frA_o$ of $\cA_o$. We may identify the derived plane $\frc_o$ of $\cC$ at $o$ with the affine plane 
$\cC\sm W$, so that  $\overline\frc_o$ is a Baer subplane of $\overline\cA_o$. Now  
$\Delta{\,\cong\,}{\rm D}_o\Delta{\,\cong\,}e^\RR{\cdot\,}\SU2\CC$ by step (a), and 
$\alpha{\,=\,}0$. \Qed
\par\medskip

{\bf 9.18 Sixteen-dimensional planes.}  {\it If $\cP$ is a $16$-dimensional smooth projective plane, then $\dim\Sigma{\,\le\,}38$ or $\cP$ is the classical octonion plane\/}, hence $\tilde c_4{\,<\,}39$  as asserted in 1.14.
\par\smallskip
{\tt Proof.} Suppose that $\cP$ is smooth. If $\dim\Sigma{\,=\,}40$, then $\cP$ has Lenz type at least V, and $\cP$ is classical by 9.14. If $\dim\Sigma{\,=\,}39$, it is only known  that 
$\cF_\Sigma$ is a flag $(v,W)$, \ see \cite{cp} 87.4 and 82.28.  Consider the derived plane 
$\cA_v$ of $\cP$ at  $v$, the stabilizer  $\Alpha$ of $\Aut\cA_v$ at~$0$, and the derivation map 
${\rm D}_v{\,:\,}\Sigma\to\Alpha$ with kernel $\Epsilon{\,=\,}\Sigma_{[v,v]}$. 
By Otte's theorem 9.14, we may assume that $\dim\Epsilon{\,<\,}16$. Now 
$39{\,-\,}15{\,\le\,}\dim\Alpha$ and  $\dim\Aut\cA_v{\,\ge\,}\dim\Alpha{\,+\,}16{\,\ge\,}40$. 
The dual  is true as well.  \\
(a) The action of $\Sigma$ on the pencil $\frL_v$ is equivalent to the action of ${\rm D}_v\Sigma$ on 
the tangent pencil in $\cA_v$.  Choose any line $K{\,\ne\,}W$ in $\frL_v$.
Assume  that $\cA_v$ is classical. Then $\Sigma_K/\Sigma_{[v]}$ can be identified with a 
subgroup of  $e^\RR{\times}\Spin8\RR$ (cf. \cite{cp} 12.17, 15.7, and 17.13). If $\Sigma$ has a compact subgroup $\Phi{\,\cong\,}\Spin8\RR$, then the center of $\Phi$ contains $3$ reflections and 
$\Phi$ fixes the centers of these reflections. Representation of $\Phi$ on the Lie algebra 
$\frl\hskip1pt\Sigma$ shows  that $\dim\Cs\Sigma\Phi{\,=\,}3$, but the stabilizer of a triangle in $\cP$ has dimension at most $30$, see \cite{cp} 83.26. Hence $\dim\Phi{\,\le\,}21$. 
As ${\rm D}_v{\,:\,}\Sigma/\Epsilon\to\Alpha$ is injective and 
$\Sigma_{[v]}/\Epsilon\to\Alpha_{[v]}{\,\cong\,}e^\RR$, we have 
$\dim\Sigma_{[v]}{\,\le\,}16$ and $\dim\Sigma_K/\Sigma_{[v]}{\,\ge\,}15$. Consequently, a maximal compact subgroup $\Phi$ of 
$\Sigma_K$ has dimension at least $14$. Either $\Phi{\,\cong\,}\Gtwo$ or $\Phi'$ is one of the groups $\Spin7\RR$ or $\SU4\CC$, cf. \cite{cp} 95.12 and 55.40. In the latter case $\cP$ is a translation 
plane by L\"uneburg \cite{Lb} V, Satz, and then $\cP$ is classical by 9.14.  Only the possibility  
$\Phi{\,\cong\,}\Gtwo$ remains. \\
(b) If $\cA_v$ is not classical, then  $\dim\Alpha{\,=\,}24$, 
$\,{\rm D}_v{\,:\,}\Sigma/\Epsilon{\,\cong\,}\Alpha$, 
and $\cA_v$ can be  coordinatized by a semi-field with multiplication 
$s\circ x{\,=\,}t{\hskip.5pt\cdot\hskip.5pt}sx{\,+\,}(1{-}t){\hskip.5pt\cdot\hskip.5pt}xs$,  
$\frac{1}{2}{\,\ne\,}t{\,\in\,}\RR$ \ ({\it mutation\/} $\OO_{(t)}$  of~$\OO$,\, see \cite{cp} 87.7). 
Obviously, $\Aut\OO_{(t)}{\,\cong\,}\Aut\OO$ and $\Alpha$ is an extension of $\RR^8$ by 
$\RR^2{\times}\Gtwo$. \\
(c) In any case, ${\rm D}_v\Phi{\,=\,}\Lambda$ is a maximal compact subgroup of  
$\nabla{\,=\,}{\rm D}_v\Sigma_K{\,\cong\,}\Sigma_K/\Epsilon$, and $\nabla$ fixes a triangle in the projective closure $\overline\cA_v$. 
Hence $39{\,-\,}\dim K^\Sigma{\,-\,}\dim\Epsilon{\,=\,}\dim\nabla{\,\le\,}16$.
It follows that $\dim\Epsilon{\,=\,}15$ and $K^\Sigma{\,=\,}\frL_v\sm\{W\}$. Dually, $\Sigma$ is transitive 
on $W\sm\{v\}$. The elation group $\Sigma_{[v,K]}{\,<\,}\Epsilon$ is $7$-dimensional, and the dual of \cite{sz11} or \cite{cp} 61.12 implies that $\Sigma_{[v,W]}$ is trans\-itive. Consequently, $\Sigma$ 
is transitive on  $\frL\sm\frL_v$. \\
(d) In $\cA_v$ there exists a one-parameter group of homologies with center~$v$. In particular, the 
pre\-image in $\Sigma$ contains a homology $\eta{\,\in\,}\Sigma\sm \Epsilon$ with center $v$ and some axis $H$. By the dual of \cite{cp} 61.20b it follows that 
$\dim H^\Sigma{\,=\,}\dim\Epsilon{\,=\,}16$, a contradiction. \Qed 
\par\smallskip
For a detailed {\tt proof} see  B\"odi \cite{B5}, Main Theorem, cf. also \cite{B6}.

\par\bigskip
{\Bf 10. Unitals of 8-dimensional planes} 
\par\medskip
Originally, a unital has been defined as the geometry of the absolute elements of a polarity (1.15). 
The term is also used more generally for subsets of a projective plane which do not correspond to 
a polarity, but have similar properties as a {\it polar\/} unital, cf. \S\hskip2pt12.  Up to conjugacy 
($\sim$), all polarities  $\rho$ of the  nearly classical planes decribed in 1.10 have been determined
by Stroppel \cite{st4}, and so have the corresponding unitals $U_{\hskip-1.3pt \rho}$ and the connected components $\Mu_\rho$ of their motion groups (cf. \cite{st3}).  Elliptic polarities (with $U{\,=\,}\emptyset$) in other planes than the classical ones seem to be rare; for examples in dimension $2$ and $4$ see \cite{pst}. Planes which are not self-dual like the near-field planes do not admit a polarity, of course. This section is based essentially on Stroppel's papers \cite{st4} and \cite{st3}.
\par\medskip
{\bf 10.1 Hughes planes.}  {\it Each polarity  of a proper  Hughes  plane $\cP$  is conjugate  to a standard polarity $\rho$ and induces on the invariant complex subplane  $\cC$ an orthogonal polarity,  $\dim U_{\hskip-1.3pt \rho}{\,=\,}5$, and  $\Mu_\rho$ is transitive on the set of {\rm outer} points of $U_{\hskip-1.3pt \rho}$\/}. 
\par\smallskip
{\tt Proof.} (a) We use the same notation as in 9.17, in particular, we  write $\Pi$ and 
$\Lambda{\,=\,}\Pi^t$ as at the end of 9.17(a). The map  $\xi\mapsto\xi^{-t}$ is an involutorial automorphism of 
$\Sigma$ as well as of $\Gamma{\,=\,}\SL3\CC$; it interchanges $\Pi\xi$ and 
$\Lambda\xi^{-t}$ and hence induces a polarity $\rho'$ on the outer subgeometry, which is dense in $\cP$. It follows that $\rho'$ extends to a polarity $\rho$ of $\cP$. \par
(b) If $\gamma{\,\in\,}\Cs{}\rho$, then 
$\Pi\xi\gamma{\,=\,}\Lambda(\xi\gamma)^{-t}{\,=\,}\Pi\xi\gamma^{-t}$,  and  
$\Gamma{\smcap}\Cs{}\rho{\,=\,}\{\gamma{\,\in\,}\Gamma\mid
\gamma\gamma^t{\,\in\,}\bigcap_{\xi{\in}\Gamma}\Pi^\xi{\,=\,}\1\}$ 
is the simple orthogonal group $\SO3\CC$, it acts faithfully on $\cC$. By construction 
$p^\rho{\,=\,}L$ and $o^\rho{\,=\,}W$. Hence $\rho|_\cC$ is the standard planar polarity, 
cf. \cite{cp} 13.18. Moreover, $\Gamma{\smcap}\Cs{}\rho$ is the connected component 
$\Mu_\rho$ of the full motion group $\Cs{}\rho$,  since 
$\Kappa{\smcap}\Gamma{\smcap}\Cs{}\rho{\,=\,}\1$, $\Kappa{\smcap}\Cs{}\rho{\,<\,}\Kappa$,  
and $\Kappa{\smcap}\Cs{}\rho$ is finite. \\
(c) Consider an arbitrary element $\sigma{\,\in\,}\Sigma{\,<\,}\GL3\CC$. 
Then $\rho\sigma$  is a polarity if, and only if, 
$(\rho\sigma)^2{\,=\,}\sigma^\rho\sigma{\,=\,}\sigma^{-t}\sigma{\,=\,}\1$ and $\sigma^t{\,=\,}\sigma$ 
is symmetric. Now $\sigma{\,=\,}\tau^2$ for a symmetric $\tau{\,\in\,}\Sigma$ \,(see \cite{st4}; instead, one may argue as in 11.3(c)\,). Consequently $(\rho\tau)^2{\,=\,}\1$ and $\rho\sigma{\,=\,}\tau^{-1}\rho\tau$ is conjugate to $\rho$ whenever $\rho\sigma$ is a polarity. \\
(d) An outer point $\Pi\xi$ belongs to $U_{\hskip-1.3pt\rho}$ if, and only if, there exists some 
$\gamma{\,\in\,}\Pi\xi{\smcap}\Lambda\xi^{-t}$. This implies $\gamma\xi^{-1}{\,\in\,}\Pi$ and 
$\gamma\xi^t{\,\in\,}\Lambda{\,=\,}\Pi^{-t}$, $\,\xi\gamma^t{\,\in\,}\Pi^{-1}{\,=\,}\Pi$, and 
$\,\gamma\gamma^t{\,\in\,}\Pi$. Conversely, it follows from $\gamma\gamma^t{\,\in\,}\Pi$ that 
$\Pi\gamma\gamma^t{\smcap}\Lambda{\,\ne\,}\emptyset$ and 
$\Pi\gamma{\,=\,}\Pi\gamma^{-t}{\,\in\,}U_{\hskip-1.3pt\rho} $.
If $\xi\xi^t{\,\in\,}\Pi$, then $\xi\xi^t{\,=\,}\pi^2$ for some $\pi{\,=\,}\pi^t{\,\in\,}\Pi$, and
$(\pi^{-1}\xi)(\pi^{-1}\xi)^t{\,=\,}\pi^{-1}\xi\xi^t\pi^{-1}{\,=\,}\1$. Thus 
$\mu{\,=\,}\pi^{-1}\xi{\,\in\,}\Mu_\rho$, $\,\Pi\mu{\,=\,}\Pi\xi$, and $\Mu_\rho$ is transitive on the set of outer points of $U_{\hskip-1.3pt\rho}$. The stabilizer $\Mu_\rho{\smcap}\Pi$ of the point $\Pi$ in $\Mu_\rho$ consists of the matrices 
$\omega{\,=\,}\left(\begin{array}{cc} C|h| & \\ \frv & h^{-2} \\ \end{array} \right)$ƒ 
with $\omega\omega^t{\,=\,}\1$ and hence $\frv{\,=\,}0, \ h{\,=\,}1, \ CC^t{\,=\,}\1$. 
Therefore $\Mu_\rho{\smcap}\Pi{\,\cong\,}\SO2\RR$, and the $\Mu_\rho$-orbit of $\Pi$ is homeomorphic to the $5$-dimensional coset space $\SO3\CC/\SO2\RR$. \Qed
\par\smallskip
{\tt Remark.} It is not known whether or not $U_\rho$ is a sphere or if $U_\rho$ is  a {\it topological\/} unital  in the sense of \S\hskip2pt12.
\par\medskip

{\bf 10.2 Mutations.} {\it The planes $\cH_{(t)}$ over the semi-fields
$\HH_{(t)}{\,=\,}(\HH,+,\circ)$ with a multiplication  
$c{\circ}z{\,=\,}t{\cdot\hskip1pt}cz{\,+\,}(1{-}t){\cdot}zc$ and 
$\frac{1}{2}{\,<\,}t{\,\in\,}\RR$ are called {\it mutations\/} of the quaternion plane~$\cH$, 
{\rm cf. 1.10}. These planes   admit exactly two classes of polarities, one with  unitals 
${\approx\,}\Ss_7$  and  $11$-dimensional  motion groups, the other with unitals 
${\approx\,}\Ss_5$ and $7$-dimensional motion groups. The unitals of the first kind intersect each secant in a  $3$-sphere\,}: {\it they are {\rm topological} unitals in the sense of \S\hskip1pt{\rm12}\/}.
\par\smallskip
{\tt Proof.} (a) {\it $\Gamma{=\,}\Aut\HH_{(t)}{\,=\,}\Aut\HH{\,\cong\,}\SO3\RR$, moreover, 
$z\mapsto z^\kappa{\,=\,}\overline z$ is an anti-automorphism of  $\HH_{(t)}$\/}: obviously, each (anti-)automorphism of $\HH$ is an (anti-)auto\-mor\-phism of  $\HH_{(t)}$. Hence we may assume that a given automorphism $\lambda{\,\in\,}\Aut\HH_{(t)}$ fixes $i,j,k$, and then~$\lambda{\,=\,}\1$. \\
(b) We use coordinates as in 1.4 and denote the line with the equation $y{\,=\,}a{\circ}x{+}b$ by
$[a,b\hskip1pt]$. As $\cH_{(t)}$ has Lenz type V, each polarity $\rho$ maps the elation center $v$ to the translation axis $uv$. Up to conjugation with a translation with center $u$,  we may always 
assume that $u^\rho{\,=\,}ov$. If $\iota$ is an an involutorial anti-automorphism of   $\HH_{(t)}$, then 
$(a,b)\mapsto[a^\iota,-b^\iota]$ defines a polarity $\rho_\iota$, which maps the line 
$c{\times}\HH$ to the slope~$(c^\iota)$. Put  $\overline\rho{\,=\,}\rho_\kappa$. 
 As $\overline a{\circ}a{\,=\,}|a|^2$, the affine part $U_{\overline\rho}\sm \{v\}$ of $U_{\overline\rho}$ consists 
 of all points   $(a,b)$ such that  $a\overline a{\,=\,}b{+}\overline b{\,=\,}2b_0$, and 
 $U_{\overline\rho}{\,\approx\,}\Ss_7$.  For  $\iota{\,=\,}(z{\,\mapsto\,}\overline z^{\hskip1pt i})$ write 
 $\pi{\,=\,}\rho_\iota$. Then 
$$a^\iota{\circ}a{\,=\,}(\overline{a'}{-}j\overline{a''})^i{\circ}(a'{+}a''j){\,=\,}
(\overline{a'}{+}a''j){\circ}(a'{+}a''j) {\,=\,}|a'|^2{\,-\,}|a''|^2{\,+\,}2(t\overline{a'}{+}(1{-}t)a')a''j\,,$$ 
and  $U_\pi\sm \{v\}{\,=\,}\{(a,b){\,\in\,}\HH^2\mid \,|a'|^2{\,-\,}|a''|^2{\,=\,}b'{+}\overline{b'}\,\and \,
(t\overline{a'}{+}(1{-}t)a')a''{\,=\,}b''\}$. Thus  $b_1$ is arbitrary, $b_0,b_2$, and $b_3$ 
depend on $a$, and $U_\pi{\,\approx\,}\Ss_5$. In particular, $\pi$ is not conjugate to 
$\overline\rho$. \\
(c) If $\rho$ is any  polarity of $\cH_{(t)}$ with  $u^\rho{\,=\,}ov$, then 
$\overline\rho\rho{\,=\,}\sigma{\,\in\,}\Sigma_{u,v,ov}$. The stabilizer $\nabla{\,=\,}\Sigma_{o,u,v}$ 
 is the direct product of $\Gamma$ and the  homology groups 
$\nabla_{[o,uv]}{\,\cong\,}\nabla_{[v,ou]}{\,\cong\,}\RR^{\times}$, the kernels of the translation group and its dual. Recall that $\dim\Sigma{\,=\,}17$.
Consequently  
$$\sigma{\,:\,}(a,b){\,\mapsto\,}(a^\gamma s, rb^\gamma s{+}t),\,
[\,c,d\,]{\,\mapsto\,}[\,rc^\gamma, rd^{\hskip.5pt\gamma}s{+}t\,]  \  \ {\rm with} \ \ 
\gamma{\,\in\,}\Gamma,\,r,s{\,\in\,}\RR^{\times},\,t{\,\in\,}\HH.$$ 
As $\rho^2{\,=\,}\1$,  we have $\overline\rho\sigma{\,=\,}\sigma^{-1}\overline\rho$. This implies first
$r\overline a^{\hskip1.5pt\gamma}{\,=\,}s^{-1}\overline a^{\hskip1.5pt\gamma^{-1}}$ for all $a{\,\in\,}\HH$  or, equivalently,  
$rs{\,=\,}1,\,\gamma^2{\,=\,}\1$, and then  $t^\gamma{\,=\,}\overline t$.
 In the case $\gamma{\,=\,}\1$ all polarities  $\rho{\,=\,}\overline\rho\sigma$ are conjugate to $\overline\rho$.
If $\gamma{\,\ne\,}\1$, then $\gamma$ is conjugate  to the map $z{\,\mapsto\,}z^i$, and $\rho$ is conjugate 
to  $\pi$. \\
(d) The elements of $\Sigma$ have the form
$$\sigma:(a,b){\,\mapsto\,}(a^hs{+}m, rb^hs{\,+\,}q{\circ}a^hs{\,+\,n)},\,
[c,d]{\,\mapsto\,}[rc^h{+}q, rd^hs{\,-\,}(rc^h{+}q){\circ}m{\,+\,}n]$$
with $h{\,\in\,}\HH'{\,=\,}\{h{\in}\HH\mid h\overline h{\,=\,}1\}, \, r,s{\,\in\,}\RR^{\times}, \, q,m,n{\,\in\,}\HH$. An easy calculation shows that $\sigma{\,\in\,}\Mu_{\overline\rho}{\;=\;}\Cs{}{\overline\rho}$ if, and only if, $r{\,=\,}s$, 
$\,q{\,=\,}\overline m$, and $|m|^2{\,=\,}n{+}\overline n$.  Hence $\dim\Mu_{\overline\rho}{\,=\,}11$. 
Similar\-ly, $\pi\sigma{\,=\,}\sigma\pi$ yields   $h{\,\in\,}\CC{\,\smor\,}h{\,\perp\,}\CC$,  
$\,r{\,=\,}s$, $\,q{\,=\,}\overline m^{\hskip1pt i}$, and $\overline m^{\hskip1pt i}{\circ}m{\,=\,}
n{+}\overline n^{\hskip1pt i}$,   and $\dim\Mu_\pi{\,=\,}7$. \\
(e) The subgroup of $\Mu_{\overline\rho}$ described in 10.3 is transitive on the set of slopes 
$\ne\infty$. It suffices therefore to consider horizontal secants and verticals. In both cases the 
secants intersect  $U_{\overline\rho}$ in a $3$-sphere. This 
 is an obvious consequence of the condition $|a|^2{\,=\,}2{\,\rm re\hskip1pt}b$. \Qed
\par\medskip
{\bf 10.3 Remark.} The solvable subgroup 
$\{(a,b){\,\mapsto\,}(a{+}m,b{+}m^\iota{\circ}a{+}n)\mid m^\iota{\circ}m{\,=\,}n{+}n^\iota\}$ of
$\Mu_\iota$ (given by the condition $h{\,=\,}1{\,=\,}r{=}s$) is sharply transitive on the affine point set 
$U_\iota\sm\{v\}$, and this is true whenever $H$ is a $4$-dimensional semi-field and $\iota$ is an involutorial anti-automorphism of $H$.
\par\medskip
The other two families of self-dual planes with a $17$-dimensional group can be treated in a similar way, see \cite{st4}, \cite{st3} for details. For $0{\,<\,}\vartheta{\,<\,}\pi$, the multiplication of the Rees algebra  $H_\vartheta{\hskip1pt=\hskip1pt}(\CC^2,+,\circ)$ is defined by 
$(a,b){\circ}(x,y){\hskip1pt=\hskip1pt}(ax{+}e^{i\vartheta}\overline by,\,bx{+}\overline ay)$. Obviously, 
$\kappa{\,:\,}(a,b){\,\mapsto\,}(a,\overline b)$ is an involutorial anti-automorphism of 
$H_\vartheta$. As before, $\kappa$ gives rise to a polarity $\hat\kappa$ of the plane 
$\cR_\vartheta$ over $H_\vartheta$, mapping the  point $(c,d)$ to the line $[c^\kappa,-d^\kappa]$.
\par\medskip
\break
{\bold 10.4  Rees planes.}  {\it All polarities of $\cR_\vartheta$ are conjugate. 
The affine part of the unital  $U_{\hat\kappa}$ is 
$U'{\,=\,}\{(\,(a,b),\,(\frac{1}{2}(a^2{\,+\,}e^{i\vartheta}b^2),{\rm re}(\overline ab){\,+\,}ir)\,) \mid 
a,b{\,\in\,}\CC, r{\,\in\,}\RR\}$ and $U_{\hat\kappa}{\,\approx\,}\Ss_5$. The $7$-dimen\-sional 
motion group is transitive on $U'$\/} \ (\cite{st4} 4.8 and \cite{st3} 5.\,5,6).
\par\medskip
{\bold 10.5 Spin planes.} The multiplication of the only semi-fields coordinatizing  self-dual planes 
$\cP$ with $\dim\Aut\cP{\,=\,}18$ has the form 
$c{\,\circ\,}z{\,=\,}cz{\,+\,}r(cz{-}zc)_1{\,=\,}cz{\,+\,}2r(c_2z_3{-}c_3z_2)$, where $r{\,>\,}0$ and the indices  refer to the coefficients of $i,\,j,\,k$. See  \cite{hl6} and Stroppel  \cite{st4}~\S\,2 and \cite{st3} \S\,3 for more details;  note that Stroppel uses the opposite multiplication of $\HH$ in the definition of 
the product $\circ$. This multiplication  defines indeed a semi-field on  $\HH{\,\cong\,}\RR^4$: 
the map $z{\,\mapsto\,}c{\,\circ\,}z$ is linear, and $c{\,\circ\,}z{\,=\,}0{\,\ne\,}c$ implies 
$cz{\,\in\,}\RR,\ \overline cczc{\,\in\,}\RR, \ zc{\,\in\,}\RR$,  hence $c{\,\circ\,}z{\,=\,}cz, \ z{\,=\,}0$, and 
multiplication by $c$ is bijective. 
Note that $z{\,\mapsto\,}z^k$ is not an automorphism of $\HH^{(r)}{\,=\,}(\HH,+,\circ)$,  neither is 
$z{\,\mapsto\,}\overline z$ an anti-automorphism, because  
$\overline{z_1\,}{\,\ne\,}\overline z_{\hskip.5pt1}$ in general. However, the product 
$\kappa{\,=\,}(z\mapsto\overline z^{\hskip.5pt k})$ of these two involutions  is an involutorial anti-automorphism of $\HH^{(r)}$.
\par\smallskip
{\tt Theorem.} {\it The Spin planes have two conjugacy classes of polarities, one has  motion groups
of dimension $7$,  the other of dimension $9$. The corresponding unitals are $5$-spheres, and the motion groups are transitive on the affine points of the unitals.\/} 
\par\smallskip
{\tt Proof.} The map 
$\hat\kappa{\,:\,}(a,b){\,\leftrightarrow\,}[a^\kappa,-b^\kappa]$ is again a polarity. The set 
$\,U_{\hat\kappa}\sm \{v\}$ of affine absolute points is 
$\{(a,b)\mid a^\kappa{\,\circ\,}a{\,=\,}b{\,+\,}b^\kappa{\,=\,}2(b_0{+}b_1i{+}b_2j)\}
\approx\HH{\times}\RR$, 
a group of motions sharply transitive on   $\,U_{\hat\kappa}\sm \{v\}\,$ is described in 10.3, 
and $\,U_{\hat\kappa}\approx\Ss_5$.  Consider the translation $\tau{\,:\,}(a,b){\,\mapsto\,}(a,b{+}n)$. The product $\hat\kappa\tau$ defines a polarity if, and only if, $n^\kappa{\,=\,}n$ or  $n_3{\,=\,}0$, and all these polarities are conjugate, since each  $\tau$ is a square. \par

The group $\Sigma{\,=\,}\Aut\cP$ is an extension of the group of  translations and elations by 
$\nabla{\,=\,}\Sigma_{o,u}$ consisting of the maps 
$$(x,y){\,\mapsto\,}(c^{-1}xa,\overline a{\,\circ\,}yad), \
[s,t]{\,\mapsto\,}[\hskip1pt\overline ascd,\overline a{\,\circ\,}tad\hskip1.5pt] \quad {\rm with} \ 
a{\,\in\,}\HH', c{\,\in\,}\CC^{\times}\hskip-2pt, \ {\rm and} \ d{\,\in\,}\RR^{\times}\hskip-2pt.$$
The motion group $\Cs\nabla{\hat\kappa}$ is given by the conditions 
$ak\overline a\overline k{\,=\,}c^2d{\,=\,}\pm1$, $a{\,\in\,}\CC{\,\smor\,}a{\,\perp\,}\CC$. 
Consequently $\dim\Mu_{\hat\kappa}{\,=\,}7$.
\par\smallskip
For $\delta{\,\in\,}\nabla$ the condition $\delta\hat\kappa\delta{\,=\,}\hat\kappa$ yields 
$\overline a(c^{-1}xa)^\kappa cd{\,=\,}x^\kappa$ or $aka{\,=\,}\pm k$. Define $\underline\delta$ by 
$a{\,=\,}k$ and $c{\,=\,}1{\,=\,}-d$.  Then $(\hat\kappa\underline\delta)^2$ is a homology with axis 
$ou$ by construction; it maps $(0,1)$ onto itself. Hence 
$\pi{\,=\,}\hat\kappa\underline\delta{\,:\,}(x,y){\,\mapsto\,}
[\hskip1pt\overline xk,-k{\circ}k\overline y\hskip2pt]$ is a polarity. 
If $\delta{\,\in\,}\Cs\nabla\pi$, the first coordinates yield  
$\overline a\hskip1pt \overline xkcd{\,=\,}\overline a\hskip1pt\overline x\hskip1pt\overline c^{\hskip1pt-1}k$ and 
$\overline c^kcd{\,=\,}c^2d{\,=\,}1$,  there is no condition on $a$. Conversely $c^2d{\,=\,}1$ implies 
that $\pi\pi^\delta$ fixes $(0,1)$ and  all points of  $ou$, and $\pi\delta{\,=\,}\delta\pi$. 
Therefore $\pi$  is not conjugate to $\hat\kappa$.  The admissible choices for $\delta$ yield all polarities of $\cP$. It follows that there are only two conjugacy classes of polarities. \\
We have 
$U_\pi\sm \{v\}{\,=\,}\{(x,y)\mid \overline xk{\circ}x{\,=\,}y{\,+\,}k{\circ}k\overline y{\,=\,}
y{\,-\,}\overline y{\,+\,}2ry_1\}{\,\approx\,}\HH{\,\times\,}\RR$, because $y{\,+\,}k{\circ}k\overline y$ 
does not involve $y_0$,  so that again $U_\pi{\,\approx\,}\Ss_5$. By \cite{st3} 3.8 there is also a sharply transitive group of motions on the affine part of $U_\pi$, and $\dim\Mu_\pi{\,=\,}9$. 
\par\medskip
Planes with a smaller group may also admit polarities:
\par\medskip
{\bold 10.6 Fixed double flag.}  Let the plane $\cP$ be coordinatized by a Cartesian field  $(\HH,+,\circ)$ of distorted quaternions as in 6.4, where $(\RR,+,\ast)$ is commutative. Suppose that $\Sigma{\,=\,}\Aut\cP$ is connected and  $13$-dimensional. {\it Then $\cP$ has exactly two conjugacy classes of polarities, one with  unitals  $\approx\Ss_7$ and $9$-dimensional motion groups, the other with unitals 
$\approx\Ss_5$ and $5$-dimensional motion groups\/}.
\par\smallskip
{\tt Proof.}  Each map $z{\,\mapsto\,}\overline z^g$ with 
$\overline g{\,=\,}\pm g$ is an involutorial anti-automorphism. Let the polarity $\rho$ be defined by 
$(x,y){\,\mapsto\,}[\hskip1pt\overline x,-\overline y\hskip1.5pt]$. The stabilizer 
$\nabla{\,=\,}\Sigma_0$
consists of the auto\-morphisms  $\delta{\,:\,}(x,y){\,\mapsto\,}(axc,byc),\,
[\hskip1pt s,t\hskip1.5pt]{\,\mapsto\,}[\hskip1pt bs\overline a,btc\hskip1.5pt]$ with 
$a,b,c{\,\in\,}\HH'{\,=\,}\{h{\,\in\,}\HH\mid h\overline h{\,=\,}1\}$, and $\Sigma{\,=\,}\nabla\Tau$ with 
$\Tau{\,=\,}\{(x,y){\,\mapsto\,}(x,y{+}n)\mid n{\,\in\,}\HH\}$.
A product $\rho\delta\tau$, $\delta{\,\in\,}\nabla,\, \tau{\,\in\,}\Tau$ is a polarity if, and only if, 
$a^2{\,=\,}b\overline c{\,=\,}\varepsilon{\,=\,}\pm~1$  and 
$\overline{n}c{\,=\,}\varepsilon\overline{c}n$.
The unital $U_{\rho\delta\tau}$ is then described  by 
$c\,{\overline x}\,{\overline a}{\,\circ\,}x{\,=\,}\varepsilon y{\,+\,}c\,\overline y\,c{\,-\,}\varepsilon\hskip1pt n$. In particular, 
$U_\rho{\,=\,}\{(x,y)\mid \overline x{\,\circ\,}x{\,=\,}2y_0\}{\smcup}\{v\}{\,\approx\,}\Ss_7$. The motion group is $\Mu_\rho{\,=\,}\{(x,y){\,\mapsto\,}(axc,y^c{+}n)\mid a,c{\,\in\,}\HH', \overline n{\,=}-n\}$, and $\dim\Mu_\rho{\,=\,}9$. \par
If $\underline\delta$ corresponds to the choice  $\varepsilon{\,=}-\hskip-2pt1,\,a{\,=\,}c{\,=\,}k$,  we have
$$\rho\underline\delta{\,=\,}\kappa{\,:\,}(x,y){\,\mapsto\,}[\hskip1pt k\overline xk,k\overline yk\hskip1pt] 
{\rm \ and \ }
U_\kappa{\,=\,}\{(x,y)\mid k\overline xk{\,\circ\,}x{\,=\,}
y{\,+\,}\overline y^k{\,=\,}2(y{-}y_3k)\}{\smcup}\{v\}{\,\approx\,}\Ss_5\,,$$
since $y_0,y_1$, and $y_2$ are uniquely determined by $x{\,\in\,}\HH$ and $y_3{\,\in\,}\RR$. The centralizer 
$\Cs\Sigma\kappa$ consists of the maps $\delta\tau$ such that 
$k\overline cx\overline ak{\,=\,}bkxk\overline a$ and 
$k\overline cy\overline bk{\,+\,}k\overline nk{\,=\,}bkykc{\,+\,}n$ for all $x,y{\,\in\,}\HH$.  This implies first 
$k^a{\,=\,}ckb$ and $x^{ckb}{\,=\,}x^k$, then $cb^k{\,=\,}e{\,=\,}\pm1$, finally $ck{\,=\,}e k\overline b$,\,
$ka{\,=\,}eak$, and $\overline nk{\,=\,}\overline kn$. The last 3 conditions suffice to characterize $\Cs{}\kappa$, and $\dim\Mu_\kappa{\,=\,}5$.
\par
In the case $a^2{\,=\,}1,\,b{\,=\,}c$, each polarity $\rho\delta\tau$ is conjugate to $\rho$:  put 
$\gamma{\,:\,}(x,y){\,\mapsto\,}(xc,ayc)$,\,
$[\hskip1pt s,t\hskip1.5pt]{\,\mapsto\,}[\hskip1pt as,atc\hskip1.5pt] $, and note that  $\tau$ is the square of a  translation $\nu$ such that $\rho\delta\nu$ is a polarity. It follows that 
$\rho\delta{\,=\,}\rho^\gamma$ and 
$\rho\delta\tau{\,=\,}\rho^\gamma\nu^2{\,=\,}\nu^{-1}\rho^\gamma\nu{\,=\,}\rho^{\gamma\nu}$.  \par
For $\varepsilon{\,=}-\hskip-2pt1$, each polarity $\rho\delta\tau$ is conjugate to  $\kappa$. 
In fact, choose $h$ such that  $\overline a^h{\,=\,} k$, and let 
$\gamma{\,:\,}(x,y){\,\mapsto\,}(hx,cky),\;
[\hskip1pt s,t\hskip1pt]{\,\mapsto\,}[\hskip1pt cks\overline h,ckt\hskip1pt]$.
Then $\gamma\rho\delta {\,:\,}(x,y){\,\mapsto\,}
[\hskip1pt -c\overline x\overline h\overline a, -c\overline y\overline k\hskip1pt]=
[\hskip1pt -c\overline xk\overline h, -c\overline y\overline k\hskip1pt]{\,=\,}(x,y)^{\kappa\gamma}$   
and  $\rho\delta\tau{\,=\,}\kappa^\gamma\nu^2{\,=\,}\nu^{-1}\kappa^\gamma\nu$. \Qed  
\par\medskip
{\bf 10.7 Summary.} {\it All  self-dual $8$-dimensional planes with a group of dimension ${\ge\,}17$ ad\-mit polarities, the corresponding unitals of the non-classical planes have dimension $7$ or~$5$; possibly excepting the Hughes planes, all unitals are spheres\/}.
\par\medskip
\begin{center}
 \begin{tabular}{|l|c|c|c|c|c|c|} 
\hline
Planes & Classical & Hughes & Mutations & Rees & Spin & Double flag  \\  \hline
Classes  & 3 & $1$ & $2$ & $1$ & $2$ &  $2$  \\
Unitals & $\ \emptyset,\;\Ss_7,\;\Ss_5$ & $\dim U{\,=\,}5$ & $\ \;\Ss_7,\;\Ss_5$ & $\Ss_5$ & 
 $\Ss_5$ & $\ \;\Ss_7,\;\Ss_5$  \\ 
 $\dim\Mu$ & $21,\,21,\,15$ & 6 & $11,\,7$ & 7 & 7, \ 9 & $\ \ 9, \ 5$\\ 
\hline
\end{tabular}
\end{center}
\par\smallskip
{\tt Remark.} According to 12.3, {\it smooth\/} polar unitals of $8$-dimensional planes  are always 
$7$-spheres or $5$-spheres. In most cases it is not easy to determine the intersections of a polar 
unital $U$ with its secants, and it is not known if $U$ is a {\it topological\/} unital. 
\par\bigskip

{\Bf 11. Unitals in other dimensions} 
\par\medskip
{\bold 11.1 Moulton planes.} {\it Each proper Moulton plane has a unique conjugacy class of polarities.
In the standard description\/} (see 1.13) {\it one polarity is given by 
$\pi{\,:\,}(x,y){\,\leftrightarrow\,}[x,-y]$; the affine part $U_\pi\sm \{\infty\}$ of the corresponding unital 
consists of  two semi-parabolas, it satisfies the equation $x{\circ}x{\,=\,}2y$\/}, cf. \cite{sz10}.
\par\smallskip
{\tt Remark.} More generally, in each plane coordinatized by a commutative Cartesian field the map
$(x,y){\,\leftrightarrow\,}[x,-y]$ defines a polarity with unital $\{(x,y)\mid x{\circ}x{\,=\,}2y\}{\smcup}\{\infty\}$.
\par\smallskip
For polarities of other flat planes see \cite{Bd} and the r\'esum\'e in \cite{cp} 38.8.
\par\medskip
{\bold 11.2 Shift planes} exist only in dimension $2\ell{\,\le\,}4$, see \cite{cp} 74.6. 
Write $z{\,=\,}(x,y){\,\in\,}\KK^2$, $\KK{\,\in\,}\{\RR,\CC\}$. The lines of the  affine shift plane 
$\cP_{\hskip-1pt f}$ are the verticals and the {\it shifts\/} $L+z$  of the graph $L$ of a suitable 
function $f{\,:\,}\KK{\to\,}\KK$. We may assume that $0{\,\in\,}L$. 
According to \cite{cp} 74.10, the automorphism group $\Sigma{\,=\,}\Aut\cP_{\hskip-1pt f}$ 
is an extension $\Theta{\rtimes}\Gamma$ of the shift group $\Theta{\,\cong\,}\RR^{2\ell}$ by a 
linear group $\Gamma$. \par\smallskip
{\tt Theorem.} {\it In each shift plane, the map $\pi{\,:\,}z{\,\leftrightarrow\,}L{\,-\,}z$ defines 
a polarity; the corresponding unital is 
$U_\pi{\,=\,}\{(x,y)\mid 2y{\,=\,}f(2x)\}{\smcup}\{(\infty)\}{\,\approx\,}\Ss_\ell$, its affine part $U'_\pi$ 
is similar to $L$.  Obviously, $\pi\delta{\,=\,}\delta^{-1}\pi$  for each shift 
$\delta{\,:\,}z{\,\mapsto\,}z{+}d{\,\in\,}\Theta$, and   $\pi\delta{\,\sim\,}\pi$. An automorphism 
$\gamma\delta$ with  $\gamma{\,\in\,}\Gamma$  is a motion in $\Mu_\pi$ if, and only if, 
$L^\gamma{\,=\,}L{\,-\,}2d$; 
hence $\delta$ is determined by $\gamma$ and  $\Mu_\pi{\,\cong\,}\Gamma$\/}.
\par\smallskip
{\tt Remarks.} (1) The statements of the Theorem can easily be verified. There may exist other polarities: 
if $\gamma{\,\in\,}\Gamma$ and $L^\gamma{\,=\,}L{\,+\,}c$, then $\pi\gamma\delta$ is a polarity if, and only if, $\gamma^2{\,=\,}\1$ and $d^\gamma{\,=\,}c{+}d$. Let 
$\kappa{\,:\,}(x,y){\,\leftrightarrow\,}(-x,y)$, and assume that $L$ is $\kappa$-symmetric, 
$L^\kappa{\,=\,}L$. Then $U_{\pi\kappa}\sm\{\infty\}{\,=\,}\KK{\times}0$ and
the polarity $\pi\kappa$ is not conjugate to $\pi$, because  $(\KK{\times}0)^\sigma{\,\ne\,}U'_\pi$ for each $\sigma{\,\in\,}\Sigma$. (In fact, $(\KK{\times}0)^\sigma$ is a straight line and $U'_\pi$ 
is similar to $L$.) Motions  $\gamma\delta{\,\in\,}\Mu_{\pi\kappa}$ are subject to the conditions 
$L^\gamma{\,=\,}L$,  $\,d^\kappa{\,=\,}{-}d$, and  $\gamma\kappa{\,=\,}\kappa\gamma$. In general (e.g., for $\ell{\,=\,}\1$ and 
$f(x){\,= \,}{\rm cosh}(x){-}1$\,), the motion group $\Mu_{\pi\kappa}$ consists only of the shifts 
preserving $U_{\pi\kappa}$. \\
(2) A detailed study of polarities in general shift planes, not necessarily topological ones, 
can be found in  \cite{ks1}.
\par\medskip

{\bold 11.3  Hughes planes.} {\it The polarities of a proper $16$-dimensional Hughes plane $\cO_r$ 
form a single conjugacy class. Each polarity $\rho$ induces on the invariant {\rm inner} quaternion subplane $\cH$ a planar polarity. The group $\SU4(\CC,1)$ is a two-fold covering of the motion group $\Mu_\rho$, the unital $U_{\hskip-1pt\rho}$ is $11$-dimensional, and $\Mu_\rho$ is transitive on the set of outer points of  $U_{\hskip-1pt\rho}$\/}.
\par\smallskip
{\tt Proof.} (a) The group $\Gamma{\,\cong\,}\SL3\HH$ acts effectively on $\cO_r$ and is flag-transitive on the {\it outer\/} geometry of all points and lines not incident with the invariant subplane 
$\cH$, see   \cite{cp} \S~86, in particular 86.34.
Analogous to 9.17 or as in \cite{cp} 86.2, there is  an outer flag $(p,L)$ such that
the stabilizers $\Pi{\,=\,}\Gamma_{\hskip-1.5pt p}$  and 
$\Lambda{\,=\,}\Gamma_{\hskip-1pt L}$ have the form \\[0pt]
\centerline{$\Pi{\,=\,} \left \{ \left ( \begin{array}{cc}  C|h|&  \\ \frv&h^{-2} \end{array} \right )  
\;  \bigg | \;  C{\,\in\,}\U2\HH, \ h{\,\in\,}\{e^{(1{+}ri)s}\mid s{\,\in\,}\RR\}, \ \frv{\,\in\,}\HH^2 \right \}$ \ and \ $\,\Lambda{\,=\,}\Pi^t$.} \\[2pt]
According to 1.13, the outer part of $\cO_r$ can be written as 
$(\Gamma/\Pi{\,,\,}\Gamma/\Lambda)$. For each 
$\gamma{\,\in\,}\Gamma$ let $\overline\gamma$ denote the {\it transpose\/} of the conjugate of 
$\gamma$. Then
$\kappa{\,:\,}\gamma{\,\mapsto\,}\overline\gamma^{\hskip1pt k}$ is an  anti-automorphism of~$\Gamma$, and 
$\gamma{\,\mapsto\,}\gamma^\ast{\,:=\,}(\gamma^\kappa)^{-1}$ is an involutorial automorphism. One readily verifies that  $\Pi^\ast{\,=\,}\Lambda$. The map $\Pi\gamma{\,\leftrightarrow\,}\Lambda\gamma^\ast$ extends  to a polarity $\rho$ of $\cO_r$. As in 10.1(d), an outer point $\Pi\xi$ belongs to the unital 
$U_{\hskip-1.3pt\rho}$ if, and only if, there exists some $\gamma{\,\in\,}\Pi\xi{\smcap}\Lambda\xi^\ast$, and this implies $\gamma,\gamma^\ast{\,\in\,}\Pi\xi$, or $\gamma\gamma^\kappa{\,\in\,}\Pi{\smcap}\Lambda$. Conversely, it follows from  $\gamma\gamma^\kappa{\,\in\,}\Pi$ that 
$\Pi\gamma{\,\in\,}U_{\hskip-1.3pt\rho}$. \\
(b) For the motion group  we find 
$\Mu_\rho{\,=\,}\Cs\Gamma\rho{\,=\,}\{\gamma{\,\in\,}\Gamma\mid\gamma^\ast{\,=\,}\gamma\}$, since $\gamma\rho{\,=\,}\rho\gamma$ is equivalent with 
$\forall_{\xi{\in}\Gamma}\,\gamma^\ast{\in\,}\Pi^\xi\gamma$, and  
$\bigcap_{\,\xi{\in}\Gamma}\Pi^\xi{\,=\,}\1$. Thus $\Mu_\rho$ is the $15$-dimensional 
{\it anti-unitary\/} group, denoted by ${\rm U}_3^\alpha\HH$ in \cite{cp} 94.32/3, cf. also \cite{cp} 13.18. It induces on $\cH$ the planar  motion group $\PSU4(\CC,1)$ mentioned in 2.2 above, and 
$\rho|_\cH$ is a planar polarity. The stabilizer $\Mu_\rho{\smcap}\Pi$ of the absolute point $\Pi$ in 
$\Mu_\rho$  
is the set $\{\gamma{\,\in\,}\Pi\mid\gamma\gamma^\kappa{\,=\,}\1\}$. In the matrix representation,  
$\gamma\gamma^\kappa{\,=\,}\1$ implies that
$h^\kappa{\,=\,}\overline h^k{\,=\,}h{\,=\,}1,\,\frv{\,=\,}0$, 
and that $C^k{\,=\,}C$ has entries from $\RR(k){\,\cong\,}\CC$. 
Hence  $\Mu_\rho{\smcap}\Pi{\,\cong\,}\U2\CC$, and the $\Mu_\rho$-orbit of the point $\Pi$ is 
$11$-dimensional. Note that $\Mu_\rho{\smcap}\Pi{\,\le\,}\Lambda$. \\
(c) As $\cH$ is a Baer subplane of $\cO_r$, there is a unique {\it inner\/} point  $o{\,\in\,}L$ and 
$p$ is on a unique line $W{\,=\,}o^\rho$ of $\cH$. Choose $o$ as origin and $W$ as line at 
infinity for homogeneous coordinates in $\cH$ and note that the cyclic permutation matrix belongs 
to $\Mu_\rho$. It follows that the base points $(1,0,0),\,(0,1,0),\,(0,0,1)$ form a polar triangle 
$\frr$ for the polarity $\rho$. Let $\pi$ be any other polarity of $\cO_r$. Up to conjugation by 
$\Gamma$, we may assume that $\frr$ is also a polar triangle for $\pi$. Consequently 
$\rho\pi{\,=\,}\sigma$ fixes $\frr$ and $\sigma$ is given by a diagonal matrix 
${\rm diag}(s_1,s_2,s_3)$. From $\pi^2{\,=\,}\rho^2{\,=\,}\1$ we conclude that 
$\sigma^\rho{\,=\,}\sigma^*{\,=\,}\sigma^{-1}$ and
$\overline\sigma^k{\,=\,}\sigma$, hence $s_\nu{\perp}k$ with respect to the usual norm in 
$\HH$. Each $s_\nu$ is contained in a complex subfield of $\HH$ (or $s_\nu{\,\in\,}\CC$ up to an automrphism of $\HH$). Therefore  $s_\nu{\,=\,}w^2_{\hskip-.5pt\nu}$ with $w_\nu{\,\perp\,}k$; 
this means $\sigma{\,=\,}\omega^2$ with 
$\rho\hskip2pt\omega{\,=\,}\omega^{-1}\hskip-1pt\rho$, and  
$\pi{\,=\,}\rho\hskip2pt\omega^2{\,=\,}\omega^{-1}\hskip-1pt\rho\hskip2pt\omega$ is  
conjugate to $\rho$. \\
(d)  The group $\Mu_\rho$ is transitive on the  set of non-absolute lines  of the inner plane~$\cH$, see
\cite{cp} 18.32 for the dual statement. In order to show that $\Mu_\rho$ is transitive on the set 
$E$ of outer points of $U_{\hskip-1pt\rho}$, it suffices therefore to prove that the stabilizer 
$\Omega{\,=\,}(\Mu_\rho)_W$ of the (unique) inner line $W$ through $p$ is transitive on 
$E{\smcap}W{\,=\,}\{\Pi\gamma\mid\gamma{\,\in\,}\Gamma_{\hskip-1pt W}{\,\and\,}
\gamma\gamma^\kappa{\,\in\,}\Pi\}$. As $\dim\Omega{\,=\,}7$ and $\dim(\Omega\smcap\Pi){\,=\,}4$
by step (b), we have $\dim p^\Omega{\,=\,}3$. 
An element  $\gamma{\,\in\,}\Gamma_{\hskip-1pt W}$ is given by a matrix of the form 
 $\bigl({D \atop \frb}{ \atop d}\bigr)$, and $\gamma\gamma^\kappa{\,\in\,}\Pi$ 
 implies $\frb{\,=\,}0$ and $\gamma{\,\in\,}\Gamma_{\hskip-2pt o,W}$. Thus we may consider 
 $E{\smcap}W$ as a subset of the coset space $\Gamma_{\hskip-2pt o,W}/\Pi_o$. 
 We will show that  $E{\smcap}W{\,=\,}\{\Pi_o\delta\mid\delta{\,\in\,}\Delta\}$, where 
$$\Delta{\,=\,}\left \{\left(\begin{array}{cc}\hskip-3pt D & \\ &\hskip-6pt d\end{array}\right) 
\,\bigg | \,D{\,=\,}\Bigl({r \atop } \ {c \atop \raise2pt\hbox{$s$}}\Bigr), \ r{\,>\,}1{\,=\,}rs, \ c{\,\in\,}\RR i{+}\RR j,\  \   c\overline c{\,=\,}r^2{-}s^2, \ d\overline d{\,=\,}1, \ d{\,\in\,}\RR(k)\right\}\,.$$ 
In fact, let $\delta \delta^\kappa{\,\in\,}\Pi_o$ and note that $\Pi_o$ is $\kappa$-invariant. 
We may assume therefore that $|d|{\,=\,}1$ and that $D$ is a triangular matrix in $\SL2\HH$. 
From $dd^\kappa{\,=\,}1$ it follows that $d^k{\,=\,}d{\,\in\,}\RR(k)$. Because $|u|{\,=\,}|v|{\,=\,}1$ 
implies  $\bigl({u \atop }{ \atop v}\bigr){\,\in\,}\U2\HH$, the diagonal elements of $D$ can be 
chosen as positive reals, and then $D{\,\in\,}\U2\HH\Leftrightarrow D{\,=\,}\1$ and 
$\Delta{\smcap}\Pi_o{\,=\,}\1$.  Now  $\delta\delta^\kappa{\,\in\,}\Pi_o$ yields 
$$DD^\kappa{\,=\,}\Bigl(\;{r^2{+}cc^\kappa \atop  \raise2pt\hbox{$sc^\kappa$}} \ 
{sc \atop \raise2pt\hbox{$s^2$}}\;\Bigr){\,\in\,}\U2\HH\,.$$
From the bottom line we get $|c|^2{\,+\,}s^2{\,=\,}s^{-2}{\,=\,}r^2$, and orthogonality of the two lines implies $(r^2{+}c\overline c^k)c^k{\,+\,}s^2c{\,=\,}0{\,=\,}r^2c^k{\,+\,}(|c|^2{+}s^2)c$ or 
$c^k{\,+\,}c{\,=\,}0$, hence  $c{\,\in\,}\RR i{\,+\,}\RR j$. Consequently $\dim\Delta{\,=\,}3$ and 
$\Omega$ is indeed transitive on the connected $3$-manifold $E{\smcap}W$. \Qed
\par\medskip

{\bold 11.4 Mutations $\cO_{(t)}$.} The definition of the mutation $\cO_{(t)}$ of the octonion plane is analogous to 10.2. {\it The planes  $\cO_{(t)}$ admit two conjugacy classes of polarities, one with {\rm topological} unitals ${\approx\,}\Ss_{15}$ and a motion group of dimension $30$, the other with unitals   ${\approx\,}\Ss_{11}$ and a motion group of dimension $18$\/}.
\par\smallskip
{\tt Proof.} (a) Octonions will be written in the form 
$c{\,=\,}c'{+}c''l{\,\in\,}\OO$ with $c',c''{\hskip1pt\in\,}\HH, \  l{\,\perp\,}\HH$, and  $l^2{=}-\nobreak1$. Again $(x,y){\,\leftrightarrow\,}[x^\iota,-y^\iota]$ defines a polarity  $\rho_\iota$, 
whenever $\iota$ is an  involutorial anti-automor\-phism of~$\OO$. Consider 
$\overline\rho{\,=\,}\rho_\kappa, \, \kappa{\,:\,}z{\,\mapsto\,}\overline z$ and 
$\pi{\,=\,}\rho_\lambda$, where  $\lambda{\,=\,}{\kappa\underline\lambda}$ and
$\underline\lambda{\,:\,}z'{+}z''l{\,\mapsto\,}z'{-}z''l$ is an automorphism. Then 
$U_{\overline\rho}{\,=\,}\{(a,b)\mid
\overline aa{\,=\,} b{+}\overline b\}{\smcup}\{(\infty)\}{\,\approx\,}\Ss_{15}$,
since $a$ and $b{-\overline b}$ are arbitrary. Similarly
$$U_\pi{\,=\,}\{(a,b)\mid a^\lambda{\circ}a{\,=\,}b{+}b^\lambda{\,=\,}
2({\rm re}\,b'{\,+\,}b''l)\}{\smcup}\{(\infty)\}{\,\approx\,}\Ss_{11}\,,$$ 
because $b''$ and the real part of $b'$ are determined by $a$. Using the fact \cite{cp} 11.31d that all involutions of 
$\Gamma{=\,}\Aut\OO{\,\cong\,}\Gtwo$ are conjugate, it follows exactly as in 10.2(c) that each polarity of  $\cO$ is conjugate either to $\overline\rho$  or to  $\pi$.
\par
(b) The elements of the $40$-dimensional group $\Sigma$ have the form
$$\sigma:(a,b){\,\mapsto\,}(a^\gamma s{+}m, rb^\gamma s{\,+\,}q{\circ}a^\gamma s{\,+\,n)},\quad 
[c,d]{\,\mapsto\,}[rc^\gamma{+}q, rd^\gamma s{\,-\,}(rc^\gamma{+}q){\circ}m{\,+\,}n]$$
with $\gamma{\,\in\,}\Gamma, \, r,s{\,\in\,}\RR^{\times}, \, q,m,n{\,\in\,}\OO$, and  
$\Mu_{\overline\rho}{\,=\,}\{\sigma{\in}\Sigma\mid r{\,=\,}s,\,q{\,=\,}\overline m,\,
n{+}\overline n{\,=\,}|m|^2\}$. Therefore $\dim\Mu_{\overline\rho}{\,=\,}30$. Similarly, 
$\Mu_\pi{\,=\,}\{\sigma{\in}\Sigma\mid r{\,=\,}s,\,\HH^\gamma{\,=\,}\HH,\,
q{\,=\,}\overline m^{\underline\lambda},\,
n{+}\overline n^{\underline\lambda}{\,=\,}\overline m^{\underline\lambda}{\,\circ\,}m\}$. We have
$\overline n^{\underline\lambda}{\,=\,}\overline{n'}{+}n''l, \ 
n{+}\overline n^{\underline\lambda}{\,=\,}2(n_0{+}n''l)$ and 
$\Gamma_{\hskip-1.5pt\HH}{\,\cong\,}\SO4\RR$ (see \cite{cp} 11.31b). Hence 
$\dim\Mu_\pi{\,=\,}18$. \par
(c) {\it Each secant intersects $U_{\overline\rho}$ in a $7$-sphere\/}:  it suffices~to  
consider horizontal secants and verticals, since he subgroup 
$\{(a,b){\,\mapsto\,}(a{+}m,b{+}\overline m{\circ}a{+}\frac{1}{2}|m|^2)\}$ of  $\Mu_{\overline\rho}$ 
is transitive on the set of  slopes $\ne\infty$.  The assertion follows now immediately from the equation $\overline aa{\,=\,} b{+}\overline b$ for the unital. \Qed
\par\medskip
{\bf 11.5 Fixed double flag.} A $16$-dimensional plane $\cP$ with an automorphism group $\Sigma$ of dimension at least $33$ such that $\Sigma$ fixes $2$ points and $2$ lines can be coordinatized by a Cartesian field $(\OO,+, \circ)$ of distorted octonions defined in the same way as the distorted quaternions  in 6.4, and $\Sigma$ has a compact subgroup $\Phi{\,\cong\,}\Spin8\RR$ and a transitive group $\Tau$ of vertical translations, see \cite{hs2}. Assume that $(\RR,+,\circ)$ is associative (and hence also commutative; equivalently, $\cP$ can be coordinatized by a Hurwitz ternary field obtained from $\OO$ by radial distorsion of the addition, cf. \cite{ps}). 
In this case, the multiplication $\circ$ can be described as follows: let 
$\mu{\,:\,}[\hskip1pt0,\infty){\,\approx\,}[\hskip1pt0,\infty)$ with $\mu(1){\,=\,}1$, put
$r{\,\ast\,}s{\,=\,}\mu^{-1}(\mu(r){\cdot}\mu(s))$ for $r,s{\,>\,}0$ and 
$a{\,\circ\,}x{\,=\,}|a|{\,\ast\,}|x|\,|ax|^{-1}\,ax$ for $a,x{\,\ne\,}0$. Then 
$\dim\Sigma{\,=\,}37$, the map $z{\,\mapsto\,}\overline z$ is an involutorial anti-automorphism, 
$(x,y){\,\leftrightarrow\,}[\hskip1.5pt\overline x,-\overline y\hskip1.5pt]$ defines a polarity  $\pi$, and 
$$U_\pi{\,=\,}\{(x,y)\mid \overline x\circ x{\,=\,}2y_0\}{\smcup}\{\infty\}{\,\approx\,}\Ss_{15}$$ 
analogous to 10.6.  The group $\Xi{\,=\,}\{(x,y)\mapsto(x,y{+}t)\mid \overline t{\,=}-t\}$ of
translations by a pure element is contained in the motion group $\Mu_\pi$, and 
$o^\Xi{\,\subseteq\,} U_\pi$; 
for each translation $\tau$ by a real element $t$, however, 
 $\pi\tau^2{\,=\,}\pi^\tau$ and $\pi\tau{\,\sim\,}\pi$. Associativity of $(\RR,\circ)$ 
 implies that $\nabla{\,=\,}\Sigma_0{\,=\,}\Rho{\times}\Phi$, where  $\Rho{\,\cong\,}e^\RR$ consists of the maps $(x,y){\,\mapsto\,}(r{\circ}x,y)$ with $r{\,>\,}0$. By stiffness (\cite{cp} 83.23) or by 
\cite{cp} 96.36, the group $\nabla$ is transitive on $ou\sm\{o,u\}$; hence $\Sigma$ is transitive on the set of all points not incident with one of the fixed lines.
The center of $\Phi$ is generated by reflections $\alpha$ with axis $ov$ and $\sigma$ with axis 
$uv$. Easy verification shows that $\alpha{\,\in\,}\Mu_\pi$ and that $\pi\alpha{\,=\,}\pi^\sigma$ is conjugate to $\pi$, whereas 
$\pi\sigma$ is not a polarity.  From  $\gamma{\,\in\,}\Rho$ it follows that 
$\pi\gamma^2{\,=\,}\pi^\gamma$, and then  $\pi\gamma{\,\sim\,}\pi$. \\
The automorphism group $\Aut\OO{\,=\,}\Gtwo$ can be identified with $\Lambda{\,=\,}\Phi_{(1,1)}$ via  $(x,y)^\lambda{=}(x^\lambda,y^\lambda)$, and $\Lambda{\,\le\,}\Mu_\pi$ because 
$\overline x^\lambda{\,=\,}\overline{x^\lambda}$, see \cite{cp} 11.28.  All involutions of $\Lambda$ are conjugate in $\Lambda$ to  $\underline\lambda{\,:\,}c'{+}c''l{\,\mapsto\,}c'{-}c''l$ 
(see \cite{cp} 11.31d), and  
$\kappa{\,=\,}\pi\underline\lambda$ is  polarity with unital \\[3pt]
\centerline{$U_{\hskip-1pt\kappa}{\,=\,}\{(x,y)\mid \overline x^{\hskip1pt\underline\lambda}{\,\circ\,}x
{\,=\,}2(y_0{+}y''l)\}{\smcup}\{\infty\}{\,\approx\,}\Ss_{11}$,} \\[3pt]
hence $\pi\lambda{\,\sim\,}\kappa$ for each involution $\lambda{\,\in\,}\Lambda$, and none of these polarities is conjugate to $\pi$. If $\lambda{\,\in\,}\Lambda$ and $\pi\lambda$ is a polarity, then
$\lambda$ is necessarily an involution.\\
Restriction of the coordinates to the distorted quaternions yields a Baer subplane $\cB{\;\ldot}\cP$.
As $\underline\lambda|_\cB{\,=\,}\1$, both $\pi$ and $\kappa$ induce the same polarity $\rho$ on 
$\cB$, the polarity with $U_{\hskip-1.5pt\rho}{\,\approx\,}\Ss_7$ described in 10.6.
\par\smallskip
{\tt Theorem.} {\it Under the assumptions made at the beginning of\/} 11.5, {\it each polarity $\rho$ 
of $\cP$ is conjugate either to $\pi$ with $U_\pi{\,\approx\,}\Ss_{15}$ or to 
$\kappa$ with $U_\kappa{\,\approx\,}\Ss_{11}$. Moreover, $\Mu_\pi{\,=\,}\Xi{\rtimes}\Phi_{(0,1)}$ so that  $\dim\Mu_\pi{\,=\,}28$, and $\Mu_\pi$ has infinitely many orbits on $U_\pi$\/}.
\par\smallskip
{\tt Proof.} (a)  {\it Up to conjugation, each polarity $\rho$ of $\cP$ satisfies $o^\rho{\,=\,}ou$\/}. 
In fact, $\rho$ induces on $o^\Tau{\,\cong\,}\RR^8$ an involutorial map 
$\omega{\,:\,}y{\,\mapsto\,}yS{\,+\,}t$ with $S{\,\in\,}\O8\RR$, and $yS^2{\,+\,}tS{\,+\,}t{\,=\,}y$ for all $y{\,\in\,}\RR^8$. Consequently, $tS{\,=\,}-t$ and $S^2{\,=\,}\1$. In suitable coordinates,  
$S{\,=\,}\bigl({-\1 \atop 0}\;{0 \atop \1}\bigr)$ and  $t$ has the form $(m,0)$; obviously  
$\omega$ has a fixed point. \\
(b) We may assume that $(1,0)^\rho$ is a line $[\hskip1pt s,0\hskip1pt]$ of slope $s{\,>\,}0$, since 
$\Phi_{(1,0)}{\,\cong\,}\Spin7\RR$ acts on the line $uv$ with orbits homeomorphic to $\Ss_7$, 
see \cite{cp} 96.36. Now we get  $(1,s)^\rho{\,=\,}[\hskip1pt s,-s\hskip1pt]$ and  
$(0,s)^\rho{\,=\,}[\hskip1pt 0,-s\hskip1pt]$. As $\rho$ induces a linear map on $\Tau{\,\cong\,}\RR^8$ 
via $(o^\tau)^\rho{\smcap}ov{\,=\,}o^{\tau^\rho}$, it follows that 
$(0,1)^\rho{\,=\,}[\hskip1pt 0,-1\hskip1pt]$. The collineation $\pi\rho{\,=\,}\sigma$ fixes the triangle $o,u,v$ and the point $(0,1)$, and $\sigma$ maps $(s,0)$ to $(1,0)$. Let $r{\,\ast\,}r{\,=\,}s$ and define the homology $\eta{\,\in\,}\Rho$ by $(x,y)^\eta{\,=\,}(r{\circ}x,y)$. Then 
$(1,0)^{\eta^2}{\,=\,}(s,0)$ and $\eta^2\sigma{\,=\,}\lambda{\,\in\,}\Lambda$. 
Put $\gamma{\,=\,}\eta^{-1}$ and note that $\Rho{\,\le\,}\Cs{}\Lambda$. 
Again $\lambda^2{\,=\,}\1\,$ (in fact,  $\eta\pi{\,=\,}\pi\eta^{-1}$,   
$\,\overline x^\lambda{\,=\,}\overline{x^\lambda}$ implies $\pi\lambda\pi{\,=\,}\lambda$,  
$\,\pi\lambda{\,=\,}\rho\eta^2{\,=\,}\eta^{-2}\rho$, and 
$\lambda^2{\,=\,}(\pi\lambda)^2{\,=\,}\rho^2{\,=\,}\1$).
Either  $\rho{\,=\,}\pi{\gamma^2}{\,=\,}\pi^\gamma$ \,(if $\lambda{\,=\,}\1$) or 
$\rho{\,=\,}\pi^\gamma\lambda{\,=\,}(\pi\lambda)^\gamma{\,\sim\,}\kappa^\gamma$. \\
(c) In order to show that $\Phi{\smcap}\Mu_\pi{\,=\,}\Phi_{(0,1)}$, we use the description of $\Phi$ given 
in \cite{hs2} p.~97 (see also \cite{cp} 12.17). Each $\phi{\,\in\,}\Phi{\,\cong\,}\Spin8\RR$ has the form 
$(x,y){\,\mapsto\,}(x\alpha,y\beta), \ (s){\,\mapsto\,}(s\gamma)$ with 
$\alpha,\beta,\gamma{\,\in\,}\SO8\RR$ and identically 
$(sx)\beta{\,=\,}s\gamma{\hskip1pt\cdot\hskip1pt}x\alpha$. 
Here $\alpha$ and $\gamma$ preserve the norm but are not necessarily linear. By definition,  
$\phi{\,\in\,}\Mu_\pi\Leftrightarrow\overline{x\alpha}=\overline x\gamma\,\and\,
\overline{y\beta}=\overline y\beta$. The second condition implies 
$1\beta{\,=\,}\varepsilon{\,=\,}\pm1$, and $s{\,=\,}\overline x$ yields 
$(\overline xx)\beta{\,=\,}\overline x\gamma{\hskip1pt\cdot\hskip1pt}x\alpha$, 
$\,\varepsilon|x|^2\hskip2pt\overline{x\alpha}{\,=\,}\overline x\gamma|x\alpha|^2$. Now  
$\varepsilon{\,=\,}1$ (because $\Mu_\pi$ 
is connected) and   $\overline{x\alpha}{\,=\,}\overline x\gamma$. 
If $(a,b){\,\in\,}U_\pi$ and $(x,y){\,\in\,}(a,b)^{\Mu_\pi}$, then $|x|{\,=\,}|a|$.
\par\medskip

{\bf 11.6 Summary: Unitals of 16-dimensional planes.}
\begin{center}
 \begin{tabular}{|l|c|c|c|c|} 
\hline
Planes & Classical & Hughes & Mutations & Double flag   \\  \hline
Classes  & $3$ & $1$ & $2$ &  $2$     \\
Unitals & $\quad \emptyset,\;\Ss_{15},\;\Ss_{11}$ & $\dim U{\,=\,}11$ & $\;\Ss_{15},\;\Ss_{11}$ &  
$\Ss_{15},\;\Ss_{11}$ \\ 
$\dim\Mu$ & $52,\,52,\,36$ & $15$ & $30,\,18$ & $28$, \hskip6pt ?  \\ 
\hline
\end{tabular}
\end{center}
\par\bigskip

{\Bf 12. Topological and smooth unitals}
\par\medskip
{\it Polar\/} unitals $U$ of the classical compact connected planes have the following properties: \par
\quad(1) for each $x{\,\in\,}U$ there is a unique line $L$, the {\it tangent\/}, such that 
$L{\smcap}U{\,=\,}\{x\}$, and \par
\quad(2) if $L$ is a line with $|L{\smcap}U|{\,>\,}1$ (a {\it secant\/}), then 
$L{\smcap}U{\,\approx\,}\Ss_k$ for a fixed number $k{\,\ge\,}0$. \\
Any closed subset $U$ of an arbitrary compact projective plane of dimension $2\ell$ satisfying (1) and (2) will be called a {\it topological\/} unital, cf. \cite{Iv2} 2.4. In general, the two notions of a unital must be distinguished.  Property (1) holds for any \emph{polar} unital of an arbitrary compact connected plane. The set of all tangents at $U$ is denoted by $U^*$. \\
A topological unital in a smooth projective plane is said to be a {\it smooth\/} unital if 
$U$ is a smooth submanifold homeomorphic to a sphere (of dimension $k{+}\ell$) and if each secant 
intersects $U$ transversally. Note that all classical polar unitals are smooth.
\par\medskip
{\bf 12.1 Proposition.} {\it Let $U$ be a topological unital in an $8$-dimensional compact plane. If 
$U$ is homeomorphic to a sphere and if there is some point $p{\,\notin\,}U$ such that each line $L$ in the pencil $\frL_p$ intersects $U$ and  $\frL_p{\smcap}U^*{\,\approx\,}\Ss_k$, then 
$U{\,\approx\,}\Ss_5$\/}, see \cite{Iv2} 3.2.
\par\medskip
{\bf 12.2 Theorem.} {\it Let $U$ be a smooth unital in the smooth plane $(P,\frL)$. Then $U^*$ is compact and the set of secants is open in $\frL$. If $\dim U{\,<\,}\dim P{\,-\,}1$, then each line intersects~$U$ and each pencil consists of secants and tangents\/}, see \cite{Iv2} 4.\hskip1pt3,4.
\par\medskip
{\bf 12.3 Theorem.} {\it A smooth \emph{polar}  unital $U$ in a $2\ell$-dimensional smooth plane is homeomorphic to a sphere  $\Ss_{\hskip1pt2\ell - 1}$ or $\,\Ss_{\hskip1.5pt(3\ell/2)-1}$; 
each secant intersects $U$ in a smooth submanifold homeomorphic to $\Ss_{\hskip1pt\ell-1}$ or 
$\,\Ss_{\hskip1pt(\ell/2)-1}$\/}, see \cite{Iv3} and \cite{Iv4}. 
\par\bigskip
\break
{\Bf 13. Appendix: Variations on a theme of Tschetweruchin} \par\medskip
In 1927 Tschetweruchin \cite{ts}  gave an example of a $2$-dimensional compact projective plane such that Desargues' theorem does not hold in any open region of the point space. His construction can easily be generalized. We discuss the properties of the corresponding planes and their duals.
\par\smallskip
In affine form, Tschetweruchin's  plane $\cT$  has the point space $\RR^2$; lines are the verticals 
$r{\times}\RR$, the ordinary straight lines $\{(x,sx{+}t) \mid x{\,\in\,}\RR\}$ with slope $s{\,\ge\,}0$, and the cubic curves with an equation $y^3{\,=\,}s^3x^3{+}t^3$ for $s{\,<\,}0$.
The dual $\cS$ of $\cT$ is obtained by inter\-changing the r\^ oles of $s$ and $x$. Its right half-plane 
$\cS^+{\,=\,}\{(x,y) \mid x{\,>\,}0\}$ is isomorphic  to a real half-plane $\cD$. The lines in the left 
half-plane $\cS^-$ satisfy linear equations in the coordinates $(\xi,\eta){\,=\,}(x^3,y^3)$. 
Hence $\cS^-$ is also isomorphic to $\cD$, and $\cS$ is composed of two copies of $\cD$ glued together at a line. In particular, the plane $\cT$ is not self-dual.
\par\smallskip
{\bf The general case.} Let $\rho$ be a homeomorphism of $\RR$ which fixes 
$0$ and $1$ (e.g. \ $x^\rho{\,=\,}x$ for $x{\,\ge\,}0$  and $x^\rho{\,=\,}kx$ for $x{\,<\,}0$ with some $k{\,>\,}1$).
Consider the ternary map $\tau:\RR^3\to\RR$ defined by \ 
$\tau(s,x,t){\,=\,}s\hskip1pt x{\,+\,}t \ (s{\,\ge\,}0)$ \  and \  
$\tau(s,x,t)^\rho{\,=\,}s^\rho x^\rho{\,+\,}t^\rho \ (s{\,<\,}0)$. \par
Lines of the geometry $\cT_\rho$ are the sets
$[s,t]{\,=\,}\{(x,y){\,\in\,}\RR^2 \mid y{\,=\,}\tau(s,x,t)\}$ with $s,t{\,\in\,}\RR$ and the verticals 
$[r]{\,=\,}r{\times \RR}$.
Add a common point $(\infty)$ to all verticals and a common point $(s)$ to all lines $[s,t]$ of 
slope $s$; these points form the line $[\infty]$ at infinity, cf. 1.4.   Then \par
$\cT_\rho$ {\it is a compact topological projective plane\/}. In fact, the lines in $\RR^2$ form an affine plane because lines of different slope intersect in a unique point and distinct lines of equal
slope are disjoint. By \cite{cp} (32.2), $\cT_\rho$ is a topological plane. 
\par\smallskip
Assume from now on  that identically $(-x)^\rho{\,=\,}-(x^\rho)$. 
In order to describe the dual $\cS_\rho$ of the plane $\cT_\rho$, we interprete the point 
$(x,y)$ of  $\cT_\rho$ as the line $[x,y]$ in $\cS_\rho$ and vice versa, 
and we represent the lines in $\cS_\rho$ by a ternary map $\tilde\tau$. 
Then $y{\,=\,}\tau(s,x,t)$  or  $-t{\,=\,}\tau(x,s,-y)$ becomes $t{\,=\,}\tilde\tau(x,s,y)$. 
Changing signs of the second coordinate and interchanging the names  $(x,y)$ and $(s,t)$
yields \vskip-15pt
$$\tilde\tau(s,x,t){\,=\,}\tau(x,s,t){\,=\,}
\begin{cases}\hskip14pt sx{\,+\,}t\hskip36pt(x{\,\ge\,}0)\cr 
                     (s^\rho x^\rho{\,+\,}t^\rho)^{\rho^{-1}}\quad(x{\,<\,}0)\cr\end{cases}\hskip-9pt .$$

As in the case $\rho:x{\,\mapsto\,}x^3$, each plane $\cS_\rho$ is composed of two real half-planes bounded by the lines $[0]$ and $[\infty]$. 
\par\smallskip
The ternary function $\tau$  determines multiplication and addition by 
$s$\smotimes$x{\,=\,}\tau(s,x,0)$ and $x$\smoplus$t{\,=\,}\tau(1,x,t)$, see  1.4 or   \cite{cp} (22.4).
The operations  \tsmotimes\ and\tsmoplus\ are related to $\tilde\tau$ in the same way.  Obviously, $x$\smoplus$t{\,=\,}x{\,+\,}t$ is the usual sum in $\RR$.
\par\smallskip
{\bf Multiplicative homeomorphisms.}  If $\rho$ is multiplicative, i.e. if    
$$x^\rho{\,=\,}\begin{cases}\hskip13pt x^r \hskip12pt (x{\,\ge\,}0)\cr -|x|^r \hskip9pt  (x{\,<\,}0)\cr\end{cases}$$ 
for  some fixed real number $r{\,>\,}0$, then both \smotimes\ and                  
\tsmotimes\ coincide with the ordinary multiplication of $\RR$. For $x{\,<\,}0$ however,                    
$x$\tsmoplus$t{\,=\,}x{\,+\,}t$ implies $(x{\,+\,}t)^\rho{\,=\,}x^\rho{\,+\,}t^\rho$, and this is true  for all $x$ only if 
$t{\,=\,}0$. These facts can be used to determine the full automorphism group
$\Gamma$ of $\cS_r$ and $\cT_r$.
\par\smallskip
{\bf The planes  $\cS_r$ with $r{\,\ne\,}1$.} Each map $(x,y)\mapsto(ax,y)$ with $a{\,>\,}0$ is a homology with axis $[0]$ and center $(0)$. The maps  $(x,y)\mapsto(x,by)$ with  $b{\,\in\,}\RR^{\times}$ form a transitive group of  $(\infty)$-$[0,0]$-homologies; in particular, 
$\beta:(x,y)\leftrightarrow(x,-y)$ is a reflection. However, there does not exist
a $(0)$-$[0]$-reflection. In fact, such a reflection has necessarily the form 
$(x,y)\mapsto(ex,y)$ with $e^2{\,=\,}1$ or $e{\,=\,}-1$, but this would imply
$(x{\,+\,}1)^r {\,=\,} x^r{\,+\,}1$ for all $x{\,>\,}0$. Similarly, the remark above on \tsmoplus\ shows that  the group of $(\infty)$-$[\infty]$-translations is trivial, and so is the group of 
$(\infty)$-$[0]$-elations. There exists a reflection $\alpha$ with axis $[1]$ and center $(-1,0)$;  in homogeneous coordinates  $\alpha$ has the form $(x,y,z)\leftrightarrow(z,y,x)$. The map
$\alpha\beta$ is a $(1,0)$-$[-1]$-reflection. We will show that $\Gamma$ is an extension of its
connected component $\Gamma^1{\,=\,}\{(x,y)\mapsto(ax,by) \mid a,b{\,>\,}0\}$ by $\alpha$ and $\beta$.
\par\smallskip
A $2$-dimensional locally compact plane is called 
locally classical, if each point has a neighbourhood which is isomorphic to a copy of the real hyperbolic plane $\cH$. By results of Polley \cite{po} and Br\"ocker \cite{br}, a locally classical affine plane is isomorphic to the real plane, cf. also Polley \cite{pl} and \cite{gl} 5.18. It follows that $\cS_r$ is not locally classical  at the lines $[0]$ and $[\infty]$. Hence each collineation fixes these two lines or interchanges them,  their intersection $(\infty)$ is a fixed point of $\Gamma$. The axis of the
$(\infty)$-$[0,0]$-reflection $\beta$ is also $\Gamma$-invariant, or the product of two conjugates of 
$\beta$ would be an elation.
\par\smallskip
Suppose that some collineation $\gamma$ interchanges the half-planes $\cS_r^\pm$. We may assume that $(1,1)^\gamma{\,=\,}(-1,1)$ and that $\gamma$ fixes $(0)$. Then we get successively
$$(0,1)^\gamma{=\,}(0,1),\,  (1)^\gamma{=\,}(-1), \, (1,0)^\gamma{=\,}(-1,0), \,
 [-1,1]^\gamma{=\,}[1,1], \, (-1)^\gamma{=\,}(1), \, (-1,1)^\gamma{=\,}(1,1).$$
It follows that $\gamma^2$ fixes a quadrangle, and then $\gamma^2{\,=\,}\1$ and $\gamma$ would be a reflection with center $(0)$ by \cite{cp} (32.\hskip1pt10,12), but we have seen that such a reflection does not exist.
\par\smallskip
{\bf Proposition.} {\it The planes $\cS_r$ and $\cS_{r'}$ are isomorphic if, and ony if, 
$r'{\,=\,}r^{\pm1}$\/}.
\par\smallskip
{\tt Proof.} 
The lines in the left half-plane $\cS_r^-$ are given by linear equations in the coordinates 
$(\xi,\eta){\,=\,}(x^\rho,y^\rho)$, and  $\cS_r^-{\,\cong\,}\cD$. If the r\^oles of the two half-planes are interchanged, it follows that $\cS_r$ is isomorphic
to $\cS_{r^{-1}}$. Hence we may assume  that $r,r'{\,>\,}1$, and that an isomorphism 
$\omega:\cS_r\to\cS_{r'}$ is the identity on the point set of  $\cS_r^+$. Now
$$[1,0]^\omega{\,=\,}[1,0],\ [0,t]^\omega{\,=\,}[0,t],\ (x,x)^\omega{\,=\,}(x,x),\ [x]^\omega{\,=\,}[x],\ 
(x,y)^\omega{\,=\,}(x,y),$$
and $\omega$ maps each point onto itself. The fact that $\omega$ maps lines of $\cS_r^-$ onto lines of $\cS_{r'}^-$ yields $(|x|^r{\,+\,}|t|^r)^{r'}{\,=\,}(|x|^{r'}{\,+\,}|t|^{r'})^r$. With
$x{\,=\,}t{\,=\,}-1$ we obtain $2^{r'}{\,=\,}2^r$ and $r'{\,=\,}r$. \Qed

\par\smallskip
{\bf The planes  $\cT_r$.} All assertions on the planes  $\cS_r$ are true  in dual form in the planes $\cT_r$. 
In particular, the group of $(0,0)$-$[\infty]$-homologies is transitive, and there exists a reflection with center 
$(-1)$ and axis $[1,0]$.  Moreover,
{\it the planes $\cT_r$ are nowhere classical\/}.
In fact,   assume that some point of $\cT_r$ has a classical neighbourhood isomorphic to $\cH$.
In $\cH$ there exists a Pappos  configuration $\cK$ such that the slopes of all $9$ lines of $\cK$ differ by at most 
$\pi/3$, and  $\cK$ can be rotated in such a way that exactly one line $L$ of $\cK$ has negative slope.   
We may suppose that $L$ does not pass through the origin.
It follows that $L$ contains $3$ points of a real line  and even a real interval (fix $2$ points of $\cK$ on $L$ and $3$ collinear points of $\cK$ outside $L$, and vary a further point in a small neighbourhood).
This implies  $r{\,=\,}1$.\Qed
\par

\end{document}